\documentclass{amsproc}
\usepackage{geometry}
\usepackage{amsmath, amsthm, amssymb}
\usepackage[utf8]{inputenc}
\usepackage[T1]{fontenc}
\usepackage[all]{xy}
\usepackage[colorlinks=true, pdfstartview=FitV, linkcolor=blue,citecolor=blue, urlcolor=blue]{hyperref}
\usepackage[abbrev,lite,nobysame]{amsrefs}
\usepackage{times}
\usepackage[usenames,dvipsnames]{color}
\usepackage{graphicx}
\usepackage{soul}  

\usepackage{dsfont}

\def\R {\mathbb{R}}

\newtheorem{proposition}{Proposition}[section]
\newtheorem{theorem}[proposition]{Theorem}
\newtheorem{corollary}[proposition]{Corollary}
\newtheorem{lemma}[proposition]{Lemma}
\theoremstyle{definition}
\newtheorem{definition}[proposition]{Definition}
\newtheorem{remark}[proposition]{Remark}

\numberwithin{equation}{section}

\title{Rigorous FEM for 1D Burgers equation}
\author{Piotr Kalita}
\address{Faculty of Mathematics and Computer Science, Jagiellonian University, ul. \L{}ojasiewicza 6, 30-348 Krak\'{o}w, Poland}
\email{piotr.kalita@ii.uj.edu.pl}
\author{Piotr Zgliczy\'{n}ski}
\address{Faculty of Mathematics and Computer Science, Jagiellonian University, ul. \L{}ojasiewicza 6, 30-348 Krak\'{o}w, Poland}
\email{umzglicz@cyf-kr.edu.pl}
\begin{document}
	\begin{abstract}
		We propose a method to integrate dissipative PDEs rigorously forward in time with the use of Finite Element Method (FEM). The technique is based on the Galerkin projection on the FEM space and estimates on the residual terms. The proposed approach is illustrated on a periodically forced one-dimensional Burgers equation with Dirichlet conditions.  For two particular  choices of the forcing we prove
the existence of  the periodic globally attracting trajectory and give precise bounds on its shape.
	\end{abstract}
\maketitle

\tableofcontents

\section{Introduction}
\label{sec:introduction}
In the study of evolutionary partial differential equations (PDEs) purely analytical methods often appear insufficient to gain full understanding
of the behaviour of solutions. While, typically, they allow us to obtain well-posedness, and, in dissipative situations, existence and some basic properties of the global attractor, the precise description of the dynamics on the attractor is often out of reach with such tools \cites{Robinson, Temam}.  On the other hand with the present day computers a lot of interesting PDEs can be numerically investigated and quite often on the heuristic
level  the dynamics of such systems can be well understood. This however does not constitute a mathematical proof and there are situations
when the numerical simulations can be misleading.
 Our goal in this paper is to take a FEM discretization method and to build rigorous numerics
around it.
By rigorous integration of dissipative PDE we understand an algorithm which produces a set which is guaranteed to contain \textit{all} solutions of a given problem originating from \textit{any} initial data in a given set. We want these bounds to be of sufficient quality to satisfy assumptions of some abstract theorem which gives us the existence of some interesting dynamical
object - in the present paper it is a periodic orbit.

The task of designing algorithms  for rigourous numerics for evolution in time of parabolic PDEs is not new. For problems with periodic boundary
conditions it was realized (using the Fourier basis) with considerable success,
\begin{itemize}
\item for the Kuramoto--Shivashinsky equation on the line where the periodic \cites{mizgli, enclosiure_ref, kuramoto_III,AK} and, more recently chaotic \cite{zgliwi} solutions were verified,
\item  for the Burgers equation where periodic orbit
attracting all solutions was shown \cite{cyranka},
\item for  the one-dimensional Ohta--Kawasaki model in \cite{cyranWanner}  some heteroclinic orbits between some fixed points  has been proven to exist.
\end{itemize}

The main point in our paper is the use of FEM basis for the algorithm.  This choice is motivated by the wide applicability of FEM in numerical
solving  of PDEs for general boundary condition. Before we outline our FEM-based approach let us highlight first some properties of the Fourier basis, which were crucial in works \cites{mizgli, enclosiure_ref, kuramoto_III,AK,zgliwi,cyranka,cyranWanner}
\begin{itemize}
\item linear differential operators with constant coefficients considered there are diagonal in the Fourier basis,
\item product of two functions from the Fourier basis is a function which is an element of the same Fourier basis.
\end{itemize}
With these properties the nonlinear terms  can be easily expressed and estimated via convolutions. The diagonality of the leading differential operator greatly helps with getting a priori bounds for short time intervals.
These properties make it possible to design an algorithm for rigorous integration of dissipative PDEs with periodic boundary conditions,
using the method of self-consistent bounds developed in \cites{mizgli,kuramoto_III}. In this approach  it is possible to obtain rigorous bounds on trajectories for all sufficiently high-dimensional Galerkin
projections and then after passing to the limit one gets rigorous bounds for the solutions of the PDE under consideration.

In the present paper  as an object to study we chose  one dimensional problem governed by the following Burgers equation with Dirichlet boundary conditions
\begin{align}
& u_t - u_{xx} + u u_x = f(x,t)\quad \textrm{for}\quad (x,t)\in (0,1) \times (t_0,\infty), \label{a}\\
& u(0,t) = u(1,t) = 0.\label{b}
\end{align}
We use the simplest FEM basis of first order Lagrange elements (intervals).  We treat this example as a toy model, which contains
probably all difficulties coming from the choice of a quite arbitrary basis, hence the experience gained studying it should be  transferable to other equations
and boundary conditions. Observe that in the FEM basis none of the above listed properties of the Fourier basis are satisfied, moreover
some additional problems arise due the following issue
\begin{itemize}
\item functions in the FEM space are typically not smooth, while solutions of the PDE, depending on the regularity of forcing term $f$, are smooth as functions of the space variable.
\end{itemize}
We demonstrate the correctness of our approach by showing that for problem \eqref{a}--\eqref{b} for two particular choices of periodic in time $f$ there exists the periodic solution $u$ in some neighborhood of numerical solution found by the standard, nonrigorous, FEM.
Let us stress that problem \eqref{a}--\eqref{b} with periodic forcing very likely can be treated using the Fourier basis (i.e. the sine Fourier series) and the technique developed
in \cite{cyranka}, but the result on attracting periodic solution of \eqref{a}--\eqref{b} is not the main goal of the paper. The goal is to outline a FEM-based method for rigorous integration of dissipative PDEs.


One time step of length $h>0$ of our method consists of the following two stages.
\begin{itemize}
\item[STAGE 1] Given a set $N$ of initial conditions taken at time $t_0$, which is  bounded, closed, and convex in $H^s$,  we construct the set $W$, such that $N \subset W$, $W$ is the bounded, closed, and convex set in $H^s$, such that all solutions starting from $N$ are defined for $t \in [t_0,t_0+h]$ and are contained in $W$.

\noindent The value of $s$ is dictated by the need to estimate the projection error of the Laplacian operator in the FEM-basis. To this end, we need $s=4$ for Burgers equation and Lagrange elements, while for other differential operators and more regular FEM bases $s$ might be bigger.

 \noindent   Effective realization of this part of the algorithm depends on the particular equation. In the case of the Burgers equation global
    in time a priori bounds based on the energy estimates (and their local in time versions) are used in our work.  These derivations are standard, but we propose several interesting tricks to make the bounds as small as possible, their derivation is presented in Appendix \ref{trapping_sets}.

\item[STAGE 2]
  We reduce in appropriate way the original PDE to the finite dimensional problem governed by a system of ODEs. The infinite dimensional residual term appears in the reduced problem so we need to estimate it by the size of the mesh with the use of the a priori estimates, $W$, obtained in the first stage. This was the reason we needed the bounds up to the order $s$. With these estimates the problem is reduced to the Ordinary Differential Inclusion. The resultant inclusion is solved rigorously using the algorithm from \cite{kapela} implemented in CAPD library for rigorous numerics.
\end{itemize}
For the purpose of the proof of the periodic solution existence,  we verify if the set obtained after the period of integration  is the subset of the set of initial data. Since this is the case, by the Schauder theorem for the mapping of the forward in time translation by the period, we obtain the periodic solution existence. Obtained periodic orbit is stable which makes it possible to prove its existence by the Schauder theorem. Proposed rigorous integration procedure can be used to construct periodic solutions which have finite number of unstable directions if the Schauder theorem is replaced with different topological tools, such as covering relations, which make it possible to deal with expanding variables. Such technique was already successfully used  with the Fourier basis in \cite{kuramoto_III, zgliwi}.

As the reader familiar with our previous work \cites{mizgli, enclosiure_ref, kuramoto_III} might notice, the basic scheme of one step of the method is the same as in the periodic boundary conditions case. Let us highlight the differences.
\begin{itemize}
\item In STAGE 1
\begin{description}
\item In the periodic case, due to the isolation property we were able to obtain the bounds for each Fourier coefficient.  This can be easily automatized and a general enclosure algorithm can be given. This is accomplished by some standard ODE-type reasonings.
\item In the FEM-case, we work with various Sobolev norms, very much in the spirit of the modern theory of PDEs \cites{Robinson, Temam}.
  This part of the algorithm relies on various tricks and is technically much more involved than in the periodic case.
\end{description}
\item In STAGE 2
\begin{description}
\item In the periodic case, it is straightforward to obtain the ordinary differential inclusion.
\item In the FEM-case, obtaining the ordinary differential inclusion  requires that  our a-priori  bounds $W$ from STAGE 1 are
  in $H^s$ for $s$ sufficiently large. This allows us treat the error contributions coming from the Laplacian and various nonlinear terms.
\end{description}
\end{itemize}
\smallskip

We note that due to $H^s$ space regularity of the solution, in order for proposed techniques to work, we need certain smoothness of $f$ with respect to the space variable. This smoothness, which is inevitable in our approach, is higher than it is typically needed for error estimates of FEM.

 The computation times for the simple problem \eqref{a}--\eqref{b} are rather long (around 1 hr). This is mainly due to the fact that we used the first order FEM-elements.  There is no doubt that using a higher order FEM will greatly improve the
performace, but this requires  local-in-time a-priori estimates in $H^s$ for $s\geq 4$. In our paper we developed such estimates by hand
for $s=1,\dots,4$, for each $s$ separately. It will be good to have an algorithm, which will do it for us for any $s>0$.     

\medskip

We stress that for considered problem \eqref{a}--\eqref{b} one can very likely provide the computer assisted proof of existence of the periodic orbit using the sine Fourier basis, which would avoid some difficulties we cope with. However, in contrast to the Fourier approach, our approach, in principle, can be generalized to problems with multiple space dimensions for some class of space domains, where the construction of FEM basis is relatively simple. Since a priori estimates depend on the results such as the Sobolev embedding theorem or interpolation inequalities, which are dimension dependent, we would need to derive different estimates to deal with multidimensional problems, however, our algorithm of rigorous integration forward in time would stay the same. We remark that, as we work with smooth solutions, we would either need the appropriate smoothness of the space domain boundary or employ some additional techniques to deal with nonsmooth domains. Since the computation time of our computer-assisted proof is quite large (around 1 hr), to deal even with 2D domains we would need to optimize the code, and here it seems particularly promising to use higher order elements. This is the goal of our further work.

\medskip

\noindent \textbf{Other known approaches to rigorous numerics for PDEs.}
There is an approach called \emph{functional  analytic} as opposed to our, which can be termed \emph{topological} or \emph{geometric}.

It is based on the Newton method and some fixed point theorem, the problem under consideration is written as
\begin{equation}\label{eq:PDEgen}
\mathcal{F}(u) = 0,
\end{equation}
where $\mathcal{F}:X\to Y$ is a mapping between Banach spaces $X$ and $Y$. The map $F$ and spaces $X,Y$ encode the boundary conditions
and when looking for periodic orbits also periodic boundary conditions in time direction.

The method has been used successfully to verify the solutions for elliptic (using FEM or Fourier basis) \cites{galemi, gale1, gale2, plum_1, plum_2, plum_3} and parabolic  \cites{mizoguchi, takayasu, mitaku, nakiki, nakao, kikina} PDEs (but this mainly for periodic boundary conditions). The up to date information about these techniques can be found in the recent monograph \cite{plum_book}. The approach for the parabolic case  presented there relies on the inversion of the space time differential operator of the problem. The technique can be used to identify the invariant objects belonging to the global attractor \cites{castelli, figueras, gameiro}. In particular it has been used recently to verify the existence of time periodic solution for the Taylor--Green problem for the 2D incompressible Navier--Stokes equations \cite{Brendan}, where the rigorous calculations using the Fourier basis have been realized.

\medskip

\noindent \textbf{The plan of the article.}
We conclude the introduction with the brief presentation of the scheme of the paper: in Section \ref{sec:2} the problem is defined and some of its basic properties such as the existence and uniqueness of the weak and strong solutions, their continuous dependence on the initial data and basic energy estimates are recalled. Section \ref{sec:trapping} is devoted to the result on the existence of trapping sets for the considered problem. Next, Section \ref{sec:gal-tail} is devoted to the algorithm of reduction of  the original PDE to the Ordinary Differential Inclusion, and the algorithm of the computer assisted proof is described. Some details concerning the implementation are presented in Section \ref{sec:4}. Section \ref{sec:alg-proof} contains the theorem on periodic orbit existence and two examples of computations. The technical mathematical part of the paper, in particular the derivation of higher order energy estimates, which is simple but cumbersome, is contained in Appendixes  \ref{sec:appendix_A} and \ref{trapping_sets}.

\section{Burgers equation: problem setting and basic properties}\label{sec:2}

In this section we provide basic facts and results concerning our model problem: the nonautonomous Burger's equation. As a space domain we always consider the interval $\Omega = (0,1)$. We will use the shorthand notation for the spaces of functions defined on the interval $\Omega$, for example we will write simply $L^2$ in place of $L^2(\Omega)$, $H^1_0$ in place of $H^1_0(\Omega)$ and so on. The norm in $H^1_0$ is defined as $\|u_x\|_{L^2} = \|u\|_{H^1_0}$. Scalar product in $L^2$ will be denoted simply by $(\cdot,\cdot)$ and duality pairing between the space $H^1_0$ and $H^{-1}$, its dual, will be denoted by $\langle\cdot,\cdot \rangle_{H^{-1}\times H^1_0}$. We will also denote in a simplified way, by dropping the time variable, the spaces of functions leading from $\mathbb{R}$ or its subinterval to some spaces of space dependent functions,  for example  $L^\infty(L^2)$ will be the abbreviation for $L^\infty(\mathbb{R};L^2)$. For $I\subset \R$ and a Banach space $X$ we will denote by $L^p_{loc}(I;X)$ the space of functions which belong to $L^p(J;X)$ for any compact $J\subset I$. We stress that all proofs of this section are standard and they use well known techniques based on the energy estimates. We include them only for the exposition completeness. We define the initial time as $t_0\in \mathbb{R}$. We are interested in solving the  problem  governed by the one-dimensional Burgers equation \eqref{a} with the boundary condition \eqref{b} and the initial condition $u(t_0) = u_0$. We will always assume that the forcing term $f$ is defined for every $t\in \mathbb{R}$.

\subsection{Assumptions on $f$ and regularity of the solution $u$.} We define the space $Y = \{ u\in H^4\cap H^1_0\,:\ u_{xx}\in H^1_0 \}$. We assume that the non-autonomous forcing term has regularity $f\in L^\infty(Y)$ and we will consider the solutions of the initial-value problems governed by \eqref{a}--\eqref{b} with the initial data belonging to some subsets of $Y$. With this regularity of $f$, the corresponding solution $u$ belongs to $L^\infty(Y)\cap L^2(H^5)$. We stress that although it appears to us that the smoothness of $f$ is inevitable in our approach, we expect that, with some modifications, we could cope with more general situation with no boundary conditions on $f$ and $f_{xx}$, cf. Remarks \ref{ref:remark_zero} and \ref{rem:45}. Key role in our argument will be played by a priori estimates in $L^2$ norm of the solution and its space derivatives up to fourth. In this section as well as in the Appendix \ref{trapping_sets} we will gradually increase the assumed conditions on smoothness of $f$ (up to $L^\infty(Y)$) and we will present  the corresponding results, starting from the existence of the weak solution, and ending with required regularity and a priori estimates.

\subsection{Weak solution and its properties.}

The weak solution for the considered problem is defined as follows
\begin{definition}\label{ref:def_burg}
The function  $u\in L^2_{loc}([t_0,\infty);H^1_0)$ with $u_t\in L^2_{loc}([t_0,\infty);H^{-1})$ is a weak solution of the Burgers equation with the initial data $u(t_0) = u_0$ if the following equation holds
\begin{equation}\label{eq:equation1}
\langle u_t , v \rangle_{H^{-1} \times H^1_0} + (u_x, v_x) + (u u_x,v) = (f(\cdot,t),v) \quad \textrm{for every} \quad v\in H^1_0\quad \textrm{for almost every}\quad t>t_0.
\end{equation}
\end{definition}

Note that in above definition the time derivative $u_t$ is understood in distributional sense. Moreover, since $u$ is continuous as the function of time with values in $L^2$, the initial condition makes sense. The proof of the following result is standard and it follows by the Galerkin method. It is omitted here, but the details of the Galerkin  technique  for semilinear problems can be found, for example, in \cite[Chapters 8 and 9]{Robinson} or \cite[Section III.1.1.4]{Temam}.

\begin{theorem}
	Suppose that $f\in L^2_{loc}(L^2)$ and $u_0\in L^2$. Then the problem given by Definition \ref{ref:def_burg} has a unique weak solution.
	\end{theorem}
We derive the energy estimate satisfied by every weak solution of the above problem.
\begin{lemma}
	Let $u_0\in L^2$ and $f\in L^2_{loc}(L^2)$ and let $u$ be the weak solution corresponding to initial data $u_0$ taken at time $t_0$ and $f$. The following bounds are valid
	\begin{align}
	& \frac{d}{dt}\|u(t)\|_{L^2}^2 + 2\pi^2\|u(t)\|_{L^2}^2 \leq  2\|f(t)\|_{L^2}\|u(t)\|_{L^2}\quad \textrm{for almost every}\quad t>t_0 ,\label{est:1}\\
& 	\int_{t_0}^t\|u_x(s)\|_{L^2}^2\, ds \leq \|u_0\|^2_{L^2} + \frac{1}{\pi^2}\int_{t_0}^t\|f(s)\|^2_{L^2}\, ds\quad \textrm{for  every}\quad t>t_0.
\label{est:2}
	\end{align}
\end{lemma}
\begin{proof}
	The proof is standard.  We test \eqref{eq:equation1} by $u$. Note that the regularity of the weak solution guarantees that $\langle u_t(t), u(t)\rangle_{H^{-1}\times H^1_0} = \frac{1}{2}\frac{d}{dt}\|u(t)\|^2$ for almost every $t> 0$, cf. \cite[Proposition 23.23]{Zeidler2A}. From Lemma \ref{lemma:zero} we obtain
	\begin{equation}
\label{eq:ener-eq}
	\frac{1}{2}\frac{d}{dt}\|u\|_{L^2}^2 + \|u_x\|_{L^2}^2 = (f(t),u).
	\end{equation}
	Now, \eqref{est:1} follows by the Schwarz and Poincar\'e inequalities (see Lemma~\ref{lem:Poincare}).
	 On the other hand, integrating \eqref{eq:ener-eq}, we obtain
	$$
	\int_{t_0}^t\|u_x(s)\|_{L^2}^2\, ds \leq \frac{1}{2}\|u_0\|^2_{L^2} + \int_{t_0}^t \|f(s)\|_{L^2}\|u(s)\|_{L^2}\, ds\ \ \textrm{for a.e.}\ \ t>0.
	$$
   From this inequality, after using  the Poincar\'{e} and Cauchy inequalities we obtain
    \begin{align*}
      &\int_{t_0}^t\|u_x(s)\|_{L^2}^2\, ds \leq \|u_0\|^2_{L^2} + 2 \int_{t_0}^t \|f(s)\|_{L^2}\|u(s)\|_{L^2}\, ds - \int_{t_0}^t\|u_x(s)\|_{L^2}^2\, ds  \\
      &\qquad \leq \|u_0\|^2_{L^2} +  \int_{t_0}^t \left(\frac{\|f(s)\|^2_{L^2}}{\pi^2} + \pi^2 \|u(s)\|^2_{L^2}\,\right) ds - \pi^2 \int_{t_0}^t\|u(s)\|_{L^2}^2\, ds \\
& \qquad \qquad       = \|u_0\|^2_{L^2} +  \frac{1}{\pi^2} \int_{t_0}^t \|f(s)\|^2_{L^2}\, ds.
    \end{align*}
	The proof is complete.
\end{proof}

The mapping that assigns to the initial data taken at time $t_0$ the value of the solution at time $t \geq t_0$ will be denoted by $S(t,t_0):L^2\to L^2$.
Clearly $S(t,t_0)$ is a process, i.e. $S(t,t_1)S(t_1,t_0) = S(t,t_0)$ for every $t_0\leq t_1\leq t$ and $S(t_0,t_0) = I$.

\subsection{Strong solution and its properties.} We give the definition of the strong solution for the considered problem.
\begin{definition}\label{ref:def_burg_strong}
	The function $u\in L^2_{loc}([t_0,\infty);H^1_0\cap H^2)$ with $u_t\in L^2_{loc}([t_0,\infty);L^2)$ is the strong solution of the Burgers equation with the initial data $u(t_0) = u_0$ if there holds
	\begin{equation}\label{eq:equation1_strong}
	u_t  - u_{xx} + u u_x = f(\cdot,t) \quad \textrm{holds in} \quad L^2\quad \textrm{for a.e.}\quad t>t_0.
	\end{equation}
\end{definition}

As in case of weak solutions time derivative in above definition is distributional, and the initial conditions makes sense due to time continuity of the solution.  The proof of the following result is standard and we omit it here \cites{Robinson, Temam}.
\begin{theorem}
	Suppose that $f\in L^2_{loc}(L^2)$ and $u_0\in H^1_0$. Then the problem given by Definition \ref{ref:def_burg_strong} has a unique strong solution.
\end{theorem}
It is clear that a strong solution is also a weak solution, so the process $S(t,t_0)$ applied to an element of $u_0\in H^1_0$ defines the value of a strong solution at time $t$ if the initial data $u_0$ is taken at time $t_0$. The following result provides the energy estimate satisfied by the strong solutions.
\begin{lemma}\label{lem:est2}
	Let $u_0\in H^1_0$ and $f\in L^2_{loc}(L^2)$ and let $u$ be the strong solution corresponding to $u_0$ taken at time $t_0$ and $f$. Let $\alpha, \beta > 0$ be two constants such that $\alpha+\beta < 2$. The following differential inequalities hold for a.e. $t > t_0$
	 	 	\begin{equation}
	 \label{eq:strong_for_trapping}
	 \frac{d}{dt}\|u_x\|_{L^2}^2  \leq -2\|u_{xx}\|_{L^2}\left(\|u_{xx}\|_{L^2} - \|f(t)\|_{L^2}  - \|u\|_{L^2}^{5/4} \|u_{xx}\|_{L^2}^{3/4} \right).
	 \end{equation}
	 \begin{equation}
	 \label{eq:strong_1}
	 \frac{d}{dt}\|u_x\|_{L^2}^2 + \pi^2(2-\alpha-\beta)\|u_{x}\|_{L^2}^2 \leq {\frac{1}{\alpha}}\|f(t)\|^2_{L^2}  + \frac{7^7}{2^{16}\beta^7}\|u\|_{L^2}^{10}.
	 \end{equation}

\end{lemma}
\begin{proof}
	We multiply \eqref{eq:equation1_strong} by $-u_{xx}$, whence we get the bound
		$$
	\frac{1}{2}\frac{d}{dt}\|u_x\|_{L^2}^2 + \|u_{xx}\|_{L^2}^2 \leq \|f(t)\|_{L^2}\|u_{xx}\|_{L^2} + \int_0^1|u||u_x||u_{xx}|\, dx.
	$$
	It follows that
	$$
	\frac{1}{2}\frac{d}{dt}\|u_x\|_{L^2}^2 + \|u_{xx}\|_{L^2}^2 \leq \|f(t)\|_{L^2}\|u_{xx}\|_{L^2} + \|u\|_{L^\infty}\|u_x\|_{L^2}\|u_{xx}\|_{L^2}.
	$$
	Using Lemma \ref{lem:int} we deduce  that
	$$
	\frac{1}{2}\frac{d}{dt}\|u_x\|_{L^2}^2 + \|u_{xx}\|_{L^2}^2 \leq \|f(t)\|_{L^2}\|u_{xx}\|_{L^2} + \|u\|_{L^2}^{5/4}\|u_{xx}\|_{L^2}^{7/4}.
	$$
	We obtain \eqref{eq:strong_for_trapping}.

 Let $\alpha,\beta>0$ be two constants. 	We use Lemma \ref{lem:young} with $p = 8$, $q = 8/7$, $\epsilon = (\beta 4/7)^{7/8}$  to estimate $\|u\|_{L^2}^{5/4}\|u_{xx}\|_{L^2}^{7/4}$ and with $p=q=2$, $\epsilon=\sqrt{\alpha}$ to estimate $\|f(t)\|_{L^2}\|u_{xx}\|_{L^2}$. We get
	$$
	\frac{1}{2}\frac{d}{dt}\|u_x\|_{L^2}^2 + \|u_{xx}\|_{L^2}^2 \leq \frac{1}{2\alpha}\|f(t)\|^2_{L^2} + \frac{\alpha}{2}\|u_{xx}\|_{L^2}^2 + \frac{7^7}{2^{17}\beta^7}\|u\|_{L^2}^{10} + \frac{\beta}{2}\|u_{xx}\|_{L^2}^2,
	$$
whence the following inequality holds
\begin{equation}\frac{d}{dt}\|u_x\|_{L^2}^2 + (2-\alpha-\beta)\|u_{xx}\|_{L^2}^2 \leq {\frac{1}{\alpha}}\|f(t)\|^2_{L^2}  + \frac{7^7}{2^{16}\beta^7}\|u\|_{L^2}^{10}.\label{eq:dux-uxx}
\end{equation} Using the Poincar\'e inequality we obtain \eqref{eq:strong_1}.
	The proof is complete.
\end{proof}
We pass to the proof of $H^1_0$ continuity of $S(t,t_0)$.
\begin{lemma}\label{lemma:110}
	Let $f\in L^2_{loc}(L^2)$ and let $t_0\in \R$. If $u_0,v_0\in H^1_0$ and $u,v$ are two strong solutions corresponding to $u_0, v_0$ taken at $t_0$, respectively, then
	$$\|u(t)-v(t)\|_{H^1_0} \leq e^{\frac{1}{2}(\|u_0\|_{L^2}^2+\|v_0\|_{L^2}^2)+\frac{1}{\pi^2}\int_{t_0}^t\|f(s)\|_{L^2}^2\, ds} \|u_0-v_0\|_{H^1_0} \quad \textrm{for every}\quad t\geq t_0.$$
\end{lemma}
\begin{proof}
	Let $u_0, v_0 \in H^1_0$ and let $u,v$ be strong solutions corresponding to $u_0, v_0$, respectively. Denoting $w=u-v$ there holds the following equation
	$$
	w_t - w_{xx}  + u u_x - vv_x  = 0 \quad  \textrm{a.e.}\quad t>t_0, x\in (0,1).
	$$
	Testing this equation by $-w_{xx}$, we obtain (using $uu_x - vv_x=uw_x + wv_x$)
	$$
	\frac{1}{2}\frac{d}{dt}\|w_x\|^2_{L^2} + \|w_{xx}\|^2_{L^2} \leq |(u w_x,w_{xx})| + |(v_xw,w_{xx})| .
	$$
	Using $\|w\|_{L^\infty} \leq \|w_x\|_{L^2}$ it follows that
	\begin{align*}
		&\frac{1}{2}\frac{d}{dt}\|w_x\|^2_{L^2} + \|w_{xx}\|^2_{L^2} \leq \|u\|_{L^\infty}\|w_x\|_{L^2}\|w_{xx}\|_{L^2} + \|w\|_{L^\infty}\|v_x\|_{L^2}\|w_{xx}\|_{L^2}  \\
      &\qquad  \leq\|u_x\|_{L^2}\|w_x\|_{L^2}\|w_{xx}\|_{L^2} + \|w_x\|_{L^2}\|v_x\|_{L^2}\|w_{xx}\|_{L^2}
	\end{align*}
	whence, as $\|w\|_{H^1_0} = \|w_x\|_{L^2}$,
	$$
		\frac{1}{2}\frac{d}{dt}\|w\|^2_{H^1_0} + \|w_{xx}\|^2_{L^2} \leq \left(\|u_x\|_{L^2}+\|v_x\|_{L^2}\right)\|w\|_{H^1_0}\|w_{xx}\|_{L^2}.
	$$
It follows that
	$$
	\frac{d}{dt}\|w\|^2_{H^1_0} \leq (\|u_x\|_{L^2}^2+\|v_x\|_{L^2}^2) \|w\|^{2}_{H^1_0}.
	$$
	The assertion follows by the Gronwall lemma and estimates \eqref{est:2}.
\end{proof}

\subsection{Asymptotic behavior of solutions.}
We start this section from the recollection of the definition of the eternal strong solution.
\begin{definition}\label{def:eternal}
	The function $u\in C(H^1_0)$ is called the eternal strong solution of the Burgers equation if there exists a constant $A>0$ such that $\|u_x(t)\|_{L^2}\leq A$ for every $t\in \R$ and for every $t_0\in \R$ the function $u$ restricted to $[t_0,\infty)$ is the strong solution given by Definition \ref{ref:def_burg_strong} with the initial data $u(t_0)$.
	\end{definition}
We recall the results of \cite{KalitaZgliczynski}.
\begin{theorem}(cf. \cite[Theorem 4.1]{KalitaZgliczynski})
	\label{thm:uniq_et}
	Let $f\in L^\infty(L^2)$. There exists a unique eternal strong solution in the sense of Definition \ref{def:eternal}.
\end{theorem}
 We use the notation $\mathcal{B}(L^2)$ for nonempty and bounded subsets of $L^2$ and $\|B\|_{L^2} = \sup_{v\in B}\|v\|_{L^2}$ for $B\in \mathcal{B}(L^2)$.
\begin{theorem}(cf. \cite[Theorem 4.7]{KalitaZgliczynski})\label{thm:attr}
	Let $f\in L^\infty(L^2)$. Let $v_0 \in B \in \mathcal{B}(L^2)$ and let $v$ be a weak solution starting from the initial data $v_0$ at time $t_0$. Let $u$ be the unique eternal solution given by Theorem \ref{thm:uniq_et}. There exists a constant $D=D(\|B\|_{L^2})>0$ (depending continuously and monotonically on $\|B\|_{L^2}$) and a constant $C>0$ dependent on $f$ but not on $v_0$ such that for every $t\geq t_0$ there holds
	$$
	\|u(t) - v(t)\|_{L^2} \leq D(\|B\|_{L^2}) e^{-C(t-t_0)}.
	$$
\end{theorem}
Since the right-hand side of the last estimate depends on $\|B\|_{L^2}$ and not on $v$ or $v_0$ we deduce that the estimate can be replaced with
\begin{equation}\label{eq:haus}
\mathrm{dist}_{L^2} (S(t,t_0)B,\{u(t)\}) \leq D(\|B\|_{L^2}) e^{-C(t-t_0)},
\end{equation}
where $\mathrm{dist}_{L^2}$ is the Hausdorff semidistance defined as
$$
\mathrm{dist}_{L^2}(A,B) = \sup_{a\in A}\inf_{b\in B} \|a-b\|_{L^2}.
$$
It is clear that the right-hand side of the estimate \eqref{eq:haus} tends to zero either as $t_0\to -\infty$ for fixed $t$ or as $t\to \infty$ as $t_0$ is fixed. This means that the unique eternal strong solution attracts all weak (and thus also strong) solutions uniformly with respect to bounded sets of the initial data. We also note that if the forcing is time-periodic, then the eternal strong solution $u$ also has to be time-periodic with the same period.

\section{Trapping sets}
\label{sec:trapping}
This section is devoted to construction of a set on which the process $S(t,t_0)$ is positively invariant (i.e. it is, in a sense, \textit{trapping} as once the trajectory enters it, it can never leave it). The results of this section will be used later when we pass to the numerical part of this article. We will later derive the error estimates of the Galerkin projection of the solution in $H^1_0$ and, to this end, we will need to control up to the \textit{fourth} space derivative of the solution. Hence, the results of this chapter provide the existence of positively invariant sets on which we control $L^2$ norms of the space derivatives
of the solution up to fourth order. Note that this requires stronger assumptions on the regularity of $f$.

  We start from the definition of the positively invariant (trapping) set.
  \begin{definition}
  	Let $X$ be a Banach space and let the family of maps $\{ S(t,t_0) \}_{t\geq t_0}$ be a process on $X$.
  	The set $B\subset X$ is said to be positively invariant (trapping) if
 for every $t_0\in \mathbb{R}$ and $t>0$ there holds $S(t_0+t,t_0)B\subset B$.
  \end{definition}
  In the following result we demonstrate the existence of a convex trapping  set being a closed and bounded subset of $H^4$. Since the proof of this result is technical, it is postponed until Appendix \ref{trapping_sets}.

  \begin{theorem}\label{thm:trapping}
  	Define $Y = \{ u\in H^4\cap H^1_0\,:\ u_{xx}\in H^1_0 \}$ endowed with the norm of $H^4$. Let $f\in L^\infty(Y)$. For the process $\{ S(t,t_0) \}_{t\geq t_0}$ defined by weak solutions of the problem \eqref{a}--\eqref{b} there exists a nonempty trapping set $\mathcal{B}_0 \subset Y$ which is convex, closed and bounded in $Y$. Moreover, it is possible to find explicitly the radii $R_1,R_2, R_3, R_4, R_5$ such that if $u\in \mathcal{B}_0$ then$$
  	\|u\|_{L^2}\leq R_1,\ \    	\|u_x\|_{L^2}\leq R_2,\ \   	\|u_{xx}\|_{L^2}\leq R_3,\ \
  	  	\|u_{xxx}\|_{L^2}\leq R_4,\ \
  	  	\|u_{xxxx}\|_{L^2}\leq R_5.
  	$$
  \end{theorem}
Note that due to assumed regularity of $f$ in above theorem one can replace weak solutions with the strong solutions.

\begin{remark}
\label{rem:r1-r5}
	The values of the radii $R_1$-$R_5$ will be calculated using the a priori estimates of the equation and then used in the computer assisted construction of the attracting trajectory. The construction will be based on the splitting of the whole space $H^1_0$ into a finite dimensional part and its infinite dimensional remainder which will be estimated using these
	radii. Then, the contribution of the remainder will be incorporated in a multivalued additive term thus leading to the need of the rigorous numerical solution of an ordinary differential inclusion.  From numerical reasons it will be crucial that the width of this inclusion is as small as possible. As this width depends on $R_1$-$R_5$ significant technical effort is put in Appendix \ref{trapping_sets} to construct the trapping set with smallest possible radii. In fact for each $R_i$ several algorithms originating from different a priori estimates are presented   in Appendix \ref{trapping_sets}. Then, in computational part for particular function $f$ we implement all algorithms and choose the smallest obtained value for each radius. For readers convenience all algorithms are summarised in Tab.\ref{tab:radii}.
\end{remark}

\begin{table}
	\begin{tabular}{|cc|}
		\hline
		Radius of trapping set & Algorithm  \\
		\hline
		\hline
		$R_1$ & Lemma \ref{lem:b1} \\
		\hline
		$R_2$ & Lemma \ref{lem:abssetL2}   \\
		& Formula \eqref{eq:alt} and Corollary \ref{cor:29}  \\
		\hline
		$R_3$ & Lemma \ref{lem:uxx_first} (two methods)\\
		& Lemma \ref{lem:xx} (three methods)\\
		\hline
		$R_4$ & Lemma \ref{lem:u3x-bnd} (two methods)\\
		& Lemma \ref{lem:xxx} (four methods)\\
		\hline
		$R_5$ & Lemma \ref{lem:u4x-bnd}\\
		& Lemma \ref{lem:u4x-bnd_2}\\
		& Lemma \ref{lem:xxxx}\\
		\hline
	\end{tabular}	
	\caption{Methods to derive the radii $R_1$--$R_5$ of the trapping sets. We always calculate the radii by all these methods (using the rigorous interval arithmetics) and choose the smallest obtained value.}\label{tab:radii}
\end{table}

\begin{remark}\label{ref:remark_zero}
	The restriction that $f$ as well as $f_{xx}$ should satisfy the Dirichlet condition comes from the fact that we need the homogeneous Dirichlet boundary conditions for $u_{xx}$ and $u_{xxxx}$, as we derive the energy estimates for the original equation to which we apply the second and fourth space derivatives. We avoid full generality which could be achieved by translating the second and fourth space derivatives of the solution $u$ by any function $a(x,t)$ which satisfies the same boundary conditions as $f$. Indeed, if $u$ solves \eqref{a}--\eqref{b}, then $v=u_{xx}$ solves
	$$
	v_t - v_{xx} + 3vu_{x}+uv_x = f_{xx},
	$$
	and $w = v+a =u_{xx}+a$ is guaranteed to satisfy the homogeneous Dirichlet condition if only $a=f$ in points $x=0$ and $x=1$. The function $w$ satisfies the following equation
	$$
	w_t - w_{xx} + 3wu_x+uw_x = f_{xx}+a_t-a_{xx}+3au_x+ua_x.
	$$
	One possible choice is $a=f$, provided $f$ is smooth enough. Then the resulting problem with homogeneous conditions has the form
	$$
	w_t - w_{xx} + 3wu_x+uw_x = f_t+3fu_x+uf_x,
	$$
	but any function which satisfies the same boundary conditions as $f$ can be used. Another possibility would be to pick $a = f(0,t)(1-x) + f(1,t)x$, which would lead us to the equation
	\begin{eqnarray*}
	w_t - w_{xx} + 3wu_x+uw_x &=& f_{xx}+f_t(0,t)(1-x) + f_t(1,t)x \\
                 &+& 3(f(0,t)(1-x) +f(1,t)x)u_x+u(f(1,t)-f(0,t)).
	\end{eqnarray*}
	The argument of finding the trapping sets and local estimates leading to Theorem \ref{thm:trapping} can be realized for the translated function $w$. Now, the trapping balls for $L^2$ norms of $w$ and $w_x$ will be centered at zero. This leads to the trapping balls centered at $-a$ and $-a_{x}$ for the derivatives $u_{xx}$ and $u_{xxx}$ of the solution of the original problem. Similar procedure can be done for $w_{xx}$ this time translating it by any function which satisfies the same boundary conditions as the right-hand side of the equation for $w$, thus leading to the equation with homogeneous Dirichlet conditions for $w_{xx}$. In the sequel, however, to avoid technicalities, we restrict ourselves to the situation where $f(0,t)=f(1,t) = f_{xx}(0,t) = f_{xx}(1,t) = 0$.

\end{remark}

Theorem \ref{thm:trapping} provides the estimates which hold on the trapping set for any time and the radii $R_1-R_5$ contained there depend on the norms of the forcing term and its space derivatives in $L^\infty(\R;L^2)$. As the non-autonomous term varies with time, and for some $t$ the norms of $f(t)$ or its space derivatives can be very small, in Appendix  \ref{trapping_sets} we also derive local estimates, which are given in the following lemma
\begin{lemma}\label{lemma_local}
	Let $u \in Y$ 
	be such that
	 $$
	 \|u\|_{L^2}\leq R^i_1,\ \    	\|u_x\|_{L^2}\leq R^i_2,\ \   	\|u_{xx}\|_{L^2}\leq R^i_3,\ \
	 \|u_{xxx}\|_{L^2}\leq R^i_4,\ \
	 \|u_{xxxx}\|_{L^2}\leq R^i_5,
	 $$
	 with $R^i_1$-$R^i_5>0,$ and let $t_{i+1}>t_i$.
	 There exist positive numbers ${M}^{i+1}_1$, ${M}^{i+1}_2$, ${M}^{i+1}_3$,
	 ${M}^{i+1}_4$, ${M}^{i+1}_5$, which can be calculated explicitly, such that
	 \begin{align*}	
	 & \|S(t,t_i)u\|_{L^2}\leq M^{i+1}_1,\ \    	\|(S(t,t_i)u)_x\|_{L^2}\leq M^{i+1}_2,\ \   	\|(S(t,t_i)u)_{xx}\|_{L^2}\leq M^{i+1}_3,\\
	  &\qquad
	 \|(S(t,t_i)u)_{xxx}\|_{L^2}\leq M^{i+1}_4,\ \
	 \|(S(t,t_i)u)_{xxxx}\|_{L^2}\leq M^{i+1}_5\quad \textrm{for every}\ \ t\in [t_i,t_{i+1}],
	 \end{align*} and positive numbers ${R}^{i+1}_1$, ${R}^{i+1}_2$, ${R}^{i+1}_3$,
	 ${R}^{i+1}_4$, ${R}^{i+1}_5$, which can be calculated explicitly, such that  \begin{align*}	
	 & \|S(t_{i+1},t_i)u\|_{L^2}\leq {R}^{i+1}_1,\ \    	\|(S(t_{i+1},t_i)u)_x\|_{L^2}\leq {R}^{i+1}_2,\ \   	\|(S(t_{i+1},t_i)u)_{xx}\|_{L^2}\leq R^{i+1}_3,\\
	 &\qquad
	 \|(S(t_{i+1},t_i)u)_{xxx}\|_{L^2}\leq R^{i+1}_4,\ \
	 \|(S(t_{i+1},t_i)u)_{xxxx}\|_{L^2}\leq R^{i+1}_5.
	 \end{align*}
	\end{lemma}
We stress that the above lemma is valid if we set as $R^i_k$, $R^{i+1}_k$ and $M_k^{i+1}$ for $k\in \{ 1,\ldots, 5\}$ the global radii $R_k$. The point of the above lemma is the construction, presented in Appendix \ref{trapping_sets} of the optimal local radii $R^{i+1}_k$ and $M^{i+1}_k$, which in practice are often much smaller than the corresponding radii on the trapping set. Such localization of the estimates  allows to make the width of the constructed differential inclusion more narrow, which is crucial from the numerical point of view.

\begin{remark}
	Local bounds $M_k^i$ and $R_k^i$ are derived for $k=1,\ldots,5$ based on energy inequalities which can be derived in several possible ways, each of them leading to the different bound. For every timestep of the simulation we calculate those values using all derived estimates and always choose the smallest obtained bound. All energy estimates used to derived the bounds are summarized, for reader's convenience in Tab.\ref{local_boundsi}. We stress that the propagation of a priori estimates could depend only on the global trapping radii $R_i$, but incorporating the local in time estimates allows us to make the width of the inclusion more narrow which is crucial for numerical reasons.
	\end{remark}
\begin{table}
	\begin{tabular}{|cc|}
		\hline
		$k$ in the local estimate & result used in the computation \\
		\hline
		\hline
		$k=1$ & inequality \eqref{eq:R1_est} \\
		\hline
		$k=2$ & Lemma \ref{lemma:young_x}   \\
		& Lemma \ref{lemma:young_x_2}  \\
		\hline
		$k=3$ & inequalities \eqref{eq:uxx_1} and \eqref{eq:uxx_2} (two methods)\\
		& Lemma \ref{lemma:young_x_3} (two methods)\\
		\hline
		$k=4$ & Lemma \ref{lem_3.18} (two methods)\\
		& Lemma \ref{lemma:young_x_4} (two methods)\\
		\hline
		$k=5$ & Lemma \ref{lem:u4}\\
		& Lemma \ref{lemma:young_x_5}\\
		\hline
	\end{tabular}	
	\caption{Methods to derive the local estimates in the time intervals and on their endpoints. Again for all time steps all estimates are calculated (using the rigorous interval arithmetics) and always the smallest one is chosen. }\label{local_boundsi}
\end{table}

\section{Rigorous integration algorithm based on the Galerkin projection in FEM space}
\label{sec:gal-tail}
\subsection{Galerkin projection of first order and its basic properties.} We define $V_k$ as the subspace of $H^1_0$ of the functions which are linear on intervals
$$
\left(\frac{i}{k},\frac{i+1}{k}\right)\quad \textrm{for}\quad i\in \{0,\ldots,k-1\}.
$$
 The dimension of the space $V_k$ is equal to $k-1$, we denote the length of the mesh interval as $h_k = \frac{1}{k}$. We split any $u \in H^1_0$ into $P_k u$  (orthogonal projection in $H^1_0$ on $V_k$) and $Q_k u = u - P_ku$. Observe, that $P_k$ coincides with the piecewise linear interpolation
 operator.

There hold the following inequalities
\begin{align}
&\|Q_k u\|_{H^1_0} \leq \frac{h_k}{\pi}\|u_{xx}\|_{L^2} \quad \textrm{for every}\quad u\in H^1_0 \cap H^2,\label{eq:estimate_1}\\
&\|Q_k u\|_{L^2} \leq \frac{h_k}{\pi}\|Q_k u\|_{H^1_0} \quad \textrm{for every}\quad u\in H^1_0.\label{eq:estimate_2}
\end{align}
The proof of \eqref{eq:estimate_1} can be found in \cite[Theorem 2.5]{Schultz}, also c.f., \cite{Nakao1}. The estimate \eqref{eq:estimate_2} is a direct consequence of the fact that the function $Q_ku$ vanishes in the nodes of the mesh and the Poincar\'e inequality given in Lemma \ref{lem:Poincare}, also c.f., \cite[Theorem 1.2]{Schultz}. 
Note, that, in fact \eqref{eq:estimate_2} holds  in subspace $Q_k H^1_0$, i.e.,
\begin{align}
&\|u\|_{L^2} \leq \frac{h_k}{\pi}\|u\|_{H^1_0} \quad \textrm{for every}\quad u\in Q_k H^1_0.\label{eq:new-estimate_2}
\end{align}

With estimates \eqref{eq:estimate_1} and \eqref{eq:estimate_2} we deduce a very simple lemma.

\begin{lemma}\label{lem:galerkin_est}
	Suppose that a function $u:[t_1,t_2]\to H^2\cap H^1_0$ is such that
	$$
	\| u_{xx}(t)\|_{L^2} \leq R\quad \textrm{for every}\quad t\in [t_1,t_2].
	$$
	Then the following estimates hold
	\begin{align}
	& \|Q_k u(t) \|_{H^1_0}\leq  \frac{h_k}{\pi} R\quad \textrm{for every}\quad t\in [t_1,t_2],\\
	& \|Q_k u(t) \|_{L^2}\leq  \frac{h_k^2}{\pi^2} R\quad \textrm{for every}\quad t\in [t_1,t_2].
	\end{align}
\end{lemma}

\subsection{Equation satisfied by $P_ku$.}\label{sec:residue}Let $u$ be the solution of the problem \eqref{a}--\eqref{b} confined in the trapping set $\mathcal{B}_0$. Splitting $u(t) = P_ku(t) + Q_ku(t)$ we obtain
\begin{align*}
&(P_ku_t, v) + (Q_ku_t, v) + ((P_ku)_x,v_x) + ((Q_ku)_x,v_x) + (P_ku(P_ku)_x,v) + (Q_ku u_x,v)\\
& \qquad \qquad  + (P_ku(Q_ku)_x,v)  = (f(t),v) \quad \textrm{for every} \quad v\in V_k.
\end{align*}
Noting that $((Q_ku)_x,v_x) =0$, we can rewrite this equation as
\begin{align}
(P_ku_t, v) + &((P_ku)_x,v_x) + (P_ku(P_ku)_x,v) = (f(t),v)  \nonumber\\
   - &\left[(Q_ku_t, v) + (Q_ku u_x,v) + (P_ku(Q_ku)_x,v)  \right]\quad \textrm{for every} \quad v\in V_k.\label{eq:Galerkin_proj}
\end{align}
The equation
$$
(P_ku_t, v) + ((P_ku)_x,v_x) + (P_ku(P_ku)_x,v) = (f(t),v)
$$
corresponds to the Galerkin scheme for the considered problem, while the expression $$\left[(Q_ku_t, v) + (Q_ku u_x,v) + (P_ku(Q_ku)_x,v)  \right]$$ is the residual error which we need to estimate.
Denote by $\{ v^i \}_{i=1}^{k-1}$ the basis functions of $V_k$ defined by the relation $v^i(jh_k)=\delta_{ij}$. Then representing $P^ku(t)$ as $P_ku(t) = \sum_{i=1}^{k-1}\alpha_i(t)v^i(x)$, we will formulate a differential inclusion satisfied by the time dependent coefficients $\alpha_i$. The system takes the form
\begin{align*}
&\sum_{i=1}^{k-1}\alpha'_i(t)(v^i, v) + \sum_{i=1}^{k-1}\alpha_i(t)(v^i_x,v_x) + \sum_{i=1}^{k-1}\sum_{j=1}^{k-1}\alpha_i(t)\alpha_j(t)(v^i v^j_x,v) = (f(t),v)  \\
& \quad \quad  - \left[(Q_ku_t, v) + (Q_ku u_x,v) +  (P_ku(Q_ku)_x,v) \right]\quad \textrm{for every} \quad v\in V_k.
\end{align*}
Let us calculate the $m$-th equation of the system ($m\in \{1,\ldots,k-1\}$) obtained by taking as the test function $v$ the $m$-th element of the basis
\begin{align*}
&\sum_{i=1}^{k-1}\alpha'_i(t)(v^i, v^m) + \sum_{i=1}^{k-1}\alpha_i(t)(v^i_x,v^m_x) + \sum_{i=1}^{k-1}\sum_{j=1}^{k-1}\alpha_i(t)\alpha_j(t)(v^i v^j_x,v^m) = (f(t),v^m)  \\
& \quad \quad  - \left[(Q_ku_t, v^m) + (Q_ku u_x,v^m) +  (P_ku(Q_ku)_x,v^m)  \right].
\end{align*}
We only provide the details of further calculations for $m\neq 1$ and $m\neq k-1$. For a given $m$ there will be three $i$-s for which the contribution to the sums will be nonzero. Their associated formula is the following
\begin{align*}
&\alpha'_m(t)(v^m, v^m) + \alpha'_{m-1}(t)(v^{m-1}, v^m)  + \alpha'_{m+1}(t)(v^{m+1}, v^m)\\
&\qquad  + \alpha_m(t)(v^m_x,v^m_x) + \alpha_{m-1}(t)(v^{m-1}_x,v^m_x) + \alpha_{m+1}(t)(v^{m+1}_x,v^m_x)\\
& \qquad + \sum_{j=1}^{k-1}\alpha_m(t)\alpha_j(t)(v^m v^j_x,v^m) +
\sum_{j=1}^{k-1}\alpha_{m-1}(t)\alpha_j(t)(v^{m-1} v^j_x,v^m)
+\sum_{j=1}^{k-1}\alpha_{m+1}(t)\alpha_j(t)(v^{m+1} v^j_x,v^m) \\
& \qquad = (f(t),v^m)   - \Big[(Q_ku_t, v^m) + (Q_ku u_x,v^m) +  (P_ku(Q_ku)_x,v^m)  \Big].
\end{align*}
We further elaborate three remaining sums by keeping the nonzero terms
\begin{align*}
&\alpha'_m(t)(v^m, v^m) + \alpha'_{m-1}(t)(v^{m-1}, v^m)  + \alpha'_{m+1}(t)(v^{m+1}, v^m)\\
&\qquad  + \alpha_m(t)(v^m_x,v^m_x) + \alpha_{m-1}(t)(v^{m-1}_x,v^m_x) + \alpha_{m+1}(t)(v^{m+1}_x,v^m_x)\\
& \qquad + \alpha^2_m(t)(v^m v^m_x,v^m) + \alpha_m(t)\alpha_{m-1}(t)(v^m v^{m-1}_x,v^m)  + \alpha_m(t)\alpha_{m+1}(t)(v^m v^{m+1}_x,v^m) \\
&\qquad  +
\alpha_{m-1}(t)\alpha_m(t)(v^{m-1} v^m_x,v^m) + \alpha^2_{m-1}(t)(v^{m-1} v^{m-1}_x,v^m)\\
& \qquad +
\alpha_{m+1}(t)\alpha_m(t)(v^{m+1} v^m_x,v^m) + \alpha^2_{m+1}(t)(v^{m+1} v^{m+1}_x,v^m) \\
& \qquad \qquad = (f(t),v^m)   - \Big[(Q_ku_t, v^m) + ((Q_ku) u_x,v^m) +  (P_ku(Q_ku)_x,v^m) \Big].
\end{align*}
It is clear that $(v^m v^m_x,v^m) = 0$. Hence, calculating all integrals in the above formula, we get
\begin{align*}
&\alpha'_m(t)\frac{2h_k}{3}  + \alpha'_{m-1}(t)\frac{h_k}{6}  + \alpha'_{m+1}(t)\frac{h_k}{6}  + \alpha_m(t) \frac{2}{h_k} - \alpha_{m-1}(t)\frac{1}{h_k} - \alpha_{m+1}(t)\frac{1}{h_k}\\
&\qquad
- \alpha_m(t)\alpha_{m-1}(t)\frac{1}{6}
+ \alpha_m(t)\alpha_{m+1}(t)\frac{1}{6}
- \alpha^2_{m-1}(t)\frac{1}{6}
+ \alpha^2_{m+1}(t)\frac{1}{6} \\
& \qquad = (f(t),v^m)   - \Big[(Q_ku_t, v^m) + (Q_ku u_x,v^m) +  (P_ku(Q_ku)_x,v^m) \Big].
\end{align*}
We can integrate by parts
$$
(P_ku(Q_ku)_x,v^m)  = -((P_ku)_x Q_ku,v^m)-(P_ku Q_ku,v^m_x),
$$
whence the system takes the form
\begin{align*}
&\alpha'_m(t)\frac{2 h_k}{3} + \alpha'_{m-1}(t)\frac{h_k}{6} + \alpha'_{m+1}(t)\frac{h_k}{6} + \alpha_m(t)\frac{2}{h_k} - \alpha_{m-1}(t)\frac{1}{h_k} - \alpha_{m+1}(t)\frac{1}{h_k}\\
&\qquad
- \alpha_m(t)\alpha_{m-1}(t)\frac{1}{6}
+ \alpha_m(t)\alpha_{m+1}(t)\frac{1}{6}
- \alpha^2_{m-1}(t)\frac{1}{6}
+ \alpha^2_{m+1}(t)\frac{1}{6} \\
& \qquad = (f(t),v^m)   - \Big[ (Q_ku_t, v^m) + (Q_ku Q_ku_x,v^m) -  (P_ku Q_ku,v^m_x)  \Big].
\end{align*}
%
We multiply this equation by $ \frac{6}{h_k}$ to get
\begin{align*}
&4\alpha'_m(t) + \alpha'_{m-1}(t) + \alpha'_{m+1}(t) +  \frac{12}{h_k^2} \alpha_m(t)  -  \frac{6}{h_k^2}\alpha_{m-1}(t) - \frac{6}{h_k^2}\alpha_{m+1}(t)\\
&\qquad
- \frac{1}{h_k} \alpha_m(t)\alpha_{m-1}(t)
+ \frac{1}{h_k} \alpha_m(t)\alpha_{m+1}(t)
- \frac{1}{h_k}\alpha^2_{m-1}(t)
+ \frac{1}{h_k}\alpha^2_{m+1}(t) \\
& \qquad = \frac{6}{h_k}(f(t),v^m)   -  \frac{6}{h_k}\left[(Q_ku_t, v^m) - (Q_ku Q_ku_x,v^m) +  (P_ku Q_ku,v^m_x) \right],
\end{align*}
for $m\in \{ 2,\ldots,k-2\}$. Together with the equations for $m=1$ and $m=k-1$ (not given here, they are derived analogously)
the system can be rewritten in the matrix form as
$$
M \alpha'= S \alpha + N(\alpha) + F(t)+ \frac{6}{h_k}((Q_ku Q_ku_x-Q_ku_t,v^m) -  (P_ku Q_ku,v^m_x) )_{m=1}^{k-1},
$$
where $F(t) = (F_m(t))_{m=1}^{k-1}$ is given by $F_m(t) = \frac{6}{h_k}(v^m, f(t))$, $M$ is the mass matrix (which multiplies all terms with time derivatives of $\alpha$), $S$ is the stiffness matrix (which constitutes the linear term multiplying $\alpha$), and $N$ is the expression coming from the quadratic terms.
We multiply this equation by $M^{-1}$. Then
$$
\alpha'= M^{-1}S \alpha + M^{-1}N(\alpha) + M^{-1}F(t)+ M^{-1}\frac{6}{h_k}((Q_ku Q_ku_x-Q_ku_t,v^m) -  (P_ku Q_ku,v^m_x) )_{m=1}^{k-1}.
$$
Now let $B$ be a nonsingular square matrix and introduce the new variables $\beta$ given by $\alpha = B\beta$. Then $B$ is the change of basis matrix.  In new variables the equation takes the form
\begin{align*}
& \qquad \beta'(t)= B^{-1}M^{-1}S B\beta(t) + B^{-1}M^{-1}N(B\beta(t)) + B^{-1}M^{-1}F(t)\\
& \qquad + B^{-1}M^{-1}\frac{6}{h_k}((Q_ku(t) Q_ku_x(t)-Q_ku_t(t),v^m) -  (P_ku(t) Q_ku(t),v^m_x) )_{m=1}^{k-1}.
\end{align*}
The matrix $B$ is found in such a way that the matrix $ B^{-1}M^{-1}S B$ is close to diagonal.
Now if we denote the coefficients of the matrix $B^{-1}M^{-1}$ by $c_{lm}$ for $l, m\in \{1,\ldots, k-1 \}$, and we denote
\begin{equation}
w^l(x) = \sum_{m=1}^{k-1} c_{lm} v^m(x), \label{eq:w^l-expansion}
\end{equation}
 then the above equation is equivalent to
\begin{align*}
& \qquad \beta'(t)= B^{-1}M^{-1}S B\beta(t) + B^{-1}M^{-1}N(B\beta(t)) + B^{-1}M^{-1}F(t)\\
& \qquad \qquad + \frac{6}{h_k}((Q_ku(t) Q_ku_x(t)-Q_ku_t(t),w^l) -  (P_ku(t) Q_ku(t),w^l_x) )_{l=1}^{k-1}.
\end{align*}
We will solve the differential inclusion
\begin{equation}\label{inclusion}
\beta'(t) \in B^{-1}M^{-1}S B\beta(t) + B^{-1}M^{-1}N(B\beta(t)) + B^{-1}M^{-1}F(t)+ G (t),
\end{equation}
where $G:\mathbb{R}\to (2^{\mathbb{R}})^{k-1}$, and  $G(t) = (G_l(t))_{l=1}^{k-1}$ with $G_l(t)$ being a subset (actually an interval) in $\R$, such that
\begin{equation}\label{G_included}
\frac{6}{h_k}\left((Q_ku(t) Q_ku_x(t)-Q_ku_t(t),w^l) -  (P_ku(t) Q_ku(t),w^l_x)\right) \in G_l(t) \quad \textrm{for}\quad l=1,\ldots, k-1.
\end{equation}

\subsection{A numerical algorithm}\label{subsec:numerical}
We divide the interval $[t_0,t_0+T]$ into subintervals
$$
 t_0 < t_1 < \ldots < t_n = t_0+T.
$$
In our algorithm we construct a sequence of sets  $\mathcal{D}_i$
\begin{align*}
\mathcal{D}_i=&\{  v\in H^4 \cap H^1_0\,:\ u_{xx}\in H^1_0,\\
&\qquad  \|v\|_{L^2}\leq R_1^i, \|v_x\|_{L^2}\leq R_2^i, \|v_{xx}\|_{L^2}\leq R_3^i, \|v_{xxx}\|_{L^2}\leq R_4^i, \|v_{xxxx}\|_{L^2}\leq R_5^i, P_kv\in P_k^i\}, \\
P_k^i &\subset V_k \quad \mbox{is nonempty, closed, convex and bounded}
\end{align*}
such that $u(t_i) \in \mathcal{D}_i$ for $i=0,1,\dots,n$. In each step of the algorithm we propagate $\mathcal{D}_i$ forward in time through the process $S(t_{i+1},t_i)$ governing the strong solutions given by Definition \ref{ref:def_burg_strong}.

Before we present the algorithm, we prove the  lemma, which allows us to improve the estimates on the $L^2$ norms of the function and its first derivative basing on the information that the projection belongs to some given set $S$.
\begin{lemma}\label{lemma:refinement}
	Let $N_m>0$ for $m\in \{ 1,\ldots,5\}$ and let $S\subset V_k$. Let the set $\mathcal{D}$ be given by
\begin{align*}	
	\mathcal{D} =&\{  v\in H^4 \cap H^1_0\,:\ u_{xx}\in H^1_0,\\
	&\qquad  \|v\|_{L^2}\leq N_1, \|v_x\|_{L^2}\leq N_2, \|v_{xx}\|_{L^2}\leq N_3, \|v_{xxx}\|_{L^2}\leq N_4, \|v_{xxxx}\|_{L^2}\leq N_5, P_kv\in S\}.
\end{align*}
Define
\begin{eqnarray}
\overline{N}_1 &=& \min \left\{ N_1, \frac{N_3 h_k^2 }{\pi^2} + \sup_{v\in S}\|v\|_{L^2} \right\}, \label{eq:R1-corr} \\
\overline{N}_2 &=& \min \left\{   N_2, \frac{N_3 h_k}{\pi} + \sup_{v\in S}\|v_x\|_{L^2} \right\}.  \label{eq:R2-corr}
\end{eqnarray}
Then
\begin{align*}	
\mathcal{D} =&\{  v\in H^4 \cap H^1_0\,:\ u_{xx}\in H^1_0,\\
&\qquad  \|v\|_{L^2}\leq \overline{N}_1, \|v_x\|_{L^2}\leq \overline{N}_2, \|v_{xx}\|_{L^2}\leq N_3, \|v_{xxx}\|_{L^2}\leq N_4, \|v_{xxxx}\|_{L^2}\leq N_5, P_kv\in S\}.
\end{align*}
\end{lemma}
\begin{proof}
	Denote
	\begin{align*}	
	\mathcal{D}_1 =&\{  v\in H^4 \cap H^1_0\,:\ u_{xx}\in H^1_0,\\
	&\qquad  \|v\|_{L^2}\leq \overline{N}_1, \|v_x\|_{L^2}\leq \overline{N}_2, \|v_{xx}\|_{L^2}\leq N_3, \|v_{xxx}\|_{L^2}\leq N_4, \|v_{xxxx}\|_{L^2}\leq N_5, P_kv\in S\}.
	\end{align*}
	It is clear that $\mathcal{D}_1 \subset D$. Let $w\in \mathcal{D}$. Then from \eqref{eq:estimate_1} and \eqref{eq:estimate_2} we deduce that
	\begin{equation}
	\|Q_k w\|_{L^2} \leq \frac{h_k^2}{\pi^2} N_3,\quad \|(Q_k w)_x\|_{L^2} \leq \frac{h_k}{\pi} N_3,\label{eq:locbd_3}
	\end{equation}
	hence
	\begin{align}
	& \|w\|_{L^2} \leq \|P_k w\|_{L^2} + \|Q_k w\|_{L^2} \leq \sup_{v\in S}\|v\|_{L^2} + \frac{h_k^2}{\pi^2} N_3,\label{314}\\
	& \|w_x\|_{L^2} \leq \|(P_k w)_x\|_{L^2} + \|(Q_k w)_x\|_{L^2} \leq \sup_{v\in S}\|v_x\|_{L^2} + \frac{h_k}{\pi} N_3,\label{315}
	\end{align}
which proves that $w\in \mathcal{D}_1$ and the proof is complete.
\end{proof}
We pass to the presentation of our algorithm.
\begin{itemize}
	\item[(i)] \textit{Initialization.} Choose $\mathcal{D}_0$, so that it contains our initial condition. Refine the radii $R^0_i$ and $R^0_2$ using Lemma \ref{lemma:refinement}.
	
	\item[(ii)] \textit{For $i=0,\ldots, n-1$ repeat steps (iii)--(v).}
	
	\item[(iii)] \textit{Computation of local integral bounds on interval $[t_i,t_{i+1}]$.} From the local a priori estimates of Lemma \ref{lemma_local} find constants $M^{i+1}_k$ such that for every solution $u$ defined on interval $[t_i,t_{i+1}]$ with the initial data $u(t_i) = u_i$ such that $u_i \in  \mathcal{D}_i$   there holds
	\begin{align}
	&\|u(t)\|_{L^2}  \leq M^{i+1}_1,
	\|u_x(t)\|_{L^2}  \leq M^{i+1}_2,
	\|u_{xx}(t)\|_{L^2}  \leq M^{i+1}_3, \|u_{xxx}(t)\|_{L^2}  \leq M^{i+1}_4,\\
	& \|u_{xxxx}(t)\|_{L^2} \leq M^{i+1}_5,\\
	&
	\|Q_ku(t)\|_{H^1_0}  \leq \frac{M^{i+1}_3 h_k}{\pi} ,
	\|Q_ku(t)\|_{L^2} \leq  \frac{M^{i+1}_3 h_k^2}{\pi^2} \quad \textrm{for}\quad t\in [t_i,t_{i+1}].\label{bounds_rem_int}
	\end{align}
	Use these estimates to find the multifunction $G$ on the interval $[t_i,t_{i+1}]$. The details of the calculation of $G$ is given in Subsection \ref{g}.

	\item[(iv)] \textit{Solving differential inclusion.} Solve rigorously numerically the inclusion \eqref{inclusion}  using the initial data in $P^i_k$ and the multivalued term, $G(t)$, calculated in step (iii). The solution of the inclusion gives the set $P^{i+1}_k \subset V_k$ which contains the projections on $V_k$ of values at $t_{i+1}$ of all trajectories such that $u(t_i) \in \mathcal{D}_i$.

	\item[(v)] \textit{Calculation of bounds at $t_{i+1}$.} From the local a priori estimates of Lemma \ref{lemma_local} find constants ${R}^{i+1}_1$, ${R}^{i+1}_2$, ${R}^{i+1}_3$, ${R}^{i+1}_4$, and  ${R}^{i+1}_5$ such that for every solution  with $u(t_i) \in \mathcal{D}_i$ there holds
	\begin{align*}
	&\|u(t_{i+1})\|_{L^2}  \leq {R}^{i+1}_1,
	\|u_x(t_{i+1})\|_{L^2}  \leq {R}^{i+1}_2,
	\|u_{xx}(t_{i+1})\|_{L^2}  \leq R^{i+1}_3, \\
	& \|u_{xxx}(t_{i+1})\|_{L^2}  \leq R^{i+1}_4, \|u_{xxxx}(t_{i+1})\|_{L^2} \leq R^{i+1}_5.
	\end{align*}
	Refine the radii $R_{1}^{i+1}$ and $R_2^{i+1}$ using Lemma \ref{lemma:refinement}.
   	Obtained radii $R^{i+1}_1$-$R^{i+1}_5$ together with the set $P^{i+1}_k$ from step (iv) define the set $\mathcal{D}_{i+1}$.
\end{itemize}


\begin{remark}
\label{rem:init-comp}
In our proofs of the existence of the periodic orbit we initialize $\mathcal{D}_0$ by taking the  global radii from Theorem~\ref{thm:trapping}, i.e.
  $R^0_m = R_m$ for $m\in \{1,\ldots,5 \}$ and
for $P_k^0$ we take some neighborhood of the numerically found  periodic point. This choice guarantees us that for every $t>t_0$ and every $w\in S(t,t_0)\mathcal{D}_0$ there hold the bounds
	$$
	\|w\|_{L^2}\leq R_1, \|w_x\|_{L^2}\leq R_2, \|w_{xx}\|_{L^2}\leq R_3, \|w_{xxx}\|_{L^2}\leq R_4, \|w_{xxxx}\|_{L^2}\leq {R_5}.
	$$
	In order to prove the periodic orbit existence we need to verify that $S(t_n,t_0)\mathcal{D}_0 \subset \mathcal{D}_0$. As the algorithm is constructed in such a way, that $S(t_n,t_0)\mathcal{D}_0 \subset \mathcal{D}_n$, it follows that $P_kS(t_n,t_0)\mathcal{D}_0\subset P^n_k$ and it suffices only to verify that $P^n_k \subset P^0_k$.
	
An alternative approach to obtain the periodic orbit would be to choose any initial radii  $R^0_m$ for $m\in \{1,\ldots,5 \}$, unrelated with the global radii from Theorem~\ref{thm:trapping}. Then, after the algorithm stops, we would need to verify that $\mathcal{D}_n\subset \mathcal{D}_0$, i.e., both that $P^n_k \subset P^0_k$ and that $R^n_m \leq R^0_m$ for $m \in \{ 1,\ldots,5\}$.

\end{remark}

	\begin{remark}
In principle in the case of the Burgers equation if our set of initial conditions is contained in the trapping set from Theorem~\ref{thm:trapping}  we could have skipped the stage (iii) and use global bounds throughout the entire simulation, i.e.
we will have $R^i_j=M^i_j=R_j$ for $j=1,\dots,5$ and all $i$. In such situation the multivalued term $G$ could be the same for all time steps.
This however will be very inefficient and will require very fine mesh ($k$ large) to obtain good rigorous bounds. In practice, in all examples which we have run the use of local estimates leads to very significant gain.
	\end{remark}

A reader might wonder whether the algorithm can be generalized to any PDE  of the form
\begin{equation}
  u_t = \Delta u + N(u,Du) + f(x,t),
\end{equation}
equipped with some boundary condition.

Let us briefly discuss the problems one can face
\begin{itemize}
\item in stage (iii) it might be impossible to obtain a local integral bounds for solution for a given time step, for example the solution
   might blow up,
\item in  stage (iv) the set found by the computer could expand too much for a given time step.
\end{itemize}
Both stages involve some heuristics - we need to obtain the a priori bounds, which could depend on the particular form of the nonlinearity and the boundary conditions.

 \subsection{Construction of $G(t)$ from local a priori bounds.}\label{g}
Following \eqref{G_included} we need to construct $G:\mathbb{R}\to (2^{\mathbb{R}})^{k-1}$ given by $G(t) = (G_l(t))_{l=1}^{k-1}$  such that
$$
\frac{6}{h_k}\left((Q_ku(t) Q_ku_x(t)-Q_ku_t(t),w^l) -  (P_ku(t) Q_ku(t),w^l_x)\right) \in G_l(t) \quad \textrm{for}\quad l=1,\ldots, k-1.
$$
The multifunction $G(t)=(G_l(t))_{l=1}^{k-1}$ will be constant on every interval $[t_i,t_{i+1})$ and
$$G_l(t) = [-\epsilon_{i+1, l},\epsilon_{i+1, l}]\qquad \textrm{for}\qquad t\in [t_i,t_{i+1}).$$
To find the concrete numerical values $\epsilon_{i+1, l}$ using the bounds \eqref{bounds_rem_int} and Lemma \ref{lem:embedding} we estimate
 \begin{align*}
  & |(Q_ku P_k u,w^l_x)| \leq \|Q_ku\|_{L^2}\|w^l_x\|_{L^2}\|P_k u\|_{L^\infty}\leq  \frac{M_3^{i+1}h_k^2}{\pi^2}\frac{\|(P_k u)_x\|_{L^2}}{2}\|w^l_x\|_{L^2}\\
  &\qquad  \leq \|u_x\|_{L^2} \frac{M_3^{i+1}h_k^{2}}{2\pi^{2}}\|w^l_x\|_{L^2} \leq  \frac{M_2^{i+1} M_3^{i+1}h_k^{2}}{2\pi^{2}}\|w^l_x\|_{L^2}.
 \end{align*}
 Now
 \begin{align*}
&
| (Q_ku(Q_ku)_x,w^l)| \leq \|Q_ku\|_{L^\infty}\|(Q_ku)_x\|_{L^2}\|w^l\|_{L^2} \\
& \qquad \leq \|Q_ku\|_{L^2}^{1/2}\|(Q_ku)_x\|_{L^2}^{3/2}\|w^l\|_{L^2} \leq \frac{(M_3^{i+1})^2 h_k^{5/2}}{\pi^{5/2}}\|w^l\|_{L^2}.
 \end{align*}
Moreover
\begin{align*}
& |(Q_ku_t,w^l)| \leq |(Q_k u_{xx},w^l)| + |(Q_k(uu_x),w^l)|+ |(Q_kf,w^l)| \\
& \qquad \leq (\|Q_k(uu_x)\|_{L^2}+\|Q_k u_{xx}\|_{L^2})\|w^l\|_{L^2} + \sup_{t\in [t_i,t_{i+1}]}|(Q_kf(t),w^l)|\\
&\qquad \leq  \frac{h_k^2}{\pi^2} (\|(uu_x)_{xx}\|_{L^2}+\| u_{xxxx}\|_{L^2})\|w^l\|_{L^2} + \sup_{t\in [t_i,t_{i+1}]}|(Q_kf(t),w^l)|.
\end{align*}
On the other hand
\begin{align*}
&  \|u_{xxxx}\|_{L^2} +  \|(uu_x)_{xx}\|_{L^2} \\
&\qquad \leq \|u_{xxxx}\|_{L^2}  + 3\|u_x\|_{L^\infty} \|u_{xx}\|_{L^2} + \|u\|_{L^\infty} \|u_{xxx}\|_{L^2}\\
& \qquad \leq \|u_{xxxx}\|_{L^2}  + 3\sqrt{2}\|u_x\|_{L^2}^{1/2} \|u_{xx}\|_{L^2}^{3/2} + \|u\|_{L^2}^{1/2}\|u_x\|_{L^2}^{1/2} \|u_{xxx}\|_{L^2}.
\end{align*}
So on interval $[t_i,t_{i+1}]$ there holds
\begin{align*}
& |(Q_ku_t,w^l)| \\
&\  \leq \frac{h_k^2}{\pi^2}\|w^l\|_{L^2}\left( M_5^{i+1}  + 3\sqrt{2}\left(M_2^{i+1}\right)^{1/2} \left(M_3^{i+1}\right)^{3/2} + \left(M_1^{i+1}\right)^{1/2}\left(M_2^{i+1}\right)^{1/2} M_4^{i+1}\right)\\
& \qquad \qquad  + \sup_{t\in [t_i,t_{i+1}]}|(Q_kf(t),w^l)|\\
& \ \ =: \frac{h_k^2}{\pi^2}C_{i+1}\|w^l\|_{L^2} + \sup_{t\in [t_i,t_{i+1}]}|(Q_kf(t),w^l)|.
\end{align*}
Summarizing the above three estimates, there holds
\begin{equation}\label{epsilon}
\epsilon_{i,l} = \frac{6h_k}{\pi^2}\left(\frac{M_2^{i+1} M_3^{i+1}}{2}\|w^l_x\|_{L^2} + \frac{(M_3^{i+1})^2 h_k^{1/2}}{\pi^{1/2}}\|w^l\|_{L^2} + C_{i+1}\|w^l\|_{L^2} \right) + \frac{6}{h_k}\sup_{t\in [t_i,t_{i+1}]}|(Q_kf(t),w^l)|.
\end{equation}
These values can be calculated effectively in step (iii) of the algorithm leading to the width of the differential inclusion in every time step.

\begin{remark}\label{rem:45}
In the above derivation we use $Q_kf$ and hence we
use the fact that $f(t) \in H^1_0$. Without this assumption one could proceed with the following bound 
\begin{align*}
& |(Q_ku_t,w^l)| = |(Q_k (uu_x+f-u_{xx}),w^l)| \\
& \qquad  \leq |(Q_k(uu_x),w^l)| + \frac{h_k^2}{\pi^2}(\|(uu_x)_{xx}\|_{L^2}+\|f_{xx}\|_{L^2}+\|u_{xxxx}\|_{L^2})\|w_l\|_{L^2}.
\end{align*}
This bound is less optimal in terms of expression with $f$, but due to presence of $h_k^2$ also leads to the admissible estimates on the remainder.
\end{remark}

The desired feature of $\epsilon_{i,l}$ should be that it decreases to zero with the decrease of $h_k$. The change of basis matrix $B$ used in the derivation of \eqref{inclusion} is not uniquely defined, hence  neither are  vectors $w^l$. In the numerical examples of Section \ref{sec:alg-proof} we have chosen $B = (b_{ij})_{i,j=1}^{k-1}$ so that its columns
are normalized as follows
\begin{equation}\label{normalization}
  \sum_{i=1}^{k-1} b_{ij}^2 =1,
\end{equation}
i.e. the euclidean norm of each column of $B$ is normalized to one.

We perform a brief analysis of the behavior of the width $\epsilon_{i,l}$ of the inclusion as $h_k$ decreases. As the  functions $w^l$ are approximately equal to the eigenfunctions of the Laplace operator, after reordering of the basis
$w^l$ so that $|\lambda_l|$ is increasing, we obtain
\begin{eqnarray*}
  \lambda_l &\approx& -\pi^2 l^2, \\
  w^l(x) &=& \sum_{m=1}^{k-1} c_{lm}v^m(x) \approx \pm \alpha(k,l)\sin(l \pi x).
\end{eqnarray*}
The normalization choice \eqref{normalization} makes the unknowns $(\beta_l)_{l=1}^{k-1}$ scale with the increase of $k$. In order to avoid this, we perform the analysis for such choice of $B$ that the $L^2$ norm of the basis vectors $w^l$ are constant. Then, approximately,
\begin{eqnarray*}
		\|w^l\|_{L_2} &\approx& 1, \\
		\|w_x^l\|_{L_2} &\approx& l\pi.
\end{eqnarray*}
The formula \eqref{epsilon} for $\epsilon_{k,l}$ contains two ingredients, the first one contains the factor $h_k$, while the other one
$\frac{1}{h_k}$. In the first ingredient given by
$$ \frac{6h_k}{\pi^2}\left(\frac{M_2^{i+1} M_3^{i+1}}{2}\|w^l_x\|_{L^2} + \frac{(M_3^{i+1})^2 h_k^{1/2}}{\pi^{1/2}}\|w^l\|_{L^2} + C_{i+1}\|w^l\|_{L^2} \right)$$
the values $M_j^i$ are the local a priori estimates on the sought trajectories which are bounded from above by constants which are independent of the size of the mesh.  In such situation
	$$ \frac{6h_k}{\pi^2}\left(\frac{M_2^{i+1} M_3^{i+1}}{2}\|w^l_x\|_{L^2} + \frac{(M_3^{i+1})^2 h_k^{1/2}}{\pi^{1/2}}\|w^l\|_{L^2} + C_{i+1}\|w^l\|_{L^2} \right) = lO(h_k) + O(h_k^{3/2}) + O(h_k).$$
 The expression is dominated by the first term which is equal to $lO(h_k)$, remaining two terms are both at least $O(h_k)$. For large $l$, and hence for coordinates corresponding to the high modes of the solution, this term will be $O(1)$. This unwelcome effect will be fixed by the fact that the diagonal entries of the linear term in the equations corresponding to the variables with large $l$ will be large negative numbers exploited in the method of dissipative modes described in Section \ref{diss}.

The second ingredient
$$\frac{6}{h_k}\sup_{t\in [t_i,t_{i+1}]}|(Q_kf(t),w^l)|$$
appears to blow-up with increasing $k$, because $1/h_k \to \infty$ with $k\to \infty$. However, from \eqref{eq:estimate_1} and \eqref{eq:estimate_2} it follows that
\begin{equation}
  |(Q_kf(t),w^l)| \leq \|Q_kf\|_{L_2} \|w^l\|_{L_2} \leq \frac{h_k^2 }{\pi^2}  \|f_{xx}\|_{L_2}\|w^l\|_{L_2}.
\end{equation}
Hence the effective size of this ingredient is equal to $O(h_k)$. 

\section{Rigorous integration of differential inclusion}\label{sec:4}

In this section we present the details concerning the technique for the rigorous numerical solution  of the inclusion \eqref{inclusion}, i.e.
the step (iii) from Section~\ref{sec:gal-tail}. 

\subsection{Solving the differential inclusion.}\label{sec:inclusion_rig}
We use the method of \cites{kapela,kuramoto_III} to rigorously solve the inclusion
\begin{equation}\label{incl}
\beta'(t) \in B^{-1}M^{-1}S B\beta(t) + B^{-1}M^{-1}N(B\beta(t)) + B^{-1}M^{-1}F(t)+ G (t).
\end{equation}
We will use the notation
\begin{equation}
f(\beta) = B^{-1}M^{-1}S B\beta + B^{-1}M^{-1}N(B\beta),\ \  h(t) = B^{-1}M^{-1}F(t).
\end{equation}
 We solve the above inclusion on the time interval $[t_i,t_{i+1}]$, and we equip it with the initial data $\beta(t_i) \in P^i_k$. By the rigorous solution we understand finding the set $P^{i+1}_k$ such that for every absolutely continuous function $\beta:[t_i,t_{i+1}]\to \mathbb{R}^{k-1}$ satisfying \eqref{incl} for almost every $t\in (t_i,t_{i+1})$ with $\beta(t_0)\in P^i_k$ there holds $\beta(t_{i+1})\in P^{i+1}_k$. Note that on $[t_i,t_{i+1})$ the set $G(t)$ is independent of time and equal to $G(t) = (G_l(t))_{l=1}^{k-1}$ with $G_l(t) = [-\epsilon_{i,l},\epsilon_{i,l}]$ where $\epsilon_{i,l}$ is given by \eqref{epsilon} and can be effectively calculated.

The set $P^{i+1}_k$ is found using the  algorithm of \cite[Lemma 5.2]{kapela} which is a part of
the rigorous numerics CAPD library \cite{CAPD}. We briefly recall the algorithm here. As the multifunction $G(t)$ is always centered at zero, in the first step, an equation
\begin{equation}\label{ode}
\beta'(t) = f(\beta(t)) + h(t),
\end{equation}
is solved with the initial data $\beta(t_i)\in P_k^i$. The rigorous numerical solution of this ODE uses the explicit Taylor scheme \cite{CAPD} (in all examples we use the fourth order scheme) and the Lohner algorithm \cite{reference_lohner} for representation of sets $P_k^i$ as parallelepipeds and their propagation in time. For a parallelepiped $P_k^i$, another parallelepiped is found which is guaranteed to contain the values at $t_{i+1}$ of all solutions of the above ODE with the initial data taken at $t_i$ in the set $P_k^i$.  Next, a correction is calculated and added to resultant set to guarantee to contain all solutions of the inclusion. This correction \cite[Sec. 6.3]{kapela} is equal to $(-d_j,d_j)_{j=1}^{k-1}$, where
$$
d_j = \left|\left(\int_{t_i}^{t_{i+1}}e^{J(t_i-s)}C\, ds\right)_j\right|,
$$
where the vector $C$ is given by
$$
C = (\epsilon_{i,l})_{l=1}^{k-1}
$$
and the matrix $J$ is given by
$$
J_{ij} = \begin{cases} \sup_{\beta \in [W]} \frac{\partial f_i(\beta)}{\partial \beta_j}\quad \textrm{if}\quad i=j,\\
\sup_{\beta \in [W]} \left|\frac{\partial f_i(\beta)}{\partial \beta_j}\right|\quad \textrm{otherwise},
\end{cases}
$$
where the set $[W]\subset\mathbb{R}^{k-1}$ is the so called enclosure, i.e. the set which is guaranteed to contain the values at all $t\in [t_i,t_{i+1}]$ of all solutions of \eqref{incl} with the initial data $\beta(t_i) \in P_k^i$. The enclosure is found using the first order rough enclosure algorithm \cites{enclosiure_ref, kuramoto_III}.

\subsection{Method of dissipative modes.} \label{diss}
Denoting the linear term in \eqref{incl} by $A = B^{-1}M^{-1}SB$, we observe that is given by the interval matrix which is close to the diagonal one, and its diagonal entries, sorted by the increasing moduli of the eigenvalues, are approximately given by $[A_{ll}] \approx -\pi^2l^2$ for $l\in \{ 1,\ldots, k-1\}$. The algorithm  of dissipative modes designed and described in \cite{kuramoto_III} allows us to use to our advantage the fact that these diagonal entries are large and negative for large $l$.
The method allows us to \textit{rigorously solve the $k-1$ dimensional differential inclusion}  by splitting all unknowns $\beta = (\beta_l)_{l=1}^{k-1}$ into two groups $\beta=(\beta_{1},\beta_{2})$.  The variables $\beta_1$ correspond to those entries $[A_{ll}]$ which have the  moduli smaller than some arbitrary cut-off value (in numerical examples we took the $8$ first coordinates of $\beta$ as $\beta_1$) while the remaining unknowns are assigned to $\beta_2$. The rigorous numerical integration which uses the Taylor scheme and Lohner algorithm is performed \textit{only for variables in $\beta_1$} while the \textit{dissipative modes algorithm allows us to treat the variables belonging to $\beta_2$}. Each of the inclusions for these variables  can be written as
\begin{equation}\label{decomp}
\beta'_{2,l}(t) + |[A_{ll}]| \beta_{2,l}(t) = f_l(\beta_1(t),\beta_2(t)) + |[A_{ll}]|\beta_{2,l}(t) + h_l(t) + G_l(t),
\end{equation}
where $\beta_{2,l}$ is the $l$-th coordinate of $\beta$ and the index $2$ has been added to indicate that this is the variable belonging to $\beta_2$. In the first step, an enclosure is found, that is the set $[W]\subset \mathbb{R}^k$ such that $\beta(t)\in [W]$ for $t\in [t_i,t_{i+1}]$. In other words all solutions of the inclusion \eqref{incl} are guaranteed to belong to this set for $t\in [t_i,t_{i+1}]$ if the initial data satisfies $\beta(t_i)\in P_k^i$. We note that in the algorithm of finding this enclosure the fact that $|[A_{ll}]|$ in equation \eqref{decomp} is large plays a crucial role, see \cite[Section 5]{kuramoto_III} for the details. Once the enclosure is found, it is possible to find the numbers $N_{l}^-$ and $N_l^+$ which are bounds from above and from below on the sets $f_l(\beta_1(t),\beta_2(t)) + [\lambda_l] \beta_{2l}(t) + h_l(t) + G_l(t)$ for $t\in [t_i,t_{i+1}]$, i.e. for every solution $\beta$ with the initial data $\beta(t_i)\in P_k^i$, every $t\in [t_{i},t_{i+1}]$, and every $y\in G_l(t)$ there holds
$$
N_l^- \leq f_l(\beta_1(t),\beta_2(t)) + |[A_{ll}]|{ \beta_{2l}(t) + h_l(t) + y \leq N_l^+.}
$$
Thus all solutions $\beta_{2,l}(t)$ satisfy the differential inequalities
$$
N_{l}^- \leq \beta'_{2,l}(t) + |[A_{ll}]| \beta_{2,l}(t) \leq N_l^+.
$$
A simple computation leads us to
\begin{align}
& (\beta_{2,l}(t_i)^-)e^{-|[A_{ll}]|(t_{i+1}-t_i)}+\frac{N_l^-}{|[A_{ll}]|}\left(1-e^{-|[A_{ll}]|(t_{i+1}-t_i)}\right)\nonumber\\
& \qquad  \leq  \beta_{2,l}(t_{i+1}) \leq (\beta_{2,l}(t_i)^+)e^{-|[A_{ll}]|(t_{i+1}-t_i)}+\frac{N_l^+}{|[A_{ll}]|}\left(1-e^{-|[A_{ll}]|(t_{i+1}-t_i)}\right),\label{incl:diss}
\end{align}
whence we obtain the interval which is guaranteed to contain $\beta_{2,l}(t_{i+1})$, cf. \cite[Theorem 23]{kuramoto_III}.
Independently, the found enclosure on variables $\beta_1$ and $\beta_2$ is used to construct the low dimensional inclusion for $\beta_1$ which is solved by the rigorous integration algorithm described in Section \ref{sec:inclusion_rig}.

Thanks to the use of the dissipative modes approach we have the following advantages.
\begin{itemize}
	\item Since the rigorous integration is performed only for fixed small number of variables in $\beta_1$, the contribution from the multivalued term $G_l(t)$ for these variables, given by $lO(h_k)$, is effectively equal to $O(h_k)$ as $l$ is small.
	
	\item In the variables belonging to $\beta_2$ the contribution from $G_l(t)$ could be equal even to $O(1)$ for $l$ close to the maximal value $k-1$. This large width of $G_l(t)$, incorporated into $N^+_l$ and $N_l^-$ is counteracted  by the division by $|[A_{ll}]|$ in \eqref{incl:diss}, making the contribution of the multivalued term effectively equal to $O(h_k)/l$ in the computation of $\beta_{2,l}(t_{i+1})$.
	
	\item The high computational cost of the Taylor integration algorithm and Lohner algorithm is significantly reduced, as the algorithm needs to be run only for low-dimensional problems.
	
	\item As the Taylor scheme is explicit, the admissible time step length required for its stability is bounded from above. This bound is, approximately, the increasing function of maximum of $1/|[A_{ll}]|$ over indexes $l$ of variables assigned to $\beta_1$. Hence, as only those $l$ for which $|[A_{ll}]|$ is small are assigned to $\beta_1$, we can perform the simulation with larger time steps, leading to the significant reduction of the computation time.
\end{itemize}
The above observations are backed by our numerical examples, the implementation of the approach by the dissipative modes turned out to be an absolutely crucial factor to prevent the uncontrolled expansion of the sets $P^i_k$ obtained during the course of the rigorous numerical simulation.

\section{Computer assisted verification of the periodic solution existence}
\label{sec:alg-proof}
In this section we present the theorem on the existence of the periodic orbit for two particular choices of the forcing term. In both examples we rigorously integrate forward in time the inclusion \eqref{inclusion} for the $1$-periodic in time forcing and for the initial data belonging to some set $\mathcal{D}_0$ given by Step (i) of the algorithm described in Section \ref{subsec:numerical}. After time $1$, the period of $f$, all solutions of the inclusion are verified to belong to the set which is a subset of $\mathcal{D}_0$, guaranteeing that all assumptions of the abstract Schauder type fixed point theorem are satisfied.
\subsection{Schauder-type theorem.} In the result of this subsection we  establish the existence of periodic orbit using the following Schauder type  theorem with $X=H^1_0$.
\begin{theorem}\label{thm:schauder}
	Let $X$ be a Banach space and let $\mathcal{B}\subset X$ be a nonempty, compact, and convex set. If the mapping $S:\mathcal{B} \to \mathcal{B}$ is continuous, then it has a fixed point $u_0\in \mathcal{B}$.
\end{theorem}

	Once the time-periodic solution exists, Theorem \ref{thm:attr} establishes that it exponentially attracts all weak solutions given by Definition \ref{ref:def_burg}. In fact Theorem \ref{thm:attr} already establishes the existence of periodic attracting orbit if the forcing term is periodic. The purpose of the present result is, on one hand, the construction, with some accuracy, of this orbit, and on the other hand the demonstration of the usefulness of the algorithm presented in the previous subsection. The key property is that the set obtained after a period of integration is the subset of the set of initial data, meaning that in the course of the integration the obtained sets shrink with respect to the initial ones.

The following theorem allows us to enclose the periodic trajectory  for $T$-periodic forcing term $f$ .
\begin{theorem}\label{thm:inclusion}
	Let $Y = \{ u\in H^4\cap H^1_0\,:\ u_{xx}\in H^1_0 \}$. Assume that $f\in L^\infty(Y)$ is $T$-periodic. Let the set $\mathcal{D}_0$ be given by
	\begin{align*}
	\mathcal{D}_0=&\{  v\in Y,\\
	&\qquad  \|v\|_{L^2}\leq R_1^0, \|v_x\|_{L^2}\leq R_2^0, \|v_{xx}\|_{L^2}\leq R_3^0, \|v_{xxx}\|_{L^2}\leq R_4^0, \|v_{xxxx}\|_{L^2}\leq R_5^0, P_kv\in P_k^0\}, \\
	P_k^0 &\subset V_k \quad \mbox{is nonempty, closed, convex and bounded}
	\end{align*}
Assume that either of the following two conditions hold
\begin{itemize}
	\item[(i)] The set $\mathcal{D}_n$ obtained in the algorithm given in Section \ref{subsec:numerical} satisfies $\mathcal{D}_n \subset \mathcal{D}_0$.
	\item[(ii)] The radii chosen in the definition of the set $\mathcal{D}_0$ are equal to the radii $R_1$-$R_5$ of the trapping set given by Theorem \ref{thm:trapping} and $P^n_k$ obtained in the algorithm given in Section \ref{subsec:numerical} satisfies $P^n_k \subset P^0_k$.
	\end{itemize}
 Then there exists a periodic and bounded trajectory $u(t)$ such that
	\begin{align*}
	&P_ku(t_i) \in P^i_k\quad \textrm{for every}\quad i=0,\ldots,n,\\
	& 	\|Q_ku(t_i)\|_{L^2} \leq \frac{h_k^2}{\pi^2} R^i_3, \|Q_ku(t_i)\|_{H^1_0} \leq \frac{h_k}{\pi} R^i_3\quad \textrm{for every}\quad i=0,\ldots,n.
	\end{align*}
\end{theorem}
\begin{proof}
	The set $\mathcal{D}_0$  is convex and nonempty. Moreover, as it is closed and bounded in $Y$, it is compact in $H^1_0$. As $S(t_0+T,t_0)\mathcal{D}_0 \subset \mathcal{D}_n$, the mapping $S(t_0+T,t_0)$ leads from this set into itself. Indeed in case (i) this directly follows from the inclusion $\mathcal{D}_n\subset \mathcal{D}_0$ and in case (ii) from  Theorem \ref{thm:trapping}, the fact that $P_kS(t_0+T,t_0)\mathcal{D}_0 \subset P^n_k$ and the inclusion $P^n_k\subset P^0_k$. Hence, by Theorem \ref{thm:schauder} it has a fixed point.  Denoting Denote this fixed point by $u_0$. Define
	$$
	u(t) = S\left(t- \left[\frac{t-t_0}{T}\right]T,t_0\right)u_0,
	$$
	where $[s] = \max\{ n\in \mathbb{Z}\,:\ n\leq s \}$ is the largest integer no greater than $s$. This function satisfies $u(t) \in \mathcal{B}_0$, it is $T$-periodic, and it is a strong solution of the Burgers equation, in agreement with Definition \ref{ref:def_burg_strong}. Moreover it must hold that $P_ku(t_i) \in P^i_k$, and, by Lemma \ref{lem:galerkin_est}
	$$
	\|Q_ku(t_i)\|_{L^2} \leq \frac{h_k^2}{\pi^2} R^i_3, \|Q_ku(t_i)\|_{H^1_0} \leq \frac{h_k}{\pi} R^i_3\quad \textrm{for every}\quad i=0,\ldots,n.
	$$
	The proof is complete.
\end{proof}

We remark that the set $P^0_k$ must be found experimentally, based on the numerical simulations. In practice it is found by classical (nonrigorous) FEM simulations and choosing some ball enclosing the found nonrigorous solution at time $t_0$. In two examples given in the following subsection we use the condition (ii) of the above theorem, which makes it necessary to numerically verify only the inclusion $P^n_k \subset P^0_k$.

\subsection{Examples}
 Our numerical experiments show that the key factor which decides if the sets $P_k^i$ shrink or expand in time, is the  width of the multivalued term $G(t)$.

 Our aim was to perform the computer-assisted  proofs in the situation when $G(t)$ is large and hence we chose two particular forms of the forcing term
\begin{align*}
& f(x,t) = 8(\sin(3\pi x)+\sin(4 \pi x))(1+\sin(2\pi t)),\\
& f(x,t) = 12 \sin(\pi x)\sin(2 \pi t).
\end{align*}
Note that the width of $G(t)$ is determined by the values $\epsilon_{i,l}$ given by \eqref{epsilon} which are monotone increasing functions of the constants $M_k^{i+1}$, for $k=1,\ldots,5$, given by Lemma  \ref{lemma_local}. These constants are in turn determined by the estimates of Appendix \ref{trapping_sets}, given in Lemmas \ref{lem:b1}, \ref{lemma:young_x}, \ref{lemma:young_x_2}, \ref{lemma:uxxauxiliary}, \ref{lemma:young_x_3}, \ref{lem_3.18}, \ref{lemma:young_x_4}, \ref{lem:u4}, and \ref{lemma:young_x_5}. The estimates of these Lemmas depend monotonically on the $L^2$ norms of $f(t), f_x(t), f_{xx}(t),$ and $f_{xxx}(t)$.  We chose for our construction the functions $f$ with amplitude $8$ and $12$, and, in the second example, frequencies $3$ and $4$, to show that we can cope with the situation where these norms are relatively large.

 \begin{figure}
	\includegraphics[width=0.8\textwidth]{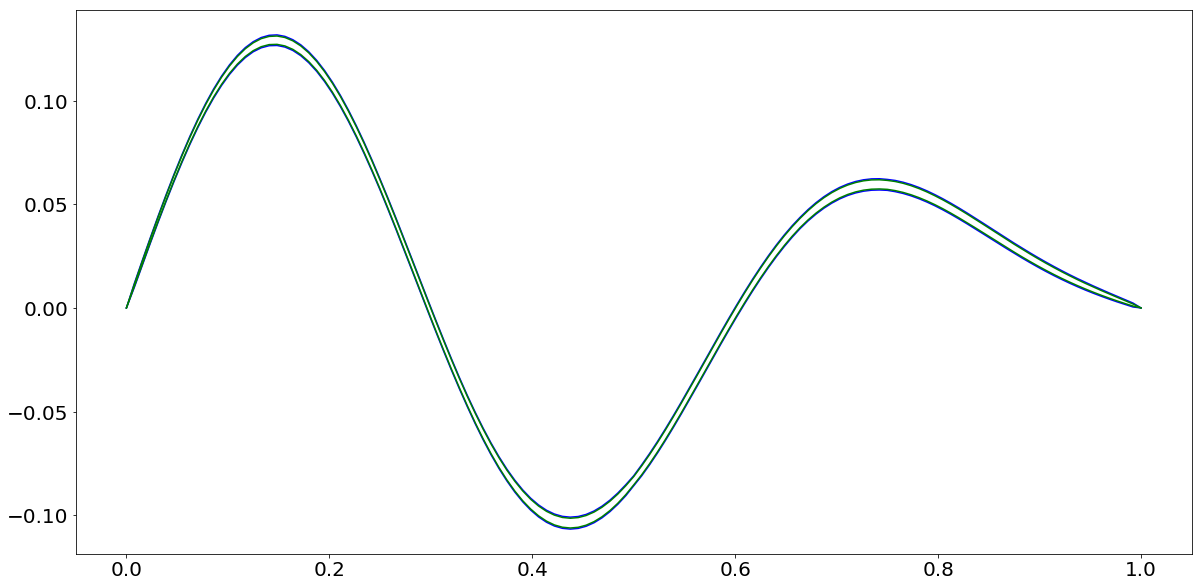}
	\caption{The initial data and simulation result after time $t=1$  for the forcing term
		\eqref{eq:example1} in the FEM basis. Initial set $P^0_k$ is enclosed between two blue curves presented in the plot. Two green curves (which almost coincide with blue ones) enclose the set $P^n_k$ obtained after  $t=1$. Green curves are contained inside the area enclosed by blue curves, and hence the set $P^n_k$ is the subset of $P^0_k$.} \label{fig:1}
\end{figure}

\begin{figure}
	\includegraphics[width=0.8\textwidth]{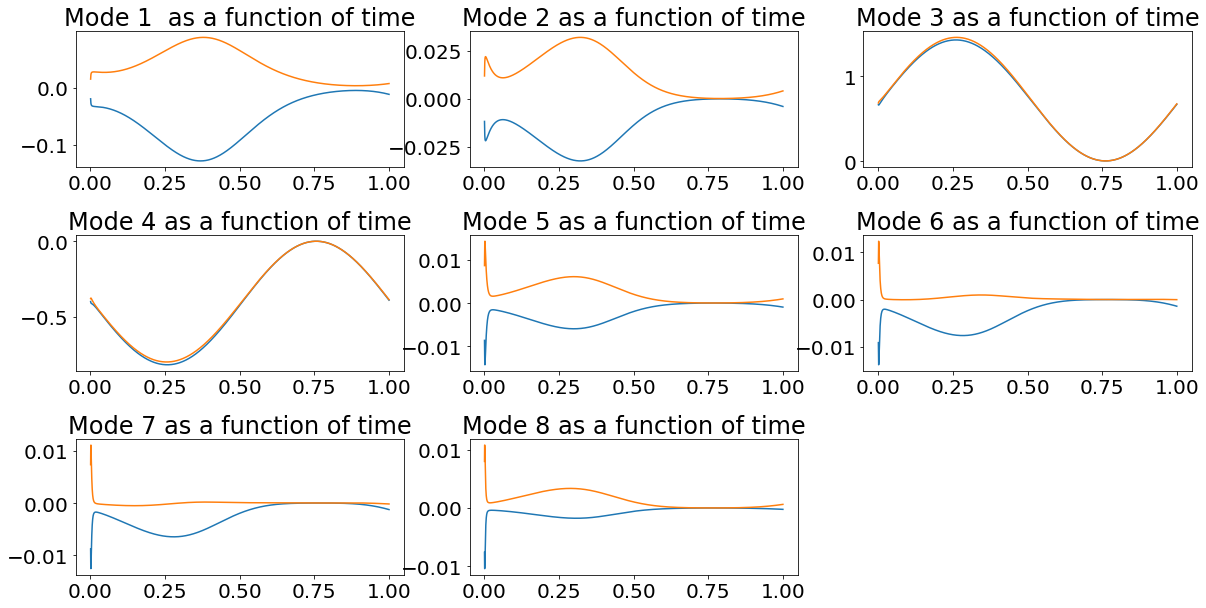}
	\caption{$8$ variables of the solution $\beta$ with the eigenvalues having the smallest modulus for the forcing \eqref{eq:example1} as functions of time. Orange line depicts the upper bound of the interval and the blue one - the lower bound.} \label{fig:2}
\end{figure}

\begin{figure}
	\includegraphics[width=0.8\textwidth]{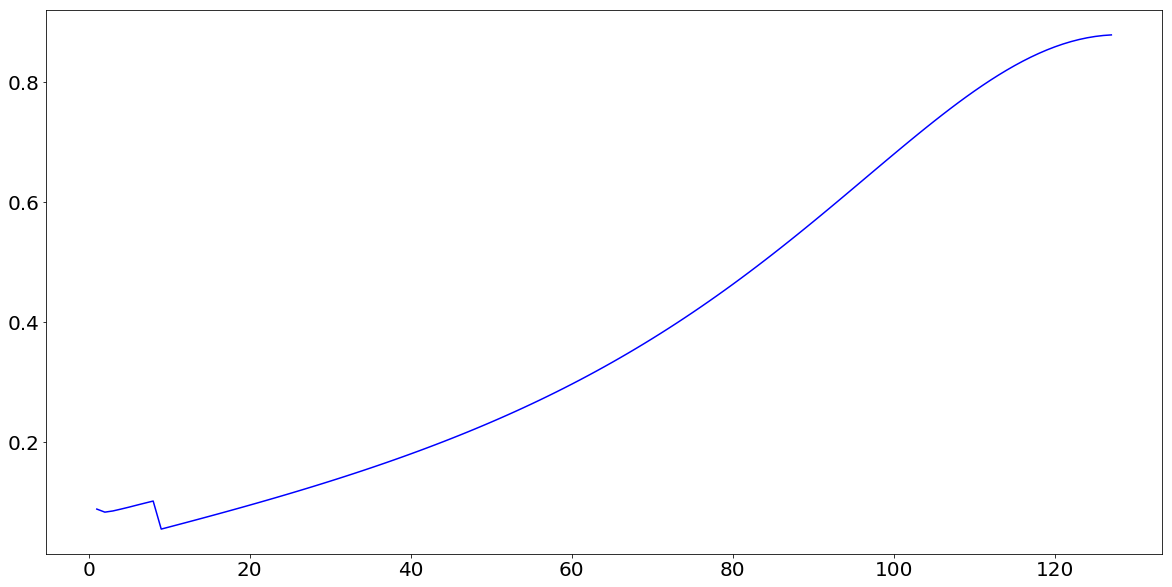}
	\caption{Perturbation for the forcing \eqref{eq:example1} as a function of the variable index ordered by descending eigenvalues for time $t=0.5$.} \label{fig:3}
\end{figure}

\subsubsection{Example 1.} We have validated the precise form of the periodic trajectory for the problem \eqref{a}--\eqref{b} with the forcing term
\begin{equation}\label{eq:example1}
f(x,t) = 8  (\sin(3\pi x) + \sin(4\pi x))(1+\sin(2\pi t)),
\end{equation}
where the space domain is equal to $(0,1)$. The simulation was performed for the mesh interval length size equal to $h_k=1/128$, i.e. the interval $(0,1)$ was divided into $128$ equidistant intervals of length $h_k$. This yields $127$ variables, which correspond to the inner mesh points in the space domain. The finite element basis has been diagonalized in such a way that the euclidean norm of each vector of the diagonal basis is equal to one. After diagonalization $8$ variables were treated as nondissipative ones, and the remaining $119$ variables were treated as dissipative. The constant time step was chosen to be equal to $1/2048$. Simulation took around $55$ minutes. First, we have run the standard FEM in order to identify the candidate for the periodic solution, and then, experimentally we found the set of initial data $P^0_k$ containing the found solution at time $t=0$. The simulation showed that the set obtained after time $1$ was the subset of $P^0_k$ and hence the periodic solution was found.  The radii $R_1$-$R_5$ are given in Tab.\ref{tab:1}.
\begin{table}
\begin{tabular}{|ccc|}
\hline
radius & value & optimal method \\
\hline
\hline
$R_1$ & 2.29264 & Lemma \ref{lem:b1}\\
$R_2$ & 13.9504 & Lemma \ref{lem:abssetL2} improved with Corollary \ref{cor:29} \\
$R_3$ & 135.816 & Lemma \ref{lem:uxx_first}, inequality \eqref{ineq_b16} with interpolation\\
$R_4$ & 1946.47& Lemma \ref{lem:u3x-bnd} with inequality \eqref{eq:u3_root}\\
$R_5$ & 130542& Lemma \ref{lem:u4x-bnd_2} with interpolation\\
\hline
\end{tabular}	
\caption{Optimal radii of trapping sets for the forcing term
	\eqref{eq:example1}. For each radius several methods to compute it were implemented and optimal radius was chosen. }\label{tab:1}
\end{table}

First $8$ coordinates of the found periodic solution in the orthogonalized basis are given in Tab.\ref{tab:2}. The corresponding coordinate of the found periodic solution at $t=0$ must belong to the found interval. Modes $3$ and $4$ are largest as they are the only two modes contained in the forcing term. It is visible that also modes $6$ and $7$ are nonzero because of the presence of the nonlinear term which induces the energy transport from modes $3$ and $4$ to those modes.
\begin{table}
	\begin{tabular}{|cc|}
		\hline
		coordinate & interval  \\
		\hline
		\hline
		$\beta_1$ & [-0.011479,0.00736611]  \\
		$\beta_2$ & [-0.00405501,0.00404116]   \\
		$\beta_3$ & [0.667299,0.672336] \\
		$\beta_4$ & [-0.39088,-0.387498]\\
		$\beta_5$ & [-0.000946911,0.000953705]\\
		$\beta_6$ & [-0.00140807,-0.0000232083]\\
		$\beta_7$ & [-0.00125524,-0.000184778]\\
		$\beta_8$ & [-0.000243156,0.000611533]\\
		\hline
	\end{tabular}	
	\caption{Values of first $8$ modes of found periodic solution at time $t=1$ for the forcing term 	\eqref{eq:example1}. For all modes higher than $8$ zero always belongs to the projection of the set $P_k^i$ onto the given mode.}\label{tab:2}
\end{table}

 Fig. \ref{fig:1} presents the initial datum and set $P^n_k$ after time $1$. The set obtained for $t=1$ is the subset of the initial datum whence all conditions of Theorem \ref{thm:inclusion} are satisfied and we obtain the sought periodic solution.

 Fig. \ref{fig:2} presents the intervals containing the first $8$ coordinates of the evolved  set as functions of time. All plots present lower (blue) and upper (orange) bound of corresponding intervals. Apart from variables $3$ and $4$ which correspond to modes of the forcing term also variables $6$ and $7$ are nonzero during the time of simulation.

 Fig. \ref{fig:3} presents the perturbation (width of the inclusion) as a function of the variable index in the diagonalized basis during one time step, for $t=0.5$. For the first $8$ nondissipative variables the perturbation contains the two components: the one coming from the estimates of the infinite dimensional remainder, and (typically much smaller)   contribution from the modes which are treated as dissipative ones.  Growth of the perturbation with increasing index of variable $l$  is the effect of the terms $\|w^l\|_{L^2}$ and $\|w^l_x\|_{L^2}$ in \eqref{epsilon}.

\begin{figure}
	\includegraphics[width=0.8\textwidth]{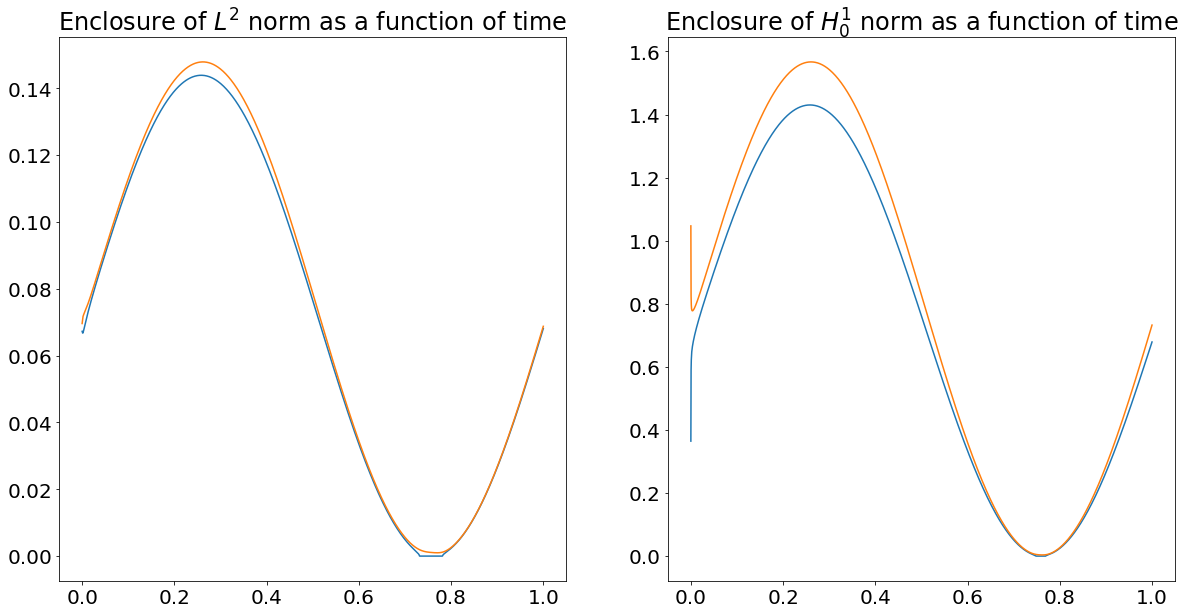}
	\caption{$L^2$ and $H^1_0$ norm enclosures for found periodic solution as the functions of time for the simulation with the forcing \eqref{eq:example1}.} \label{fig:4}
\end{figure}

Finally, Fig.\ref{fig:4} presents the intervals which are guaranteed to contain $L^2$ and $H^1_0$ norm of the solution as the function of time. The sets are constructed by algebraic adding of the sets $P^i_k$ obtained during the simulation (taking into account both dissipative and leading modes which contribute into $P_k^i$) and local $L^2$ and $H^1_0$ estimates for the infinite dimensional remainder $Q_k$.  We formulate the result of this simulation as the following theorem.

\begin{theorem}
	Consider the unique eternal bounded periodic solution $u(t)$ of problem \eqref{a}--\eqref{b} with the forcing term \eqref{eq:example1}. Let $k=128$ and let $P_k u(t)$ be the $H^1_0$ orthoprojection on the space $V_k$. In the diagonalized basis of $V_k$, the $8$ coordinates of $P_ku(0)$ corresponding to the highest (least negative) $8$ eigenvalues, belong to the intervals given in the Tab. \ref{tab:2}. The $L^2$ and $H^1_0$ norm of the periodic solution belongs to the intervals depicted in Fig. \ref{fig:4}.
\end{theorem}

\subsubsection{Example 2.} In the second simulation we have validated the precise form of the periodic trajectory for the problem \eqref{a}--\eqref{b} with the forcing
\begin{equation}\label{eq:example2}
f(x,t) = 12 \sin(\pi x) \sin (2\pi t).
\end{equation}

\begin{figure}
	\includegraphics[width=0.8\textwidth]{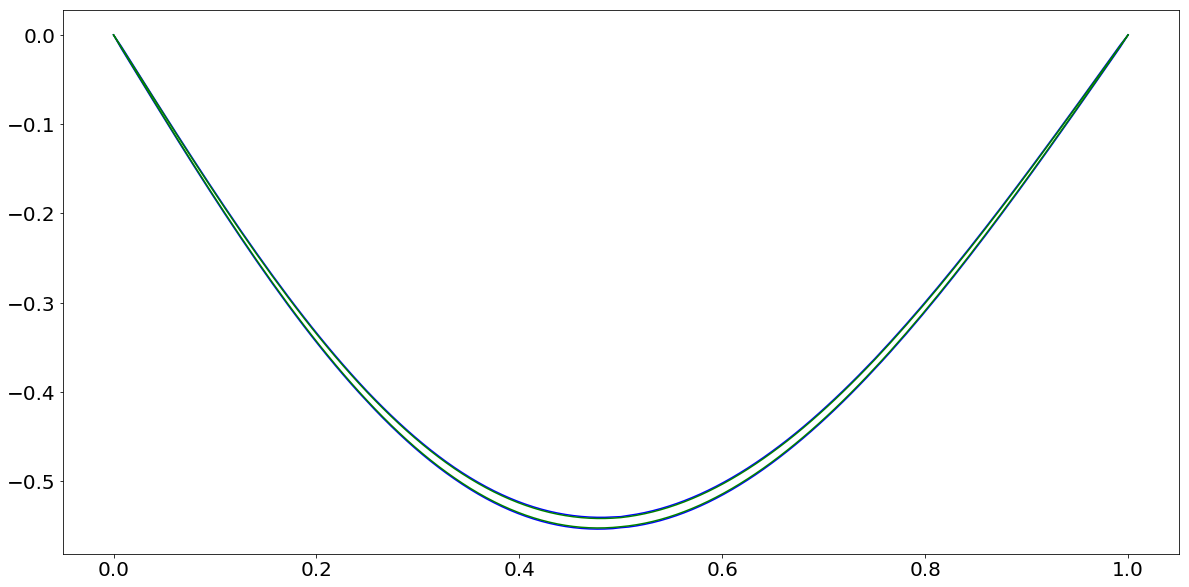}
	\caption{The initial data and simulation result after time $t=1$  for the forcing term
		\eqref{eq:example2} in the FEM basis. Initial data $P_k^0$ is enclosed between two blue curves presented in the plot. Two green curves (which almost coincide with blue ones) enclose the set $P^n_k$ obtained after  $t=1$. Green curves are contained inside the area enclosed by blue curves, and hence the set $P^n_k$ is the subset of $P^0_k$.} \label{fig:5}
\end{figure}

\begin{figure}
	\includegraphics[width=0.8\textwidth]{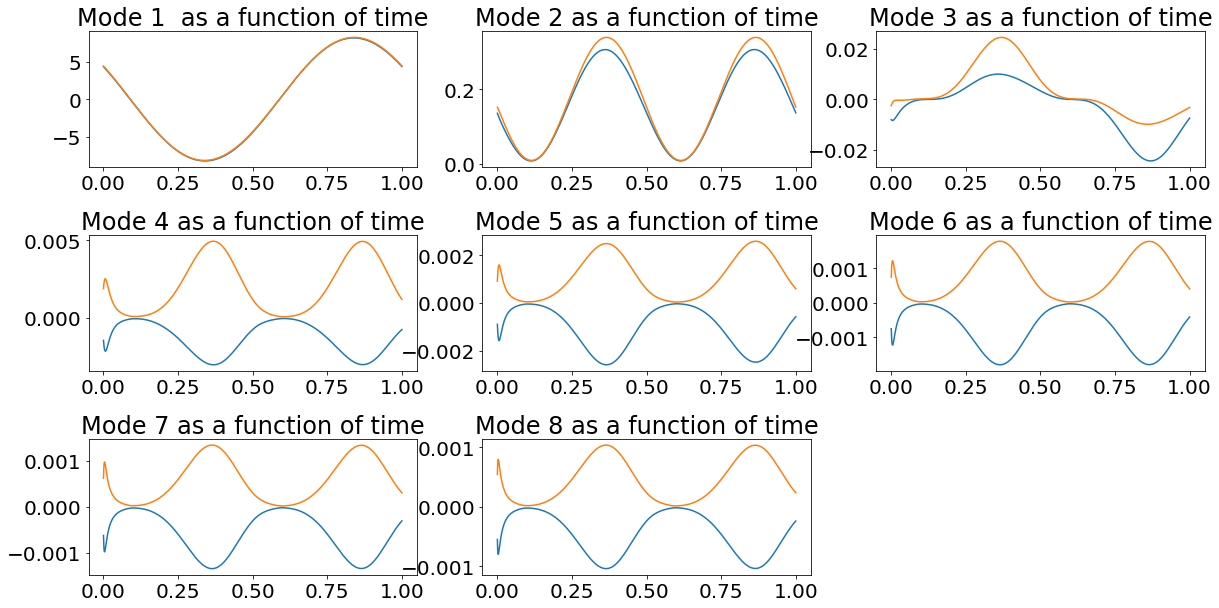}
	\caption{$8$ variables of the solution $\beta$ with the eigenvalues having the smallest modulus for the forcing \eqref{eq:example2} as functions of time. Orange line depicts the upper bound of the interval and the blue one - the lower bound.} \label{fig:6}
\end{figure}

Similar as in the first test, we divided the time interval $(0,1)$ into the equidistant $2048$ time steps, and the space interval $(0,1)$ into $128$ intervals. We treated $8$ modes corresponding to the highest (least negative) eigenvalues as the nondissipative ones. The computation time, as in the first example, was around 55 minutes. Intervals containing the values of first $8$ modes of the found periodic solution are presented in Tab. \ref{tab:3}. Clearly, the first mode has the  highest amplitude, as the energy is inserted into the problem by the forcing term on this mode, but also modes $2$ and $3$ are nonzero indicating the occurrence of energy transport from the first mode to the higher modes via the nonlinearity. The set of initial data as well as the set obtained for $t=1$ in the FEM basis  are presented in Fig. \ref{fig:5}. The plots of intervals containing the first $8$ modes of the found periodic solution as functions of time are presented in Fig. \ref{fig:6}. The result can be formulated as follows.

\begin{theorem}
	Consider the unique eternal bounded periodic solution $u(t)$ of problem \eqref{a}--\eqref{b} with the forcing term \eqref{eq:example2}. Let $k=128$ and let $P_k u(t)$ be the $H^1_0$ orthoprojection on the space $V_k$. In the diagonalized basis of $V_k$, the $8$ coordinates of $P_ku(0)$ corresponding to the highest (least negative) $8$ eigenvalues, belong to the intervals given in the Tab. \ref{tab:3}.
\end{theorem}

\begin{table}
	\begin{tabular}{|cc|}
		\hline
		coordinate & interval  \\
		\hline
		\hline
		$\beta_1$ & [4.33795,4.41266]  \\
		$\beta_2$ & [0.136948,0.152245]   \\
		$\beta_3$ & [-0.00751169,-0.00328565] \\
		$\beta_4$ & [-0.000773302,0.00117999]\\
		$\beta_5$ & [-0.000581862,0.00059719]\\
		$\beta_6$ & [-0.000409425,0.000408848]\\
		$\beta_7$ & [-0.000304782,0.00030476]\\
		$\beta_8$ & [-0.000238547,0.000238549]\\
		\hline
	\end{tabular}	
	\caption{Values of nondissipative $8$ modes of found periodic solution at time $t=1$ for the forcing term 	\eqref{eq:example2}. }\label{tab:3}
\end{table}

\section*{Acknowledgement}The authors wish to thank anonumous referees for insighful remarks and careful reading of our manuscript. They also thank D. Wilczak and T. Kapela for stimulating discussions and great help concerning the C++ programming in CAPD. Work was supported by  the National Science Center (NCN) of the Republic of Poland by the project no UMO-2016/22/A/ST1/00077. Work of PK has also been partially supported by NCN of the Republic of Poland by the grant no DEC-2017/25/B/ST1/00302.

\appendix
\section[Appendix A.]{Auxiliary inequalities and lemmas} \label{sec:appendix_A}


\subsection{Some basic inequalities.} We first recall several basic inequalities which will be used in the paper. We stress that thay are well known and we provide the proofs (with optimal constants) only for the exposition completeness. First of all, we will frequently use the following inequality
\begin{lemma}[Young inequality with $\epsilon$]\label{lem:young}
	If $a \geq 0, b\geq 0$, $\epsilon>0$, and $p> 1$ then
	\begin{equation}\label{eq:YoungIneq}
	ab \leq \frac{a^p}{\epsilon^p p} + \frac{\epsilon^q b^q}{q}, \quad \frac{1}{p} + \frac{1}{q}=1.
	\end{equation}
	Usually in this work \eqref{eq:YoungIneq} is applied with $p=q=2$ (we will say in such case that we use the Cauchy inequality), but when different $p,q$ are used then we will list the values of $p,q,\epsilon$.
\end{lemma}
We recall the following inequality.
\begin{lemma}[Poincar\'{e} inequality]\label{lem:Poincare}
	For every $u\in H^1_0$ and for every $u\in H^1$ with zero mean there holds
	\begin{equation}\label{eq:Poinc-inequality}
	\|u\|_{L^2}\leq \frac{1}{\pi}\|u_x\|_{L^2},
	\end{equation}
	where the constant $1/\pi$ in is optimal for both classes of functions.
\end{lemma}
The proofs of the following interpolation inequalities are well known. We recall them only for the completeness of the exposition.
\begin{lemma}[Embedding constant $H^1_0\subset L^\infty$]\label{lem:embedding}
	We have the inequality
	$$\|u\|_{L^\infty}\leq \frac{1}{2}\|u_x\|_{L^2}\quad \textrm{for}\quad u\in H^1_0$$
\end{lemma}
\begin{proof}
	For $x_0\in [0,1]$ there hold the bounds
	\begin{align*}
	& |u(x_0)| \leq \int_0^{x_0} |u_x(x)|\, dx,\\
	& |u(x_0)| \leq \int_{x_0}^{1} |u_x(x)|\, dx,
	\end{align*}
	and the proof follows easily.
\end{proof}	
\begin{lemma}[Interpolation inequalities]\label{lem:int}
	We have the following inequalities.
	\begin{align}
	&\|u\|_{L^\infty} \leq \|u_x\|_{L^2}^{1/2}\|u\|_{L^2}^{1/2}\quad \textrm{for}\quad u\in H^1_0,\label{eq:inter_1}\\
	&\|u\|_{L^\infty} \leq \sqrt{2}\|u_x\|_{L^2}^{1/2}\|u\|_{L^2}^{1/2}\quad \textrm{for}\quad u\in H^1\quad \int_0^1u(x)\, dx = 0,\label{eq:inter_1.5}\\
	&\|u_x\|_{L^2} \leq \|u\|_{L^2}^{1/2}\|u_{xx}\|_{L^2}^{1/2}\quad \textrm{for}\quad u\in H^2\cap H^1_0 \quad \textrm{or} \quad  u\in H^2, u_x \in H^1_0.\label{eq:inter_2}
	\end{align}
\end{lemma}
\begin{proof}
	We first prove \eqref{eq:inter_1}. For a smooth function $u:[0,1]\to \mathbb{R}$
	$$
	\frac{d}{dx}|u(x)|^2 = 2u(x)u_x(x) \quad \textrm{for}\quad x\in [0,1].
	$$
	Assume that $|u(x_0)| = \sup_{u\in[0,1]} |u(x)|$.
	Hence, as $u(0) = 0$,
	$$
	|u(x_0)|^2 = 2\int_0^{x_0} u(y) u_x(y)\, dy \leq 2\int_0^{x_0} |u(y)| |u_x(y)|\, dy.
	$$
	In a similar way, as $u(1) = 0$,
	$$
	|u(x_0)|^2 = - 2\int_{x_0}^1 u(y) u_x(y)\, dy \leq 2\int_{x_0}^1 |u(y)| |u_x(y)|\, dy.
	$$
	This means that
	$$
	2\sup_{x\in [0,1]} |u(x)|^2 = 2 |u(x_0)|^2 \leq  2\int_{0}^1|u(y)||u_x(y)|\, dy\leq 2\|u\|_{L^2}\|u_x\|_{L^2},
	$$
	whence we get the assertion. If $u$ is not smooth the assertion \eqref{eq:inter_1} follows by density. Inequality \eqref{eq:inter_1.5} follows the similar proof that uses the fact that the function which is mean free on $(0,1)$ must have a root in this interval.
	We pass to the proof of \eqref{eq:inter_2}. There holds
	$$
	\|u_x\|_{L^2}^2 = \int_0^1 u_x u_x \, dx= -\int_0^1 u_{xx}u \, dx \leq \|u_{xx}\|_{L^2}\|u\|_{L^2},
	$$
	and the proof is complete.
\end{proof}

Now we prove the important property of the trilinear term which appears in the Burgers equation.
\begin{lemma}\label{lemma:zero}
	If $u\in H^1$, satisfies $u(0) = u(1)$ then
	$$
	\int_0^1 u u_x u\, dx = 0.
	$$
	In particular, the above equality holds for $u\in H^1_0$.
\end{lemma}	
\begin{proof}
	For a smooth function $u$ defined on $[0,1]$ such that $u(0) = u(1)$ there holds the relation
	\begin{equation} \label{eq:u2ux}
	\int_0^1 u u_x u \, dx = \int_0^1 \frac{d}{dx}\left(\frac{u^3}{3}\right) \, dx = \frac{u(1)^3-u(0)^3}{3} = 0.
	\end{equation}
	By density, this relation holds also for functions from $H^1$, which have the same value on both endpoints of the interval.
\end{proof}

\subsection{A polynomial equation.} We will use several times the following lemma.
\begin{lemma}
	\label{lem:unique-sol-eq}  Let $d_i >0$ for $i=0,\dots,s$ and $0< p_i < 1$ for $i=1,\dots,s$.
	The equation \begin{equation}
	x=h(x)=d_0 + \sum_{k=1}^s d_k x^{p_i} , \label{eq:sol-eq}
	\end{equation}
	has a unique positive solution $A$, moreover $x < h(x)$ for $x<A$ and $x > h(x)$ for $x > A$. 
\end{lemma}
\begin{proof}
	Observe that $x<h(x)$ for $x \to 0$ and $x > h(x)$ for $x \to \infty$. 	Let us set
	\begin{equation}
	x_0=\sup\{x \geq 0 \, : \ t < h(t) \quad \textrm{for every}\quad t \in (0,x) \}.
	\end{equation}
	Obviously $x_0 < \infty$ and it is the smallest solution of (\ref{eq:sol-eq}). Observe that $h'(x_0)\leq 1$, because otherwise in the neighborhood of $x_0$ the function $h(x)$ will be growing faster then $x$, hence we will have $h(x_0-\delta) < x_0-\delta$ for some small $\delta>0$. This implies that there exists $\bar{x}_0 < x_0$ such that $\bar{x}_0=h(\bar{x}_{0})$. We obtain a contradiction. Hence $h'(x_0)\leq 1$.
 Since  $h''(x) <0$ for $x >0$, it follows that
	\begin{equation*}
	h'(x) < 1 , \quad x > x_0.
	\end{equation*}
	Now we rule out the existence of other solutions of equation (\ref{eq:sol-eq}). Let us take $x >x_0$, then
	\begin{eqnarray*}
		h(x)=h(x_0) + \int_{x_0}^x h'(t)dt < x_0 + \int_{x_0}^x 1 dt=x.
	\end{eqnarray*}
	The proof is complete.
\end{proof}

\subsection{A trick of X. Wang.}

The following observation  inspired by \cite{Wang2007} is quite simple, but for the role it plays in our developments we elevate it to the status
of the theorem.
\begin{theorem}\label{thm:bound_xx}
	Let  $A,B,C,D,E$ be real constants, such that $A \geq 0$ and $C>0$. Assume that we have absolutely continuous functions
	$g:[t_0,\infty) \to [0,A]$ (i.e. $g$ is bounded) and $v:[t_0,\infty) \to \mathbb{R}$   such that
	\begin{eqnarray}\label{eq:first-bound}
	\frac{d g}{dt} + Cv(t) \leq B\quad \textrm{for almost every} \quad  t \geq t_0, \\
	\label{eq:sec-bound}
	\frac{d v}{dt} \leq D + E v(t)\quad \textrm{for almost every} \quad  t \geq t_0
	\end{eqnarray}
	Then for every $0 \geq \lambda $ such that $\lambda+E \geq 0$ and for every $t \geq t_0$ there holds
	\begin{eqnarray}
		v(t)  + \frac{E+\lambda}{C} g(t) &\leq&  \left(v(t_0) + \frac{E + \lambda}{C}g(t_0) \right) e^{-\lambda (t-t_0)}  \label{eq:A-Wang-trick}  \\
		& & + \left(\frac{D}{\lambda}+\frac{E+ \lambda}{C}\left( A + \frac{B}{\lambda}\right)\right)\left(1-e^{-\lambda (t-t_0)}\right). \nonumber
	\end{eqnarray}
\end{theorem}
\begin{proof}
	We will show that for sufficiently large $F\geq 0$ the value $v + Fg$ is bounded from above.   Indeed,	if only $FC - E \geq 0$, then
	\begin{eqnarray*}
		\frac{d}{dt} (v + Fg) &\leq& D + E v + F(B - Cv)=-(FC - E)v + (D + FB) \\
		&=&-(FC - E)(v+Fg) + (D + FB) + (FC - E)Fg\\
		&\leq& -(FC - E)(v+Fg) + (D + FB) + (FC - E)FA.
	\end{eqnarray*}
	From the above inequality it follows that  $v + Fg$ is bounded from above and there holds
	\begin{align*}
		v(t) + Fg(t) &\leq (v(t_0) + Fg(t_0))e^{ -(FC - E)(t-t_0)}  \\
		&+ ((D + FB) + F(FC - E)A)\frac{ \left(1-e^{-(FC-E)(t-t_0)} \right)}{FC-E}.
	\end{align*}	
	Let us set
	\begin{equation*}
	\lambda=FC-E.
	\end{equation*}
	Then $F=\frac{\lambda +E}{C}$ and we require that $\lambda + E \geq 0$ in order to have $F\geq 0$.
	After this substitution we obtain our assertion.
\end{proof}

\begin{corollary}\label{cor:v-bounded}
  Assume that assumptions of Theorem \ref{thm:bound_xx} are satisfied.
	Then $v$   is bounded from above.
\end{corollary}
Usually the above corollary and theorem will be applied to function $v$, which is nonnegative, hence the lower bound will be automatic.

\begin{corollary}\label{cor:29}
	To find $\lambda$ such that in \eqref{eq:A-Wang-trick} the constant in front of $(1-e^{-\lambda(t-t_0)})$ is minimal consider the following cases
	\begin{itemize}
		\item Either $CD+BE \leq 0$ or
		($CD+BE > 0$ and $\sqrt{\frac{CD+BE}{A}} \leq -E$) (note that $E$ must be a negative number). Then one needs to take $\lambda = -E$ to minimize the constant in front of $\left(1\right)$, whence it follows that
		\begin{eqnarray*}
			v(t)  &\leq&  v(t_0) e^{-\lambda (t-t_0)}   -\frac{D}{E}\left(1-e^{-\lambda (t-t_0)}\right),
		\end{eqnarray*}
		Observe that the above inequality follows directly from (\ref{eq:sec-bound}).

		\item $CD+BE > 0$ and $\sqrt{\frac{CD+BE}{A}} > -E$. Then the minimal value of the constant
		in front of $\left(1-e^{-\lambda (t-t_0)}\right)$ is obtained for
		\begin{equation}
		\lambda = \overline{\lambda}=\sqrt{\frac{CD + BE}{A}}.  \label{eq:lambda-bar}
		\end{equation}
		We get the bound
		\begin{eqnarray*}
			v(t)  + \frac{E+\overline\lambda}{C} g(t) &\leq&  \left(v(t_0) + \frac{E + \overline\lambda}{C}g(t_0) \right) e^{-\overline\lambda (t-t_0)}  \\
			& & + \frac{1}{C}\left(EA+B+2 \sqrt{A(CD+BE)}\right)\left(1-e^{-\overline\lambda (t-t_0)}\right),
		\end{eqnarray*}		
	\end{itemize}
	
\end{corollary}

\begin{remark}
	One is tempted to think, that the above Corollary can improve  the estimate $g(t)\leq A$. This, however, is not the case. Indeed, after an easy calculation we obtain
	$$
	g(t) \leq \frac{\frac{1}{C}\left(EA+B+2\sqrt{A(CD+BE)}\right)}{\frac{E+\overline{\lambda}}{C}} = A + \frac{B\sqrt{A} + A\sqrt{CD+BE}}{E\sqrt{A}+\sqrt{CD+BE}}.
	$$
	The last constant is greater than or equal to $A$.
\end{remark} 

 \section{Derivarion of trapping sets radii and local estimates.}
\label{trapping_sets}

This appendix contains five subsections devoted, respectively, to the calculation of the radii of the trapping sets for the $L^2$ norm of the function and its space derivatives up to fourth. The results derived here depend heavily on the particular form of the Burgers equations, but we believe that some of the techniques will be also transferable to other problems.
 Apart from the global estimates also the local estimates are derived. For each quantity several techniques to compute the estimates are presented and in the numerical realization all estimates are computed and the smallest one is always chosen. The computations are rather standard, although sometimes technically cumbersome. We present all derivations for the exposition completeness. Summary of all results of thos section used in our computational code can be found in Tables \ref{tab:radii} and \ref{local_boundsi}.

\subsection{Trapping sets in $L^2$.}
Let $R_1 >0$ be a number. We define
	$$
W_{L^2}(R_1) = \{  v\in L^2\, :\ \|v\|_{L^2} \leq R_1 \}.
$$
We prove the following lemma on the existence of $L^2$-trapping set. Note that we need the regularity of $f$ to be only $L^\infty(L^2)$.
\begin{lemma}\label{lem:b1}
Let $f\in L^\infty(L^2)$. There exists an $L^2$-trapping set which is nonempty and bounded in $L^2$. In fact if only
	\begin{equation}\label{eq:R1}
R_1 \geq \frac{\|f\|_{L^\infty(L^2)}}{\pi^2}.
\end{equation}
then $W_{L^2}(R_1)$ is $L^2$-trapping. Moreover for every $t_1\in \mathbb{R}$ and every $t>t_1$ there holds the local estimate
	\begin{equation}\label{eq:R1_est}
 \|u(t)\|_{L^2} \leq \|u_0\|_{L^2} e^{-\pi^2(t-t_1)} + \frac{\|f\|_{L^\infty(t_1,t;L^2)}}{\pi^2}\left(1-e^{-\pi^2(t-t_1)}\right)
\end{equation}
\end{lemma}
\begin{proof}
	By the comparison of \eqref{est:1} with the solution of the ODE
	$$v'(s) = -2\pi^2 v(s) + 2\|f\|_{L^\infty(t_1,t;L^2)}\sqrt{v(t)}$$
	we obtain \eqref{eq:R1_est}.
Now, to prove that $W_{L^2}(R_1)$ is trapping let us assume that $\|u_0\|_{L^2} \leq R_1$. Then
	$$
		\|u(t)\|_{L^2} \leq
	R_1 e^{-\pi^2(t-t_1)} +  \frac{	\|f\|_{L^\infty(L^2)}}{\pi^2}\left(1-e^{-\pi^2(t-t_1)}\right)  \leq R_1,
	$$ and the proof is complete.
\end{proof}

\subsection{Trapping sets in $H^1_0$.}
In this subsection we will use the notation
$$
W_{H^1}(R_1,R_2) = \{  v\in H^1_0(\Omega)\, :\ \|v\|_{L^2} \leq R_1, \|v_x\|_{L^2} \leq R_2 \}.
$$
	\subsection{A priori estimates leading to the trapping set on $\|u_x\|_{L^2}$.} The radius of the trapping set for $\|u_x\|_{L^2}$ can be effectively calculated due to the following Lemma.
\begin{lemma}
\label{lem:abssetL2}
	Let $f\in L^\infty(L^2)$.
	There exists an $H^{1}_0$-trapping  set which is nonempty and bounded in $H^1_0$. In fact if only
	\begin{equation}
	R_1 \geq \frac{\|f\|_{L^\infty(L^2)}}{\pi^2}, \label{eq:R1-bound}
	\end{equation}
	and $R_2 = \min\left\{A/\pi,(A R_1)^{1/2}\right\}$,
 where $A$ is greater than  or  equal to  the positive root of the equation (whose existence and uniqueness follow from Lemma~\ref{lem:unique-sol-eq})
	\begin{equation}
	x - \|f\|_{L^\infty(L^2)} - R_1^{5/4}x^{3/4} = 0,  \label{eq:A-H1bound}
	\end{equation}
	then $W_{H^1}(R_1,R_2)$ is $H^1_0$-trapping.
\end{lemma}
\begin{proof}
	
Since the ball $W_{L^2}(R_1)$ is $L^2$-trapping we can assume that if $u_0\in W_{H^1}(R_1,R_2)$ then $u(t) \in W_{L^2}(R_1)$, i.e. $\|u(t)\|_{L^2} \leq R_1$ for every $t \geq t_0$. To obtain the bounds for $\|u_x\|_{L^2}$ we will use \eqref{eq:strong_for_trapping}. As the factor $\|u_{xx}\|_{L^2} - \|f(t)\|_{L^2}  - \|u\|_{L^2}^{5/4} \|u_{xx}\|_{L^2}^{3/4}$ at the right-hand side of the estimate \eqref{eq:strong_for_trapping}  depends on the norm of $\|u_{xx}\|_{L^2}$, we will get the bound in terms of $\|u_x\|_{L^2}$ using the Poincar\'e inequality or the interpolation inequality \eqref{eq:inter_2}.

By \eqref{eq:strong_for_trapping}, it is enough to find
$R_2$, such that if $\|u\|_{L_2} \leq R_1$ and $\|u_x\|_{L_2} \geq R_2$, then
\begin{equation}
  \frac{d \|u_x\|^2_{L_2}}{dt}=-2\|u_{xx}\|_{L^2}\left(\|u_{xx}\|_{L^2} - \|f(t)\|_{L^2}  - \|u\|_{L^2}^{5/4} \|u_{xx}\|_{L^2}^{3/4} \right) \leq 0.
\end{equation}
It is easy to see that if $\|u_{xx}\|_{L_2} \geq A$ and $\|u\|_{L_2} \leq R_1$, then
\begin{equation}
  \frac{d \|u_x\|^2_{L_2}}{dt} \leq 0
\end{equation}
Now, from the interpolation inequality \eqref{eq:inter_2} it follows that
\begin{eqnarray*}
  \|u_{xx}\|_{L_2} \geq \frac{\|u_x\|^2_{L_2}}{\|u\|_{L_2}}.
\end{eqnarray*}
Hence in order to have $\|u_{xx}\|_{L_2}\geq A$ it is enough to satisfy
\begin{equation*}
  \frac{\|u_x\|^2_{L_2}}{R_1} \geq A.
\end{equation*}
Hence we obtain $\|u_x\|_{L_2} \geq (AR_1)^{1/2}$.

If instead of the interpolation inequality we use the Poincar\'e inequality $
  \|u_{xx}\|_{L_2} \geq \pi \|u_x\|_{L_2},
$ then in order to have $\|u_{xx}\|_{L_2}\geq A$ it is enough to satisfy
\begin{equation*}
  \|u_x\|_{L_2} \geq \frac{A}{\pi}.
\end{equation*}
 Therefore $W_{H^1}(R_1,R_2)$ is a trapping region, and the proof is complete.
\end{proof}

\begin{remark}
We will prove that in the above lemma, if only $R_1$ is optimal, i.e.,
$$
R_1=\frac{F}{\pi^2}, \quad F=\|f\|_{L^\infty(L^2)},
$$
then the approach by the interpolation inequality yields always the better result then the approach by the Poincar\'{e} inequality. Indeed, we will prove that
\begin{equation}
(A R_1)^{1/2} < \frac{A}{\pi}, \label{eq:H1-ineq}
\end{equation}
which is equivalent to
\begin{equation}
A > \pi^2 R_1 = F.
\end{equation}
 The equation  \eqref{eq:A-H1bound} which defines $A$ is given as follows.
	\begin{equation}
	x - F - \frac{F^{5/4}}{\pi^{5/2}}x^{3/4} = 0,  \quad . \label{eq:sol-A-h1bound}
	\end{equation}
Observe that for $x=F$ there holds $	x - F - \frac{F^{5/4}}{\pi^{5/2}}x^{3/4} <0$, hence indeed the root of the above equation satisfies $A >F$. This establishes \eqref{eq:H1-ineq}.
\end{remark}

	\begin{table}[htb]
	\centering
	\begin{tabular}{|c|l|l|}
		\hline
		$\|f\|_{L^\infty(L^2)}$ & $\frac{A}{\pi}$ & $(AR_1)^{1/2}$ \\
		\hline
		$0.1$ & 0.032013 & $0.032013$ \\
		$1$ & $0.33733$ & $0.327675$ \\
		$10$ & 6.17 & $4.43$ \\
		$100$ & 34117 & 1042 \\
		\hline
	\end{tabular}
	\caption{The minimal values of  $R_2$ computed for various $\|f\|_{L^\infty(L^2)}$. Observe that always $\frac{A}{\pi} >(AR_1)^{1/2}$, the difference increases with $\|f\|_{L^\infty(L^2)}$}\label{tab:R2}
\end{table}

\subsubsection{Using the Wang's trick}
The bound of Lemma \ref{lem:abssetL2} can be possibly improved using the ideas of Theorem \ref{thm:bound_xx} and Corollary \ref{cor:29}. Indeed we have two estimates, the first one following from \eqref{eq:ener-eq}, and another one from \eqref{eq:strong_1}
\begin{align*}
&\frac{d}{dt}\|u\|_{L^2}^2 + 2\|u_x\|_{L^2}^2 \leq 2\|f\|_{L^\infty(L^2)}R_1,\\
&\frac{d}{dt}\|u_x\|_{L^2}^2 + (2-\alpha-\beta)\pi^2\|u_{x}\|_{L^2}^2 \leq {\frac{1}{\alpha}}\|f\|^2_{L^\infty(L^2)}  + \frac{7^7}{2^{16}\beta^7}R_1^{10},\\
\end{align*}
Setting
\begin{align}
&A = R_1^2,\ \  B = 2\|f\|_{L^\infty(L^2)}R_1,\ \  C = 2,\label{eq:alt}\\
& D = {\frac{\|f\|^2_{L^\infty(L^2)}}{\alpha}}  + \frac{7^7 R_1^{10}}{2^{16}\beta^7 },\ \  E = - (2-\alpha-\beta)\pi^2,\nonumber
\end{align}
one can use Corollary \ref{cor:29} and get (possibly) lower value of $R_2$, the radius of the trapping set for $\| u_x\|_{L^2}$. Note that in such case the trapping set would have the form
$$
W_{H^1}(R_1,R_2,S) = \{  v\in H^1_0(\Omega)\, :\ \|v\|_{L^2} \leq R_1, \|v_x\|_{L^2}^2 + S \|v\|_{L^2}^2\leq R_2^2 \}.
$$
This set is convex and closed and bounded in  $H^1$. The values $\alpha$ and $\beta$ can be chosen to minimize $R_2$. In all numerical examples we calculate $R_2$ according to Lemma \ref{lem:abssetL2} and to  Corollary \ref{cor:29} and we choose the lowest obtained value. We also perform the search over the discrete set of possible values of $\alpha, \beta$ which are positive to get the possibly best estimate for $\|u_x\|_{L^2}$.  In the sequel we will simply denote the found trapping set with the smallest radii by $W_{H^1}(R_1,R_2)$, remembering that it is possible that the method based on Corollary  \ref{cor:29} can produce the smaller radius $R_2$, and then the set would depend on the constant $S$.

\subsubsection{Local in time a priori estimates of $\|u_x\|_{L^2}$}

\begin{lemma}\label{lemma:young_x}
Let $f\in L^\infty(t_1,t;L^2)$ and let $\alpha,\beta > 0$ be such that $\alpha+\beta < 2$. Assume that the solution of the Burgers equation $u:[t_1,t]\to H^1_0$ satisfies the estimate
$$
\|u(s)\|_{L^2} \leq R_1 \quad \textrm{for}\quad s\in [t_1,t].
$$
Then
	\begin{align}
&  \|u_x(t)\|_{L^2}^2 \leq \|(u_0)_x\|_{L^2}^2 e^{-\pi^2(2-\alpha-\beta)(t-t_1)}  \label{ux:estimate}\\
& \qquad \qquad + \left(\frac{\|f\|_{L^\infty(t_1,t;L^2)}^2}{\alpha}+\frac{R_1^{10}7^7}{2^{16}\beta^7}\right)\frac{1}{\pi^2(2-\alpha-\beta)}(1-e^{-\pi^2(2-\alpha-\beta)(t-t_1)}).\nonumber
\end{align}
The minimal value of the constant in front of
$(1-e^{-\pi^2(2-\alpha-\beta)(t-t_0)})$ is obtained by taking
$\alpha$ as the positive root of the equation
$$
\alpha + \sqrt[4]{\frac{R_1^5}{\|f\|_{L^\infty(t_1,t;L^2)}}} \alpha^{1/4} - 1 = 0,
$$
and $\beta = 7/4 - 7\alpha/4$.
	\end{lemma}
\begin{proof}
	Application of the Gronwall lemma in \eqref{eq:strong_1} implies
		\begin{align*}
	& \nonumber \|u_x(t)\|_{L^2}^2 \leq \|(u_0)_x\|_{L^2}^2 e^{-\pi^2(2-\alpha-\beta)(t-t_1)}  \\
	& \qquad \qquad + \left(\frac{\|f\|_{L^\infty(t_1,t;L^2)}^2}{\alpha}+\frac{7^7 R_1^{10}}{2^{16}\beta^7}\right)\int_{t_1}^t e^{-\pi^2(2-\alpha-\beta)(t-s)}\, ds.
	\end{align*}
	Calculating the integral, we obtain \eqref{ux:estimate}. We minimize the expression $$\left(\frac{\|f\|_{L^\infty(t_1,t;L^2)}^2}{\alpha}+\frac{7^7 R_1^{10}}{2^{16}\beta^7}\right)$$ over the set $\{ (\alpha,\beta)\in \R^2\, :\ \alpha>0, \beta>0, \alpha+\beta < 2 \}$. A simple but cumbersome computation shows that this expression is minimal, if $\alpha$ is the positive root of the equation
$$
\alpha + \sqrt[4]{\frac{R_1^5}{\|f\|_{L^\infty(t_1,t;L^2)}}} \alpha^{1/4} - 1 = 0,
$$
and
$$
\beta = \frac{7}{4}\sqrt[4]{\frac{R_1^5}{\|f\|_{L^\infty(t_1,t;L^2)}}}\alpha^{1/4} = \frac{7}{4} - \frac{7\alpha}{4},
$$
and the assertion is proved.
	\end{proof}

\begin{remark}
	It is possible to minimize with respect to $\alpha$ and $\beta$ the whole expression on the right-hand side of \eqref{ux:estimate} (including the exponents and not just the constant in the second term).
	Let us observe that
	the righthand side of \eqref{ux:estimate} attains its minimum for
	$$
	\beta = \frac{7}{4}\sqrt[4]{\frac{R_1^5}{\|f\|_{L^\infty(t_1,t;L^2)}}}\alpha^{1/4}.
	$$
	Hence it is sufficient to minimize the radius over the set
	$$
	\left\{  (\alpha,\beta)\in \mathbb{R}^2\, :\ \alpha>0, \beta > 0, \alpha+\beta<2, \beta = \frac{7}{4}\sqrt[4]{\frac{R_1^5}{\|f\|_{L^\infty(t_1,t;L^2)}}}\alpha^{1/4} \right\},
	$$
	which is a one dimensional numerical procedure. The approximate solution of a minimization problem gives slightly better estimate of the right-hand side of \eqref{ux:estimate} than $\alpha$ and $\beta$ given in Lemma \ref{lemma:young_x}.
	
	
\end{remark}

The following Lemma provides the alternative estimate based on the interpolation inequality.
\begin{lemma}\label{lemma:young_x_2}
	Let $f\in L^\infty(t_1,t;L^2)$ and let $\alpha,\beta > 0$ be such that $\alpha+\beta < 2$. Assume that the solution of the Burgers equation $u:[t_1,t]\to H^1_0$ satisfies the estimate
	$$
	\|u(s)\|_{L^2} \leq R_1 \quad \textrm{for}\quad s\in [t_1,t].
	$$
	Then
	\begin{align}
	&  \|u_x(t)\|_{L^2}^2 \leq  \frac{D \tanh(\sqrt{CD}  (t-t_1)) + \sqrt{CD} \|u_x(t_1)\|_{L^2}^2 }{C\tanh(\sqrt{CD}  (t-t_1)) \|u_x(t_1)\|_{L^2}^2 + \sqrt{CD} }
		\end{align}
	with
	$$
	C = \frac{2-\alpha-\beta}{R_1^2} \quad \textrm{and} \quad D = {\frac{\|f\|^2_{L^2(t_1,t;L^2)}}{\alpha}}  + \frac{7^7 R_1^{10}}{2^{16}\beta^7}.
	$$
\end{lemma}
\begin{proof}
	
Using the interpolation inequality \eqref{eq:inter_2} in \eqref{eq:dux-uxx}, we obtain
\begin{equation*}
	 \frac{d}{dt}\|u_x\|_{L^2}^2 + \frac{2-\alpha-\beta}{R_1^2}\|u_{x}\|_{L^2}^4 \leq {\frac{1}{\alpha}}\|f(t)\|^2_{L^2}  + \frac{7^7}{2^{16}\beta^7}\|u\|_{L^2}^{10}.
	 \end{equation*}

Hence $z(t)=\|u_x(t)\|_{L^2}^2$ satisfies the following differential inequality
\begin{equation}
  z'(t) \leq -C z^2(t) + D,
\end{equation}
with
$$
C = \frac{2-\alpha-\beta}{R_1^2} \quad \textrm{and} \quad D = {\frac{1}{\alpha}}\|f\|^2_{L^2(t_1,t;L^2)}  + \frac{7^7}{2^{16}\beta^7}R_1^{10}
$$
It is easy to check that the solution of equation $z'(t)=-C z^2(t) +D$ with initial condition $y(t_0)=z(t_0)$
is given by
\begin{equation}
  y(t) =   \frac{D \tanh(\sqrt{CD}  (t-t_0)) + \sqrt{CD} z(t_0)}{C\tanh(\sqrt{CD}  (t-t_0)) z(t_0) + \sqrt{CD} }
\end{equation}

Hence
\begin{equation}
  z(t) \leq y(t),
\end{equation}
and the proof is complete.

\end{proof}

\begin{remark}
	Again, it is possible to minimize the resulting estimate over constants $\alpha, \beta$ in the set
	$\{  (\alpha,\beta)\in \mathbb{R}^2\, :\ \alpha>0, \beta>0, \alpha+\beta < 2  \}$ to obtain the smallest possible upper bound.
	\end{remark}

\subsection{Trapping set for $u_{xx}$ in $L^2$.}
\label{subsec:estm-uxx}

From now on we assume that $f(0,t)=f(1,t) = 0$. Then, provided $u_{xx}\in H^1$, the boundary condition $u_{xx}(0,t) = u_{xx}(1,t) = 0$ holds.
We differentiate the original equation twice with respect to $x$ and denote $v=u_{xx}$. This function satisfies the following system
\begin{align}
& v_t - v_{xx} + 3u_x v + uv_x = f_{xx},\label{uxx_eq}\\
& v(0,t) = v(1,t) = 0.\label{uxx_id}
\end{align}
If the condition $f(0,t)=f(1,t) = 0$ is not satisfied, the method that we use still works, the assumption is made to avoid technicalities, see Remark \ref{ref:remark_zero}.

 The aim of this subsection is to derive the optimal estimates for $\|u_{xx}\|_{L^2}$. The following result follows by the Galerkin method, similar as the existence and uniqueness result for the original problem given in Definitions \ref{ref:def_burg} and \ref{ref:def_burg_strong}, see \cites{Robinson, Temam}. We skip the proof for the sake of the article brevity.
\begin{lemma}
	Assume that $f\in L^\infty(H^1_0)$ (which, in particular, implies that $f(1,t)=f(0,t) = 0$), and that $u\in L^\infty(H^1_0)$. If $v_0 \in L^2$ then the problem governed by \eqref{uxx_eq} with the boundary data \eqref{uxx_id} and the initial condition $v_0$ taken at time $t_0$ has a unique weak solution with the regularity $v\in C([t_0,\infty);L^2)$.  Moreover, if $u_0\in H^2\cap H^1_0$ then, in distributional sense, $v = u_{xx}$, where $u$ is the weak solution of the problem given by Definition  \ref{ref:def_burg} with the initial data $u_0$ taken at $t_0$. If $v_0 \in H^1_0$ and $f\in L^\infty(H^2\cap H^1_0)$ then the weak solutions of the problem governed by \eqref{uxx_eq} with the boundary data \eqref{uxx_id} and the initial condition $v_0$ taken in time $t_0$ is also its strong solution with the regularity $v\in C([t_0,\infty);H^1_0)$.
\end{lemma}

\begin{lemma}\label{lem:uxx_first}
	Let $f\in L^\infty(H^1_0)$.
	There exists an $H^{2}\cap H^1_0$-trapping  set which is nonempty and bounded in $H^2$. In fact if only $R_1$ and $R_2$ are given by Lemma \ref{lem:abssetL2}
	and $R_3 = \min\left(A/\pi,(A R_2)^{1/2}\right)$, where $A$ is greater than  or  equal to the smaller number of the positive roots of the equations
	\begin{equation}
	x -  5 R_1^{1/2} R_2 x^{1/2} - \|f_x\|_{L^\infty(L_2)} = 0
	\end{equation}
	\begin{equation}\label{ineq_b16}
	x -  \frac{5 \sqrt{2}}{4} R_2^{7/4}  x^{1/4} - \|f_x\|_{L^\infty(L_2)} = 0
		\end{equation}
	then the set $W_{H^2}(R_1,R_2,R_3) = \{ u\in H^2\cap H^1_0\;:\ \|u\|_{L^2}\leq R_1, \|u_x\|_{L^2}\leq R_2, \|u_{xx}\|_{L^2}\leq R_3 \}$ is $H^2\cap H^1_0$-trapping.
	\end{lemma}
\begin{proof}	
We take the scalar product in $L^2$ of \eqref{uxx_eq} with $v$, whence
$$
\langle v_t,v\rangle - (v_{xx},v) + (3u_xv+uv_x,v) = -(f_{x},v_x)
$$
Keeping in mind that $3u_xv+uv_x = (uv+u_x^2)_x$, and both $u$ and  $v$ satisfy the Dirichlet condition at the boundary, it follows that
\begin{equation}\label{eq:bounduxx_1}
\frac{1}{2} \frac{d }{dt} \|v\|^2_{L^2} +\|v_x\|_{L^2}^2  - (uv+u_x^2,v_x)=- (f_x,v_x).
\end{equation}
Now note that
$$
(u_x^2,v_x) = -((u_x)^2_x,v) = -2(u_x v,v) = 4(uv,v_x).
$$
This means that \eqref{eq:bounduxx_1} has the following two alternative representations
$$
\frac{1}{2} \frac{d }{dt} \|v\|^2_{L^2} +\|v_x\|_{L^2}^2  - 5(uv,v_x)=- (f_x,v_x).
$$
$$
\frac{1}{2} \frac{d }{dt} \|v\|^2_{L^2} +\|v_x\|_{L^2}^2  - \frac{5}{4}(u_x^2,v_x)=- (f_x,v_x).
$$
Now
$$
  |(uv,v_x)| \leq \|u\|_{L^\infty}\|v\|_{L^2} \|v_x\|_{L_2},\qquad
  |(u_x^2,v_x)| \leq \|u_x\|_{L^\infty} \|u_x\|_{L^2} \|v_x\|_{L_2}
$$
By Lemma \ref{lem:int} we deduce that
$$
|(uv,v_x)| \leq  \|u\|_{L^2}^{1/2} \|u_x\|_{L_2}^{1/2}\|u_{xx}\|_{L^2} \|v_x\|_{L_2} \leq \|u\|_{L^2}^{1/2} \|u_x\|
\|u_{xxx}\|_{L^2}^{1/2} \|v_x\|_{L_2}  = \|u\|_{L^2}^{1/2} \|u_x\|_{L_2} \|v_x\|_{L_2}^{3/2}.
$$
$$
  |(u_x^2,v_x)|\leq  \sqrt{2}\|u_x\|_{L^2}^{3/2} \|u_{xx}\|_{L^2}^{1/2}\|v_x\|_{L_2} \leq \\
\sqrt{2}\|u_x\|_{L^2}^{3/2} \|u_x\|_{L^2}^{1/4} \|u_{xxx}\|_{L^2}^{1/4} \|v_x\|_{L_2} = \\
\sqrt{2}\|u_x\|_{L^2}^{7/4}   \|v_x\|_{L_2}^{5/4}.
$$
We obtain the following two estimates
\begin{align}
& \frac{d }{dt} \|v\|^2_{L^2} = - 2\|v_x\|_{L^2}^2  + 10(uv,v_x) - 2(f_x,v_x)\nonumber \\
& \qquad \leq - 2\|v_x\|_{L^2}^2 + 10 \|u\|_{L^2}^{1/2} \|u_x\|_{L_2} \|v_x\|_{L_2}^{3/2} +2 \|f_x(t)\|_{L_2}\|v_x\|_{L^2}\nonumber\\
&\qquad = -2 \|v_x\|_{L^2} \left( \|v_x\|_{L^2} -  5 \|u\|_{L^2}^{1/2} \|u_x\|_{L_2} \|v_x\|_{L_2}^{1/2} - \|f_x(t)\|_{L_2}\right).\label{eq:second_1}
\end{align}

\begin{align}
&\frac{d }{dt} \|v\|^2_{L^2} = - 2\|v_x\|_{L^2}^2  + \frac{10}{4}(u_x^2,v_x) - 2(f_x,v_x) \nonumber \\
& \qquad \leq - 2\|v_x\|_{L^2}^2 + \frac{10 \sqrt{2}}{4} \|u_x\|_{L^2}^{7/4}   \|v_x\|_{L_2}^{5/4} +2 \|f_x(t)\|_{L_2}\|v_x\|_{L^2}\nonumber\\
&\qquad = -2 \|v_x\|_{L^2} \left( \|v_x\|_{L^2} -  \frac{5 \sqrt{2}}{4} \|u_x\|_{L^2}^{7/4}   \|v_x\|_{L_2}^{1/4} - \|f_x(t)\|_{L_2}\right).\label{eq:second_2}
\end{align}
The rest of the proof follows the same argument as Lemma \ref{lem:abssetL2}.
\end{proof}

\subsubsection{Bounds based on Wang's trick}

In the next argument we use Corollary \ref{cor:29} to combine the above estimates of $\|v\|_{L^2} = \|u_{xx}\|_{L^2}$ with the estimates of $\|u_x\|_{L^2}$ obtained  in Lemma \ref{lem:est2} to get possibly better bound of $R_3$. To this end let us first recall the equation \eqref{eq:strong_1}
\begin{equation}\label{eq:est_uxx2}
\frac{d}{dt}\|u_x\|_{L^2}^2 + (2-\gamma-\delta)\|u_{xx}\|_{L^2}^2 \leq {\frac{1}{\gamma}}\|f(t)\|^2_{L^2}  + \frac{7^7}{2^{16}\delta^7}\|u\|_{L^2}^{10},
\end{equation}
where $\gamma, \delta>0$ are arbitrary constants such that $\gamma+\delta < 2$. As $u_{xx}=v$ and $\|u(t)\|_{L^2} \leq R_1$, this leads us to
\begin{equation*}
\frac{d}{dt}\|u_x\|_{L^2}^2 + (2-\gamma-\delta)\|v\|_{L^2}^2 \leq {\frac{\|f\|^2_{L^\infty(L^2)}}{\gamma}}  + \frac{7^7 R_1^{10}}{2^{16}\delta^7}.
\end{equation*}
The above bound together with either of the bounds of Lemma \ref{lemma:uxxauxiliary}
allow us to use Corollary \ref{cor:29}. We use this result with the following three sets of parameters
\begin{align}
&A = R_2^2,\ \  B ={\frac{\|f\|^2_{L^\infty(L^2)}}{\gamma}}  + \frac{7^7R_1^{10}}{2^{16}\delta^7},\ \  C = 2-\gamma-\delta,\label{ABCD_1}\\
& D = \frac{5^4 3^3 R_2^4 R_1^2 }{2^4\beta^3}  +\frac{\|f_x\|^2_{L^\infty(L^2)}}{\alpha},\ \  E = - (2-\alpha-\beta)\pi^2,\nonumber
\end{align}

\begin{align}
&A = R_2^2,\ \  B ={\frac{\|f\|^2_{L^\infty(L^2)}}{\gamma}}  + \frac{7^7R_1^{10}}{2^{16}\delta^7},\ \  C = 2-\gamma-\delta,\label{ABCD_2}\\
& D = \frac{3 \cdot 5^{13/3}R_2^{14/3}}{2^{28/3}\beta^{5/3}}  +\frac{\|f_x\|^2_{L^\infty(L^2)}}{\alpha},\ \  E = - (2-\alpha-\beta)\pi^2,\nonumber
\end{align}

\begin{align}
	&	 A = R_2^2,\ \  B =  {\frac{\|f\|^2_{L^\infty(L^2)}}{\gamma}}  + \frac{7^7R_1^{10}}{2^{16}\delta^7},\ \  C = 2-\gamma-\delta,\label{ABCD_3}\\
&	 D = \frac{\|f_x\|_{L^\infty(L^2)}^2}{\alpha},\ \  E= \frac{25R_1R_2}{\beta}-\pi^2(2-\alpha-\beta).\nonumber
\end{align}

We define the set
\begin{align*}
& W_{H^2}(R_1,R_2,R_3,S) \\
&\ \ = \{ v\in H^2\cap H^1_0\,:\ \|v\|_{L^2}\leq R_1, \|v_x\|_{L^2}\leq R_2, \|v_{xx}\|^2_{L^2}   + S\|v_x\|^2_{L^2} \leq R_3^2\}.
\end{align*}
Corollary \ref{cor:29} implies the following result, that states that, for appropriate choice of constants, this set is $H^2\cap H^1_0$ trapping.
\begin{lemma}\label{lem:xx}
	Assume that $f\in L^\infty(H^1_0)$. If only $R_1$ and $R_2$ are taken as in Lemma \ref{lem:abssetL2}, and
	$$
	R_3^2 \geq F(\alpha,\beta,\gamma,\delta),\ \  S = G(\alpha,\beta,\gamma,\delta)
	$$
	for some $\alpha, \beta,\gamma,\delta>0$ such that $\alpha+\beta\leq 2$ and $\gamma+\delta < 2$, where
	$$
	F(\alpha,\beta,\gamma,\delta) = \begin{cases}
	-\frac{D}{E} \ \textrm{when}\ CD+BE\leq 0 \ \textrm{or}\ CD+BE > 0 \ \textrm{and}\ \sqrt{\frac{CD+BE}{A}} \leq -E,\\
	\frac{1}{C}\left(EA+B+2 \sqrt{A(CD+BE)}\right)  \ \textrm{when}\ CD+BE > 0 \ \textrm{and}\ \sqrt{\frac{CD+BE}{A}} > -E,
	\end{cases}
	$$
	$$
	G(\alpha,\beta,\gamma,\delta ) = \begin{cases}
	0 \ \textrm{when}\ CD+BE\leq 0 \ \textrm{or}\ CD+BE > 0 \ \textrm{and}\ \sqrt{\frac{CD+BE}{A}} \leq -E,\\
	\frac{1}{C}\left(E+\sqrt{\frac{CD+BE}{A}}\right) \ \textrm{when}\ CD+BE > 0 \ \textrm{and}\ \sqrt{\frac{CD+BE}{A}} > -E,
	\end{cases}
	$$
	and
	$A, B(\gamma, \delta), C(\gamma, \delta), D(\alpha, \beta), E(\alpha,\beta)$ are given by either of three possibilities \eqref{ABCD_1}--\eqref{ABCD_3},
	then the set $W_{H^2}(R_1,R_2,R_3,S)$, is $H^2\cap H^1_0$ trapping.
\end{lemma}

\begin{remark}
	For sets of parameters \eqref{ABCD_1} and \eqref{ABCD_2} the constant $E$ is negative and hence it is sufficient to use Lemma \ref{lem:uxx_first} to get the existence of trapping set with the radius $-D/E$. It turns out, however, that Lemma \ref{lem:xx} can yield better bounds than merely Lemma \ref{lem:uxx_first}. If we use \eqref{ABCD_3}, then $E$ can be positive and then the argument of type as in  Lemma \ref{lem:uxx_first} does not work, however, \ref{lem:xx} gives us the trapping set.
\end{remark}

\begin{remark}\label{remark:min}
	The minimization of $F(\alpha,\beta,\gamma,\delta)$ can be performed numerically. The set is trapping for every
	$\alpha,\beta,\gamma,\delta > 0$ such that $\alpha+\beta \leq 2$ and $\gamma+\delta < 2$. So to obtain the smallest radius of the trapping set, one needs to search the parameter space and use the approximate minimizer in the estimates.	
\end{remark}

%

\subsubsection{Summary for global bounds on $\|u_{xx}\|_{L_2}$}
Summarizing, we obtained five possibilities to get trapping set
on $\|u_{xx}\|_{L_2}$: two by Lemma \ref{lem:uxx_first} and three by Lemma \ref{lem:xx}. In practice we compute all five radii and choose the smallest one. This smallest radius will be denoted by $R_3$ and the resulting trapping set by $W_{H^2}(R_1,R_2,R_3)$. To simplify the notation, the additional constant $S$, which enters the trapping set definition if the method of Lemma \ref{lem:xx} yields the optimal estimate, is neglected in further notation. Note that always the trapping set is convex, and, due to the bound on $\|u_{xx}\|_{L^2}$,  compact in $H^1_0$.

\subsubsection{Local estimates on $\|u_{xx}\|_{L^2}$}.
\begin{lemma}\label{lemma:uxxauxiliary}
	Assume that for $t\in [t_1,t_2]$ there hold the bounds $\|u(t)\|_{L_2}\leq R_1$ and $\|u_x(t)\|_{L_2}\leq R_1$.  Then for every $\alpha,\beta >0$ such that $\alpha+\beta\leq 2$ and for a.e. $t\in (t_1,t_2)$ there hold the following estimates:
	\begin{align}
	& \frac{d }{dt} \|v\|^2_{L^2} + (2-\alpha-\beta)\pi^2\|v\|_{L^2}^2 \leq \frac{5^4 3^3 R_2^4 R_1^2}{2^4\beta^3}  +\frac{\|f_x\|^2_{L^\infty(t_1,t_2;L^2)}}{\alpha},\label{eq:uxx_1}\\
	&\frac{d }{dt} \|v\|^2_{L^2} + (2-\alpha-\beta)\pi^2\|v\|_{L^2}^2 \leq \frac{3 \cdot 5^{13/3} R_2^{14/3} }{2^{28/3}\beta^{5/3}}  +\frac{\|f_x\|^2_{L^\infty(t_1,t_2;L^2)}}{\alpha},\label{eq:uxx_2}\\
	&\label{eq:210}
	\frac{d}{dt}\|v\|_{L^2}^2  \leq \frac{\|f_x\|_{L^\infty(t_1,t_2;L^2)}^2}{\alpha}  + \left(\frac{25R_1R_2}{\beta}+\pi^2(\alpha+\beta-2)\right)\|v\|^2_{L^2}.
		\end{align}
	\end{lemma}
\begin{proof}
Estimates \eqref{eq:second_1} and \eqref{eq:second_2} can be rewritten as
\begin{align*}
& \frac{d }{dt} \|v\|^2_{L^2} + 2\|v_x\|_{L^2}^2 \leq 10 \|u\|_{L^2}^{1/2} \|u_x\|_{L_2} \|v_x\|_{L_2}^{3/2} +2 \|f_x\|_{L^\infty(t_1,t_2;L^2)}^2\|v_x\|_{L^2},\\
&\frac{d }{dt} \|v\|^2_{L^2} + 2\|v_x\|_{L^2}^2 \leq \frac{5 \sqrt{2}}{2} \|u_x\|_{L^2}^{7/4}   \|v_x\|_{L_2}^{5/4} +2 \|f_x\|_{L^\infty(t_1,t_2;L^2)}^2\|v_x\|_{L^2}.
\end{align*}
After using the Young inequality they get the form
\begin{align}
& \frac{d }{dt} \|v\|^2_{L^2} + (2-\alpha-\beta)\|v_x\|_{L^2}^2 \leq \|u_x\|_{L^2}^4 \|u\|_{L^2}^2\frac{5^4 3^3 }{2^4\beta^3}  +\frac{\|f_x\|_{L^\infty(t_1,t_2;L^2)}^2}{\alpha}\label{est_uxx_aux_1},\\
&\frac{d }{dt} \|v\|^2_{L^2} + (2-\alpha-\beta)\|v_x\|_{L^2}^2 \leq\|u_x\|_{L^2}^{14/3} \frac{3 \cdot 5^{13/3}}{2^{28/3}\beta^{5/3}}  +\frac{\|f_x\|_{L^\infty(t_1,t_2;L^2)}^2}{\alpha}\label{est_uxx_aux_2}.
\end{align}
After using the previously obtained radii of the trapping sets for $\|u\|_{L^2}$ and $\|u_{x}\|_{L^2}$, and the Poincar\'{e} inequality this yields \eqref{eq:est_uxx1} and \eqref{eq:est_uxx2}.
	On the other hand, multiplying \eqref{uxx_eq} by $v$ and integrating over $(0,1)$ we obtain
$$
\frac{1}{2}\frac{d}{dt}\|v\|_{L^2}^2 + \|v_x\|_{L^2}^2 + 3\int_0^1u_xvv\, dx + \int_0^1uv_xv\, dx = - \int_0^1 f_x v_x\, dx.
$$
Integrating by parts and using the Schwarz inequality, we obtain
$$
\frac{1}{2}\frac{d}{dt}\|v\|_{L^2}^2 + \|v_x\|_{L^2}^2  \leq \|f_x(t)\|_{L^2}\|v_x\|_{L^2} + 5\int_0^1uv_xv\, dx.
$$
It follows that
$$
\frac{1}{2}\frac{d}{dt}\|v\|_{L^2}^2 + \|v_x\|_{L^2}^2  \leq \|f_x(t)\|_{L^2}\|v_x\|_{L^2} + 5\|u\|_{L^\infty}\|v\|_{L^2}\|v_x\|_{L^2}.
$$
From Lemma \ref{lem:int}, and the Young inequality, we have
\begin{eqnarray*}
	\|f_x(t)\|_{L^2}\|v_x\|_{L^2} &\leq&  \frac{\|f_x(t)\|_{L^2}^2}{2\alpha} + \frac{\alpha}{2}\|v_x\|_{L^2}^2, \\
	5\|u\|_{L^\infty}\|v\|_{L^2}\|v_x\|_{L^2} &\leq& \frac{\left( 5\|u\|^{1/2}_{L^2}  \|u_x\|_{L^2}^{1/2} \|v\|_{L^2}\right)^2}{2\beta} + \frac{\beta}{2} \|v_x\|_{L^2}^2,
\end{eqnarray*}
where $\alpha, \beta$ are positive constants. Hence we deduce
\begin{equation}\label{eq:est_uxx1}
\frac{d}{dt}\|v\|_{L^2}^2 + (2-\alpha-\beta)\|v_x\|_{L^2}^2  \leq \frac{\|f_x(t)\|_{L^2}^2}{\alpha}  + \frac{25\|u\|_{L^2}\|u_x\|_{L^2}\|v\|^2_{L^2}}{\beta}.
\end{equation}
Since $\|u(t)\|_{L^2} \leq R_1$ and $\|u_x(t)\|_{L^2} \leq R_2$ for every $t\geq t_1$ and $v=u_{xx}$, assuming that $\alpha+\beta \leq 2$, by the Poincar\'{e} inequality it follows that
$$
\frac{d}{dt}\|v\|_{L^2}^2 + \pi^2(2-\alpha-\beta)\|v\|_{L^2}^2  \leq \frac{\|f_x\|_{L^\infty(t_1,t_2;L^2)}^2}{\alpha}  + \frac{25R_1R_2\|v\|^2_{L^2}}{\beta},
$$
and the assertion \eqref{eq:210} follows.
\end{proof}

We pass to the result which gives alternative to the local in time estimates on the quantity $\|u_{xx}\|_{L^2}^2$  of Lemma \ref{lemma:uxxauxiliary}. The lemma is based on the estimates \eqref{est_uxx_aux_1} and \eqref{est_uxx_aux_2}, similar as Lemma \ref{lemma:uxxauxiliary} and is analogous to Lemma \ref{lemma:young_x_2}, namely interpolation inequality is used in place of the Poincar\'{e} inequality. Indeed, using the interpolation inequality $\|u_{xx}\|_{L^2}^2 \leq \|u_x\|_{L^2}\|u_{xxx}\|_{L^2} \leq R_2 \|u_{xxx}\|_{L^2}$ in \eqref{est_uxx_aux_1} and \eqref{est_uxx_aux_2} we obtain

\begin{align}
& \frac{d }{dt} \|v\|^2_{L^2} + \frac{2-\alpha-\beta}{R_2^2}\|v\|_{L^2}^4 \leq \frac{5^4 3^3 R_2^4 R_1^2}{2^4\beta^3}  +\frac{\|f_x\|^2_{L^\infty(t_1,t_2;L^2)}}{\alpha},\\
&\frac{d }{dt} \|v\|^2_{L^2} +  \frac{2-\alpha-\beta}{R_2^2}\|v_x\|_{L^2}^4 \leq  \frac{3 \cdot 5^{13/3} R_2^{14/3}}{2^{28/3}\beta^{5/3}}  +\frac{\|f_x\|^2_{L^\infty(t_1,t_2;L^2)}}{\alpha}.
\end{align}

\begin{lemma}\label{lemma:young_x_3}
	Let $f\in L^\infty(t_1,t_2;H^1_0)$ and let $\alpha,\beta > 0$ be such that $\alpha+\beta \leq 2$. Assume that the solution of the Burgers equation $u:[t_0,t]\to H^1_0$ satisfies the estimates
	$$
	\|u(s)\|_{L^2} \leq R_1,\,\, \|u_x(s)\|_{L^2} \leq R_2 \quad \textrm{for}\quad s\in [t_1,t_2].
	$$
	Then
	\begin{align}
	&  \|u_{xx}(t)\|_{L^2}^2 \leq  \frac{D \tanh(\sqrt{CD}  (t-t_1)) + \sqrt{CD} \|u_{xx}(t_1)\|_{L^2}^2 }{C\tanh(\sqrt{CD}  (t-t_1)) \|u_{xx}(t_1)\|_{L^2}^2 + \sqrt{CD} } \ \ \textrm{for}\ \ t\in (1_1,t_2),
	\end{align}
	with
	\begin{align*}
	&	C = \frac{2-\alpha-\beta}{R_2^2} \quad \textrm{and either} \quad D = \frac{5^4 3^3 R_2^4 R_1^2}{2^4\beta^3}  +\frac{\|f_x\|^2_{L^\infty(t_1,t_2;L^2)}}{\alpha} \\
		&\qquad \qquad \qquad  \textrm{or} \quad D= \frac{3 \cdot 5^{13/3} R_2^{14/3}}{2^{28/3}\beta^{5/3}}  +\frac{\|f_x\|^2_{L^\infty(t_1,t_2;L^2)}}{\alpha}.
	\end{align*}
\end{lemma}

\subsection{Trapping set for $u_{xxx}$ in $L^2$.}
\label{subsec:u3x-bnd}
In this section we establish global and local estimates for $\|u_{xxx}\|_{L^2}$.  We will use the notation
$$X = \{ u\in H^3\cap H^1_0\,:\ u_{xx}\in H^1_0\}.$$ Define the set
\begin{align*}
& W_{H^3}(R_1,R_2,R_3,R_4) \\
&\ \ = \{ v\in X\,:\  \|v\|_{L^2}\leq R_1, \|v_x\|_{L^2}\leq R_2, \|v_{xx}\|_{L^2}  \leq R_3, \|v_{xxx}\|_{L^2}\leq R_4\}.
\end{align*}

\begin{lemma}
	\label{lem:u3x-bnd}
	Assume that $f\in L^\infty(H^2 \cap H^1_0)$. There exists the $X$-trapping set which is nonempty and bounded in $H^3$. In fact if only $R_1, R_2, R_3$, are as in Section \ref{subsec:estm-uxx}
	and
	$$
	R_4 \geq \min\left(A/\pi,(A R_3)^{1/2}\right),
	$$
	where $A$ is the smaller of the positive roots of two equations
\begin{equation}\label{eq:u3_root}
	x  - \frac{7}{2}R_2 R_3^{3/4}x^{1/4} - \|f_{xx}\|_{L^\infty(L^2)} = 0,
	\end{equation}
\begin{equation}\label{eq:u3_root_2}
x - 7R_1^{1/2}R_2^{1/2} R_3^{1/2}x^{1/2}  - \|f_{xx}\|_{L^\infty(L^2)} = 0,
\end{equation}
	then the set $W_{H^3}(R_1,R_2,R_3,R_4)$, is $X$-trapping.
\end{lemma}
\begin{proof}
	Multiplying \eqref{uxx_eq} by $-v_{xx}$ and integrating over the space interval $(0,1)$ we get the bound
	$$
	\frac{1}{2}\frac{d}{dt}\|v_x\|_{L^2}^2 + \|v_{xx}\|_{L^2}^2 - 3(u_xv,v_{xx}) - (uv_x,v_{xx}) = -(f_{xx},v_{xx}).
	$$
	Let us integrate by parts
	\begin{align*}
		&(uv_x,v_{xx}) = \left(u,\frac{1}{2}\frac{d}{dx}v_x^2\right) = -\frac{1}{2}(u_xv_x,v_x) = \frac{1}{2}(v,v_x,v) + \frac{1}{2}(u_x,v_{xx}v)\\
		&\qquad  = \frac{1}{6}\int_0^1\frac{d}{dx}v^3\, dx + \frac{1}{2}(u_xv,v_{xx}) = \frac{1}{2}(u_x,v_{xx}v).
	\end{align*}
	This means that the above equation can be rewritten in the following two possible ways
	\begin{align}
	& 	\frac{1}{2}\frac{d}{dt}\|v_x\|_{L^2}^2 + \|v_{xx}\|_{L^2}^2 - \frac{7}{2}(u_xv,v_{xx}) = -(f_{xx},v_{xx}),\label{eq:aux_11}\\
	& 	\frac{1}{2}\frac{d}{dt}\|v_x\|_{L^2}^2 + \|v_{xx}\|_{L^2}^2 - 7(uv_x,v_{xx}) = -(f_{xx},v_{xx}).
	\end{align}
	Using Lemma~\ref{lem:int} we estimate the scalar products above as follows (we want to get rid of $\|v_x\|_{L^2}$, but we are happy
with $\|v_{xx}\|_{L^2}^p$ as long as $p < 2$)
	\begin{align*}
	& |(u_x v, v_{xx})|\leq \|u_x\|_{L^2} \|v\|_{L^\infty} \|v_{xx}\|_{L^2} \leq \|u_x\|_{L^2}  \|v_x\|_{L^2}^{1/2}  \|v\|_{L^2}^{1/2} \|v_{xx}\|_{L^2}\\
	&\qquad  \leq
	\|u_x\|_{L^2}  \|v\|_{L^2}^{1/4}\|v_{xx}\|_{L^2}^{1/4} \|v\|_{L^2}^{1/2} \|v_{xx}\|_{L^2}=
	\|u_x\|_{L^2}  \|v\|_{L^2}^{3/4}\|v_{xx}\|_{L^2}^{5/4}
	\end{align*}
	and
	\begin{align*}
	&|(uv_x,v_{xx})| \leq \|u\|_{L^\infty} \|v_x\|_{L^2} \|v_{xx}\|_{L^2}  \\
	&\leq \|u\|_{L^2}^{1/2} \|u_x\|_{L^2}^{1/2}  \|v\|_{L_2}^{1/2} \|v_{xx}\|_{L^2}^{1/2} \|v_{xx}\|_{L^2} =
	\|u\|_{L^2}^{1/2} \|u_x\|_{L^2}^{1/2}  \|v\|_{L_2}^{1/2} \|v_{xx}\|_{L^2}^{3/2}
	\end{align*}
	It follows that
	\begin{align*}
	&\frac{1}{2}\frac{d}{dt}\|v_x\|_{L^2}^2 + \|v_{xx}\|_{L^2}^2 \leq  \frac{7}{2}\|u_x\|_{L^2}\|v\|_{L^2}^{3/4}\|v_{xx}\|_{L^2}^{5/4} +
	\|f_{xx}\|_{L^2}\|v_{xx}\|_{L^2}.
	\end{align*}
		\begin{align*}
	&\frac{1}{2}\frac{d}{dt}\|v_x\|_{L^2}^2 + \|v_{xx}\|_{L^2}^2  \leq  7\|u\|_{L^2}^{1/2}\|u_x\|_{L^2}^{1/2}\|v\|_{L^2}^{1/2}\|v_{xx}\|_{L^2}^{3/2} +
	\|f_{xx}\|_{L^2}\|v_{xx}\|_{L^2}.
	\end{align*}
	Using the fact that  $\|u(t)\|_{L^2} \leq R_1$, $\|u_x(t)\|_{L^2} \leq R_2$, $\|v\|_{L^2}=\|u_{xx}\|_{L^2} \leq R_3$, we obtain
	\begin{align*}
		&\frac{1}{2}\frac{d}{dt}\|v_x\|_{L^2}^2 \leq - \|v_{xx}\|_{L^2}^2
		+ \frac{7}{2}R_2 R_3^{3/4}\|v_{xx}\|_{L^2}^{5/4} +
		\|f_{xx}\|_{L^\infty(L^2)}\|v_{xx}\|_{L^2}  \\
		& \qquad
		= - \|v_{xx}\|_{L^2} \left( \|v_{xx}\|_{L^2} - \frac{7}{2}R_2 R_3^{3/4}\|v_{xx}\|_{L^2}^{1/4} -
		\|f_{xx}\|_{L^\infty(L^2)} \right).
	\end{align*}
	
		\begin{align*}
	&\frac{1}{2}\frac{d}{dt}\|v_x\|_{L^2}^2 \leq - \|v_{xx}\|_{L^2}^2
	+ 7	R_1^{1/2}R_2^{1/2} R_3^{1/2}\|v_{xx}\|_{L^2}^{3/2} +
	\|f_{xx}\|_{L^\infty(L^2)}\|v_{xx}\|_{L^2}  \\
	& \qquad
	= - \|v_{xx}\|_{L^2} \left( \|v_{xx}\|_{L^2} - 7R_1^{1/2}R_2^{1/2} R_3^{1/2}\|v_{xx}\|_{L^2}^{1/2} -
	\|f_{xx}\|_{L^\infty(L^2)} \right).
	\end{align*}
	Now, we need to find the positive numbers which are roots of the equations \eqref{eq:u3_root} and \eqref{eq:u3_root_2} (by Lemma~\ref{lem:unique-sol-eq} these roots are unique).
	As $v \in H^1_0$, i.e. then $v_x$ has zero mean, the Poincar\'e equality $\pi\|v_x\|_{L^2} \leq \|v_{xx}\|$ holds. Moreover we have the interpolation inequality $\|v_x\|_{L^2}\leq \|v\|_{L^2}^{1/2}\|v_{xx}\|_{L^2}^{1/2}\leq R_3^{1/2}\|v_{xx}\|_{L^2}^{1/2}$. Proceeding exactly as in the proof of Lemma \ref{lem:abssetL2} we obtain the assertion of the Lemma.
\end{proof}

\subsubsection{Using Wang's trick}

Let us rewrite the equations \eqref{est_uxx_aux_1} and \eqref{est_uxx_aux_2} as
\begin{align}
& \frac{d }{dt} \|v\|^2_{L^2} + (2-\delta-\gamma)\|v_x\|_{L^2}^2 \leq \frac{5^4 3^3R_2^4 R_1^2 }{2^4\delta^3}  +\frac{\|f_x\|^2_{L^\infty(L^2)}}{\gamma},\label{eq_uxxx_aux_1}\\
&\frac{d }{dt} \|v\|^2_{L^2} + (2-\delta-\gamma)\|v_x\|_{L^2}^2 \leq  \frac{3 \cdot 5^{13/3} R_2^{14/3}}{2^{28/3}\delta^{5/3}} +\frac{\|f_x\|^2_{L^\infty(L^2)}}{\gamma}\label{eq_uxxx_aux_2}.
\end{align}
We can combine either of the above two estimates with either of estimates \eqref{eq:uxxx} and \eqref{eq:uxxx_2} and use Corollary \ref{cor:29} to get possibly smaller value of $R_4$. This allows us to define the following four sets of parameters
\begin{align}
&A = R_3^2,\ \ B =\frac{5^4 3^3 R_2^4 R_1^2}{2^4\delta^3}  +\frac{\|f_x\|^2_{L^\infty(L^2)}}{\gamma},\ \  C = 2-\delta-\gamma,\label{ABCD_4}\\
& D = \frac{3 \cdot 7^{8/3} 5^{5/3}R_2^{8/3}R_3^{2}}{2^{25/3}\alpha^{5/3}}+\frac{\|f_{xx}\|_{L^\infty(L^2)}^2}{\beta},\ \  E = - \pi^2(2-\alpha-\beta),\nonumber
\end{align}

\begin{align}
&A = R_3^2, B = \frac{3 \cdot 5^{13/3} R_2^{14/3}}{2^{28/3}\delta^{5/3}} +\frac{\|f_x\|^2_{L^\infty(L^2)}}{\gamma},\ \  C = 2-\delta-\gamma,\label{ABCD_5}\\
& D =\frac{3 \cdot 7^{8/3} 5^{5/3}R_2^{8/3}R_3^{2}}{2^{25/3}\alpha^{5/3}}+\frac{\|f_{xx}\|_{L^\infty(L^2)}^2}{\beta},\ \  E = - \pi^2(2-\alpha-\beta),\nonumber
\end{align}

\begin{align}
&A = R_3^2,\ \  B =\frac{5^4 3^3 R_2^4 R_1^2}{2^4\delta^3}  +\frac{\|f_x\|^2_{L^\infty(L^2)}}{\gamma},\ \  C = 2-\delta-\gamma,\label{ABCD_6}\\
& D =\frac{7^43^3R_1^2R_2^2R_3^2}{2^4\alpha^3}+\frac{\|f_{xx}\|_{L^\infty(L^2)}^2}{\beta}, \ \ E = - \pi^2(2-\alpha-\beta),\nonumber
\end{align}

\begin{align}
&A = R_3^2,\ \  B = \frac{3 \cdot 5^{13/3} R_2^{14/3}}{2^{28/3}\delta^{5/3}} +\frac{\|f_x\|^2_{L^\infty(L^2)}}{\gamma},\ \  C = 2-\delta-\gamma,\label{ABCD_7}\\
& D =\frac{7^43^3R_1^2R_2^2R_3^2}{2^4\alpha^3}+\frac{\|f_{xx}\|_{L^\infty(L^2)}^2}{\beta},\ \  E = - \pi^2(2-\alpha-\beta),\nonumber
\end{align}
We define the set
\begin{align*}
& W_{H^3}(R_1,R_2,R_3,R_4,S) \\
&\ \ = \{ v\in X\,:\ \|v\|_{L^2}\leq R_1, \|v_x\|_{L^2}\leq R_2, \|v_{xx}\|_{L^2}  \leq R_3, \|v_{xxx}\|_{L^2}^2+S\|v_{xx}\|_{L^2}^2\leq R_4\}.
\end{align*}
Similar as in the estimate on $\|u_{xx}\|_{L^2}$ Corollary \ref{cor:29} implies the following result which states that this set is $X$ trapping for appropriate $R_4$ and $S$.
\begin{lemma}\label{lem:xxx}
	Assume that $f\in L^\infty(H^2\cap H^1_0)$. There exists the $X$ trapping set which is nonempty and bounded in $H^3$. In fact if only $R_1$, $R_2$, $R_3$ are taken as in Section  \ref{subsec:estm-uxx}, and
	$$
	R_4^2 \geq F(\alpha,\beta,\gamma,\delta),\ \  S = G(\alpha,\beta,\gamma,\delta)
	$$
	for some $\alpha, \beta,\gamma,\delta>0$ such that $\alpha+\beta\leq 2$ and $\gamma+\delta < 2$, where
	$$
	F(\alpha,\beta,\gamma,\delta) = \begin{cases}
	-\frac{D}{E} \ \textrm{when}\ CD+BE\leq 0 \ \textrm{or}\ CD+BE > 0 \ \textrm{and}\ \sqrt{\frac{CD+BE}{A}} \leq -E,\\
	\frac{1}{C}\left(EA+B+2 \sqrt{A(CD+BE)}\right)  \ \textrm{when}\ CD+BE > 0 \ \textrm{and}\ \sqrt{\frac{CD+BE}{A}} > -E,
	\end{cases}
	$$
	$$
	G(\alpha,\beta,\gamma,\delta ) = \begin{cases}
	0 \ \textrm{when}\ CD+BE\leq 0 \ \textrm{or}\ CD+BE > 0 \ \textrm{and}\ \sqrt{\frac{CD+BE}{A}} \leq -E,\\
	\frac{1}{C}\left(E+\sqrt{\frac{CD+BE}{A}}\right) \ \textrm{when}\ CD+BE > 0 \ \textrm{and}\ \sqrt{\frac{CD+BE}{A}} > -E,
	\end{cases}
	$$
	and
	$A, B(\delta, \gamma), C(\delta, \gamma), D(\alpha, \beta), E(\alpha,\beta)$ are given by either of four possibilities \eqref{ABCD_4}--\eqref{ABCD_7},
	then the set $W_{H^3}(R_1,R_2,R_3,R_4,S)$, is $X$ trapping.
\end{lemma}

\subsubsection{Summary of global bounds}
Similar as in the case of $\|u_{xx}\|_{L^2}$ we get six trapping set bounds on $\|u_{xxx}\|_{L^2}$ by Lemmas \ref{lem:u3x-bnd} and \ref{lem:xxx}. Since we are interested to know the best possible bound, we calculate all of them, optimizing all four bounds of Lemma \ref{lem:xxx} with respect to constants $\alpha, \beta, \gamma, \delta$ and we choose the best obtained bound for further computations.

\subsubsection{Local estimates on $\|u_{xxx}\|_{L^2}^2$}

In the following Lemma we obtain the differential inequality on $\|v_x\|^2$, which will be useful to get the localized in time estimates.

\begin{lemma}\label{lem_3.18}
	Suppose that on interval $[t_1,t_2]$ there hold bounds $\|u(t)\|_{L^2}\leq R_1$, $\|u_x(t)\|_{L^2}\leq R_2$, $\|u_{xx}(t)\|_{L^2}=\|v(t)\|_{L^2}\leq R_3$. Then for every $\alpha, \beta> 0$ such that $\alpha+\beta\leq 2$ and for a.e. $t\in (t_1,t_2)$ there hold the bounds
	\begin{equation}\label{eq:uxxx}
		\frac{d}{dt}\|v_x\|_{L^2}^2 + \pi^2(2-\alpha-\beta)\|v_{x}\|_{L^2}^2 \leq\frac{3 \cdot 7^{8/3} 5^{5/3}R_2^{8/3}R_3^{2}}{2^{25/3}\alpha^{5/3}}+\frac{\|f_{xx}\|_{L^\infty(t_1,t_2;L^2)}^2}{\beta}.
	\end{equation}
		\begin{equation}\label{eq:uxxx_2}
	\frac{d}{dt}\|v_x\|_{L^2}^2 + \pi^2(2-\alpha-\beta)\|v_{x}\|_{L^2}^2 \leq \frac{7^43^3R_1^2R_2^2R_3^2}{2^4\alpha^3}+\frac{\|f_{xx}\|_{L^\infty(t_1,t_2;L^2)}^2}{\beta}.
	\end{equation}
	\end{lemma}
\begin{proof}
		Similar as in the proof of Lemma \ref{lem:u3x-bnd} we multiply \eqref{uxx_eq} by $-v_{xx}$ and integrate over the space interval $(0,1)$. This gives us the bounds
	\begin{align*}
	&\frac{d}{dt}\|v_x\|_{L^2}^2 + 2\|v_{xx}\|_{L^2}^2 \leq  7R_2R_3^{3/4}\|v_{xx}\|_{L^2}^{5/4} +
	2\|f_{xx}\|_{L^\infty(t_1,t_2;L^2)}\|v_{xx}\|_{L^2} = I_1+I_3,
	\end{align*}
		\begin{align*}
	&\frac{d}{dt}\|v_x\|_{L^2}^2 + 2\|v_{xx}\|_{L^2}^2 \leq  14 R_1^{1/2}R_2^{1/2}R_3^{1/2}\|v_{xx}\|_{L^2}^{3/2} +
	2\|f_{xx}\|_{L^\infty(t_1,t_2;L^2)}\|v_{xx}\|_{L^2} = I_2+I_3.
	\end{align*}.
Observe that the terms on the right-hand side of both estimates contain $\|v_{xx}\|$ to power lower than $2$, therefore  we would like use the Young inequality to majorize them by $\|v_{xx}\|^2$.
	\begin{itemize}
		\item[(i)] \textit{Estimate of $I_1$.} We use the Young inequality \eqref{eq:YoungIneq} with $a=7R_2R_3^{3/4}$, $b=\|v_{xx}\|_{L^2}^{5/4}$, $q=8/5$, $p=8/3$ and $\epsilon=\left(\alpha 8/5\right)^{5/8}$. This yields
		\begin{equation}\label{eq:uxxx-fterm}
		I_1 = 7R_2R_3^{3/4}\|v_{xx}\|_{L^2}^{5/4} \leq \frac{3 \cdot 7^{8/3} 5^{5/3}R_2^{8/3}R_3^{2}}{2^{25/3}\alpha^{5/3}}   + \alpha \|v_{xx}\|_{L^2}^2.
		\end{equation}
		
		\item[(ii)] \textit{Estimate of $I_2$.} Again we use the Young inequality \eqref{eq:YoungIneq} with  $a=14R_1^{1/2}R_2^{1/2}R_3^{1/2}$,  $b=\|v_{xx}\|_{L^2}^{3/2}$, $q=4/3$, $p=4$ and $\epsilon=\left(\alpha4/3\right)^{3/4}$. We obtain
		$$
		I_2= 14R_1^{1/2}R_2^{1/2}R_3^{1/2}\|v_{xx}\|_{L^2}^{3/2} \leq \frac{7^43^3R_1^2R_2^2R_3^2}{2^4\alpha^3}    + \alpha\|v_{xx}\|_{L^2}^2
		$$
		
		\item[(iii)]   \textit{Estimate of $I_3$.} This time we set $p=q=2$, $\epsilon=\sqrt{2\gamma}$  and we obtain
		\begin{equation*}
		I_3 = 2\|f_{xx}\|_{L^\infty(t_1,t_2;L^2)}^2\|v_{xx}\|_{L^2} \leq  \frac{\|f_{xx}\|_{L^\infty(t_1,t_2;L^2)}^2}{\beta} + \beta\|v_{xx}\|_{L^2}^2
		\end{equation*}
	\end{itemize}
	Hence, by the Poincar\'{e} inequality we deduce the assertion of the Lemma.
	\end{proof}

Analogously to Lemmas \ref{lemma:young_x_2} and \ref{lemma:young_x_3} we formulate a result which gives local in time estimates on the quantity $\|u_{xxx}\|_{L^2}^2$ alternative with respect to the ones which follow from Lemma \ref{lem_3.18}. The lemma is based on the estimates of Lemma \ref{lem_3.18} but we estimate the term $\|v_{xx}\|_{L_2} = \|u_{xxxx}\|_{L^2}$ from below not by the Poinvar\'{e} inequality, but by interpolation inequalities. Hence, using the interpolation inequality $\|u_{xxx}\|_{L^2}^2 \leq \|u_{xx}\|_{L^2}\|u_{xxxx}\|_{L^2} \leq R_3 \|u_{xxxx}\|_{L^2} = R_3 \|v_{xx}\|_{L^2}$ we obtain

	\begin{equation}\label{eq:uxxx_tan}
\frac{d}{dt}\|v_x\|_{L^2}^2 + \frac{2-\alpha-\beta}{R_3^2}\|v_{x}\|_{L^2}^4 \leq\frac{3 \cdot 7^{8/3} 5^{5/3}R_2^{8/3}R_3^{2}}{2^{25/3}\alpha^{5/3}}+\frac{\|f_{xx}\|_{L^\infty(t_1,t_2;L^2)}^2}{\beta}.
\end{equation}
\begin{equation}\label{eq:uxxx_2_tan}
\frac{d}{dt}\|v_x\|_{L^2}^2 + \frac{2-\alpha-\beta}{R_3^2}\|v_{x}\|_{L^2}^4 \leq \frac{7^43^3R_1^2R_2^2R_3^2}{2^4\alpha^3}+\frac{\|f_{xx}\|_{L^\infty(t_1,t_2;L^2)}^2}{\beta}.
\end{equation}

\begin{lemma}\label{lemma:young_x_4}
	Let $f\in L^2(t_1,t_2;H^2\cap H^1_0)$ and let $\alpha,\beta > 0$ be such that $\alpha+\beta < 2$. Assume that the solution of the Burgers equation $u:[t_1,t_2]\to X$ satisfies the estimates
	$$
	\|u(s)\|_{L^2} \leq R_1, \|u_x(s)\|_{L^2} \leq R_2, \|u_{xx}(s)\|_{L^2}\leq R_3 \quad \textrm{for}\quad s\in [t_1,t_2].
	$$
	Then
	\begin{align}
	&  \|u_{xxx}(t)\|_{L^2}^2 \leq  \frac{D \tanh(\sqrt{CD}  (t-t_1)) + \sqrt{CD} \|u_{xxx}(t_1)\|_{L^2}^2 }{C\tanh(\sqrt{CD}  (t-t_1)) \|u_{xxx}(t_1)\|_{L^2}^2 + \sqrt{CD} }
	\end{align}
	with
	\begin{align}
	&C = \frac{2-\alpha-\beta}{R_3^2} \quad \textrm{and either} \quad D = \frac{3 \cdot 7^{8/3} 5^{5/3}R_2^{8/3}R_3^{2}}{2^{25/3}\alpha^{5/3}}+\frac{\|f_{xx}\|_{L^\infty(t_1,t_2;L^2)}^2}{\beta}\\
	& \qquad \qquad  \quad \textrm{or} \quad D=\frac{7^43^3R_1^2R_2^2R_3^2}{2^4\alpha^3}+\frac{\|f_{xx}\|_{L^\infty(t_1,t_2;L^2)}^2}{\beta}.
	\end{align}
\end{lemma}

\subsection{Trapping set for $u_{xxxx}$ in $L^2$} The last estimates will be the ones of $\|u_{xxxx}\|_{L^2}$. Similar as in previous situations we will get a local and a global estimates of this quantity. Let us differentiate the original equation four times with respect to the space variable and denote $w = u_{xxxx}$ and $v=u_{xx}$. This procedure is valid provided we reinforce the previous assumptions by $f_{xx}(0,t) = 0$ and   $f_{xx}(1,t) = 0$. After differentiation we obtain the following equation.
$$
w_t - w_{xx} + 5 u_x w + 10 vv_x + uw_x = f_{xxxx}
$$
with the boundary conditions
$$
w(0,t) = w(1,t) = 0.
$$
We test this equation with $w$ which yields
$$
\frac{1}{2}\frac{d}{dt}\|w\|_{L^2}^2 + \|w_x\|_{L^2}^2 + 5 (u_x w, w) + 10 (v_xv,w) + (uw_x,w) = (f_{xxxx},w).
$$
Performing integration by parts in the last term on the left-hand side, we deduce
\begin{equation}\label{eq:u4}
\frac{1}{2}\frac{d}{dt}\|w\|_{L^2}^2 + \|w_x\|_{L^2}^2 + \frac{9}{2} (u_x w, w) + 10 (v_xv,w) = (f_{xxxx},w).
\end{equation}
By the Poincar\'{e} inequality
$$
\frac{1}{2}\frac{d}{dt}\|w\|_{L^2}^2 + \pi^2\|w\|_{L^2}^2 \leq \frac{9}{2} \|u_x\|_{L^\infty} \|w\|_{L^2}^2 + 10\|v\|_{L^\infty}\|v_x\|_{L^2}\|w\|_{L^2}+ \|f_{xxxx}(t)\|_{L^2}\|w\|_{L^2}.
$$
Multiplying by $2$ and using the Cauchy inequality with $\epsilon$, we obtain for any $\beta,\gamma >0$
$$
\frac{d}{dt}\|w\|_{L^2}^2 + (2\pi^2-\beta-\gamma - 9 \|u_x\|_{L^\infty})\|w\|_{L^2}^2 \leq  \frac{1}{\gamma}100\|v\|^2_{L^\infty}\|v_x\|^2_{L^2}+ \frac{1}{\beta}\|f_{xxxx}(t)\|^2.
$$
On the other hand, \eqref{eq:aux_11} implies
$$
\frac{1}{2}\frac{d}{dt}\|u_{xxx}\|_{L^2}^2 + \|w\|_{L^2}^2 = \frac{7}{2}(u_{xxx}, u_x u_{xxx}) + (f_{xx},w).
$$
After multiplication by two and some simple computations
$$
\frac{d}{dt}\|u_{xxx}\|_{L^2}^2 + (2-\alpha)\|w\|_{L^2}^2 \leq  7\|u_x\|_{L^\infty}\|u_{xxx}\|_{L^2}^2+ \frac{1}{\alpha}\|f_{xx}(t)\|^2.
$$
Now suppose that evolution is inside the trapping set
\begin{align*}
& W_{H^3}(R_1,R_2,R_3,R_4) \\
&\ \ = \{ v\in X\,:\  \|v\|_{L^2}\leq R_1, \|v_x\|_{L^2}\leq R_2, \|v_{xx}\|_{L^2}   \leq R_3, \|v_{xxx}\|_{L^2}\leq R_4\}.
\end{align*}
Then there hold two differential inequalities
$$
\frac{d}{dt}\|w\|_{L^2}^2 + (2\pi^2-\beta-\gamma - 9 \sqrt{2}R_2 R_3)\|w\|_{L^2}^2 \leq  \frac{1}{\gamma}100R_3R_4^3+ \frac{1}{\beta}\|f_{xxxx}\|^2_{L^\infty(L^2)}.
$$
$$
\frac{d}{dt}\|u_{xxx}\|_{L^2}^2 + (2-\alpha)\|w\|_{L^2}^2 \leq  7\sqrt{2}\sqrt{R_2 R_3}R_4^2+ \frac{1}{\alpha}\|f_{xx}\|_{L^\infty(L^2)}^2.
$$
These two inequalities, by Corollary \ref{cor:29} and the earlier obtained bound on $\|u_{xxx}\|_{L^2}$ allow us to find the trapping set for
$\|w\|_{L^2}$. Let
$$Y = \{ u\in H^4\cap H^1_0\,:\ u_{xx}\in H^1_0\}.$$ Define the set
\begin{align*}
& W_{H^4}(R_1,R_2,R_3,R_4,R_5,S) \\
&\ \ = \{ v\in Y\,:\  \|v\|_{L^2}\leq R_1, \|v_x\|_{L^2}\leq R_2, \|v_{xx}\|^2_{L^2}   \leq R_3,\\
& \qquad \qquad \|v_{xxx}\|_{L^2}\leq R_4,  \|v_{xxxx}\|^2_{L^2}   + S\|v_{xxx}\|^2_{L^2} \leq R_5^2 \}.
\end{align*}
Now, we can use Corollary \ref{cor:29} taking
\begin{align}
&\label{abcd_4}A = R_3^2,\ \  B=7\sqrt{2}\sqrt{R_2 R_3}R_4^2+ \frac{\|f_{xx}\|_{L^\infty(L^2)}^2}{\alpha},\ \  C = 2-\alpha,\\
&\nonumber D= \frac{100R_3R_4^3}{\gamma}+ \frac{\|f_{xxxx}\|^2_{L^\infty(L^2)}}{\beta},\ \  E=\beta+\gamma + 9 \sqrt{2}R_2 R_3 - 2\pi^2,
\end{align}
which leads us to the following result

\begin{lemma}\label{lem:u4x-bnd}
	Assume that $f\in L^\infty(Y)$. There exists the $Y$-trapping set which is nonempty and bounded in $H^4$. In fact if only $R_1$--$R_4$ are taken as in Section  \ref{subsec:u3x-bnd}, and
	$$
	R_5^2 \geq F(\alpha,\beta,\gamma),\ \  S = G(\alpha,\beta,\gamma)
	$$
	for some $\alpha, \beta,\gamma>0$, where
	$$
	F(\alpha,\beta,\gamma) = \begin{cases}
	-\frac{D}{E} \ \textrm{when}\ CD+BE\leq 0 \ \textrm{or}\ CD+BE > 0 \ \textrm{and}\ \sqrt{\frac{CD+BE}{A}} \leq -E,\\
	\frac{1}{C}\left(EA+B+2 \sqrt{A(CD+BE)}\right)  \ \textrm{when}\ CD+BE > 0 \ \textrm{and}\ \sqrt{\frac{CD+BE}{A}} > -E,
	\end{cases}
	$$
	$$
	G(\alpha,\beta,\gamma) = \begin{cases}
	0 \ \textrm{when}\ CD+BE\leq 0 \ \textrm{or}\ CD+BE > 0 \ \textrm{and}\ \sqrt{\frac{CD+BE}{A}} \leq -E,\\
	\frac{1}{C}\left(E+\sqrt{\frac{CD+BE}{A}}\right) \ \textrm{when}\ CD+BE > 0 \ \textrm{and}\ \sqrt{\frac{CD+BE}{A}} > -E,
	\end{cases}
	$$
	and
	$A, B(\alpha), C(\alpha), D(\beta,\gamma), E(\beta,\gamma)$ are given by \eqref{abcd_4},
	then the set $W_{H^3}(R_1,R_2,R_3,R_4,R_5,S)$, is $Y$-trapping.
\end{lemma}

The following results give alternative bounds for the radius of positively invariant set for $\|u_{xxxx}\|_{L^2}$.
\begin{lemma}\label{lem:u4x-bnd_2}
	If $R_1$--$R_4$ are  as in Section \ref{subsec:u3x-bnd} then there holds the bound
	\begin{equation}\label{eq:uxxxx}
	\frac{1}{2}\frac{d}{dt}\|w\|_{L^2}^2 \leq \|w_x\|_{L^2}\left(\|f_{xxx}(t)\|_{L^2} - \|w_x\|_{L^2} + 11R_1^{1/2}R_2^{1/2}R_4^{1/2}\|w_x\|_{L^2}^{1/2} + 10\sqrt{2}R_2^{1/2}R_3^{1/2}R_4\right).
	\end{equation}
	Hence if only $A$ is the positive root of
	$$
	x - \|f_{xxx}\|_{L^\infty(L^2)} - 11R_1^{1/2}R_2^{1/2}R_4^{1/2}x^{1/2} - 10\sqrt{2}R_2^{1/2}R_3^{1/2}R_4 = 0,
	$$
	and
	$$
	R_5 \geq \min\left(A/\pi,(A R_4)^{1/2}\right),
	$$
	then the set
	\begin{align*}
	& W_{H^4}(R_1,R_2,R_3,R_4,R_5)  = \{ v\in X\,:\  \|v\|_{L^2}\leq R_1, \|v_x\|_{L^2}\leq R_2, \|v_{xx}\|^2_{L^2}   \leq R_3,\\
	& \qquad \qquad \|v_{xxx}\|_{L^2}\leq R_4,  \|v_{xxxx}\|_{L^2}\leq R_5 \}
	\end{align*}
	is $Y$-trapping.
	\end{lemma}
 \begin{proof}
We rewrite \eqref{eq:u4} as
$$\frac{1}{2}\frac{d}{dt}\|w\|_{L^2}^2 + \|w_x\|_{L^2}^2 -9 (u w_x, w) + 10 (v_xv,w) = -(f_{xxx},w_x).
$$
As
$$
(vv_x,w) = - (u_xw,w)-(u_xv_x,w_x) = 2(uw,w_x) - (u_xv_x,w_x),
$$
we can rewrite the above equality as
$$\frac{1}{2}\frac{d}{dt}\|w\|_{L^2}^2 + \|w_x\|_{L^2}^2 + 11(uw,w_x) - 10(u_xv_x,w_x) = -(f_{xxx},w_x).
$$
We deduce
$$\frac{1}{2}\frac{d}{dt}\|w\|_{L^2}^2 \leq \|f_{xxx}\|_{L^2}\|w_x\|_{L^2} - \|w_x\|_{L^2}^2 + 11\|u\|_{L^\infty}\|w\|_{L^2}\|w_x\|_{L^2} + 10\|u_x\|_{L^\infty}\|v_x\|_{L^2}\|w_x\|_{L^2}.
$$
By interpolation inequalities it follows that
\begin{align*}
& \frac{1}{2}\frac{d}{dt}\|w\|_{L^2}^2 \leq \|w_x\|_{L^2} \cdot \\
& \ \  \cdot \left(\|f_{xxx}\|_{L^2} - \|w_x\|_{L^2} + 11\|u\|_{L^2}^{1/2}\|u_x\|_{L^2}^{1/2}\|w_x\|_{L^2}^{1/2}\|v_x\|^{1/2} + 10\sqrt{2}\|u_x\|^{1/2}_{L^2}\|u_{xx}\|^{1/2}_{L^2}\|v_x\|_{L^2}\right),
\end{align*}
 and the proof is complete.
 \end{proof}

\subsubsection{Local bounds}
By the Young inequality applied to \eqref{eq:uxxxx} we deduce the following Lemma which is useful to get the local estimates for $\|u_{xxxx}\|_{L^2}$
\begin{lemma}\label{lem:u4}
	Suppose that on interval $[t_1,t_2]$ there hold bounds $\|u(t)\|_{L^2}\leq R_1$, $\|u_x(t)\|_{L^2}\leq R_2$, $\|u_{xx}(t)\|_{L^2}\leq R_3$, $\|u_{xxx}(t)\|_{L^2}\leq R_4$. Then for every $\alpha,\beta,\gamma>0$ such that $\alpha+\beta+\gamma \leq 2$ there holds the estimate
$$
\frac{d}{dt}\|w\|_{L^2}^2 + (2-\alpha-\beta-\gamma)\pi^2\|w\|_{L^2}^2 \leq  \frac{\|f_{xxx}(t)\|_{L^2}^2}{\alpha} + \frac{200R_2R_3R_4^2}{\beta}   + \frac{3^3 11^4 R_1^2R_2^2R_4^2}{2^4\gamma^3}.
$$
\end{lemma}
Now, similar as in Lemma \ref{lem_3.18} there hold the following bounds with $\delta, \epsilon > 0$, and $\delta+\epsilon < 2$
$$
\frac{d}{dt}\|v_x\|_{L^2}^2 + (2-\epsilon-\delta)\|w\|_{L^2}^2 \leq\frac{3 \cdot 7^{8/3} 5^{5/3}R_2^{8/3}R_3^{2}}{2^{25/3}\gamma^{5/3}}+\frac{\|f_{xx}\|_{L^\infty(L^2)}^2}{\epsilon}.
$$
$$
\frac{d}{dt}\|v_x\|_{L^2}^2 + (2-\epsilon-\delta)\|w\|_{L^2}^2 \leq \frac{7^43^3R_1^2R_2^2R_3^2}{2^4\gamma^3}+\frac{\|f_{xx}\|_{L^\infty(L^2)}^2}{\epsilon}.
$$
Hence, we can derive two more algorithms to find the radius of the positively invariant set for $\|u_{xxxx}\|_{L^2} = \|w\|_{L^2}$ using Corollary \ref{cor:29} with the following two sets of parameters
\begin{align}
& A = R_4^2,\ \  B=\frac{3 \cdot 7^{8/3} 5^{5/3}R_2^{8/3}R_3^{2}}{2^{25/3}\delta^{5/3}}+\frac{\|f_{xx}\|_{L^\infty(L^2)}^2}{\epsilon},\ \  C = 2-\epsilon-\delta,\label{abcd_41}\\
& D= \frac{\|f_{xxx}\|_{L^\infty(L^2)}^2}{\alpha} + \frac{200R_2R_3R_4^2}{\beta}   + \frac{3^3 11^4 R_1^2R_2^2R_4^2}{2^4\gamma^3},\ \  E=(\alpha+\beta+\gamma-2)\pi^2,\nonumber\\
&A = R_4^2,\ \  B=\frac{7^43^3R_1^2R_2^2R_3^2}{2^4\delta^3}+\frac{\|f_{xx}\|_{L^\infty(L^2)}^2}{\epsilon},\ \  C = 2-\epsilon-\delta,\label{abcd_42}\\
&D= \frac{\|f_{xxx}\|_{L^\infty(L^2)}^2}{\alpha} + \frac{200R_2R_3R_4^2}{\beta}   + \frac{3^3 11^4 R_1^2R_2^2R_4^2}{2^4\gamma^3},\ \  E=(\alpha+\beta+\gamma-2)\pi^2.\nonumber
\end{align}

\begin{lemma}\label{lem:xxxx}
	Assume that $f\in L^\infty(H^3\cap H^1_0)$ satisfies $f_{xx}\in L^\infty(H^1_0)$. If $R_1$--$R_4$ are taken as in Section  \ref{subsec:u3x-bnd}, and
	$$
	R_5^2 \geq F(\alpha,\beta,\gamma,\delta),\ \  S = G(\alpha,\beta,\gamma,\delta)
	$$
	for some $\alpha, \beta,\gamma,\delta,\epsilon>0$ such that $\alpha+\beta+\gamma \leq 2$ and $\epsilon+\delta < 2$, where
	$$
	F(\alpha,\beta,\gamma,\delta,\epsilon) = \begin{cases}
	-\frac{D}{E} \ \textrm{when}\ CD+BE\leq 0 \ \textrm{or}\ CD+BE > 0 \ \textrm{and}\ \sqrt{\frac{CD+BE}{A}} \leq -E,\\
	\frac{1}{C}\left(EA+B+2 \sqrt{A(CD+BE)}\right)  \ \textrm{when}\ CD+BE > 0 \ \textrm{and}\ \sqrt{\frac{CD+BE}{A}} > -E,
	\end{cases}
	$$
	$$
	G(\alpha,\beta,\gamma,\delta , \epsilon) = \begin{cases}
	0 \ \textrm{when}\ CD+BE\leq 0 \ \textrm{or}\ CD+BE > 0 \ \textrm{and}\ \sqrt{\frac{CD+BE}{A}} \leq -E,\\
	\frac{1}{C}\left(E+\sqrt{\frac{CD+BE}{A}}\right) \ \textrm{when}\ CD+BE > 0 \ \textrm{and}\ \sqrt{\frac{CD+BE}{A}} > -E,
	\end{cases}
	$$
	and
	$A, B(\delta, \epsilon), C(\delta, \epsilon), D(\alpha, \beta,\gamma), E(\alpha,\beta,\gamma)$ are given by either of two possibilities \eqref{abcd_41}--\eqref{abcd_42},
	then the set $W_{H^3}(R_1,R_2,R_3,R_4,R_5,S)$, is $Y$-trapping.
\end{lemma}
Using \eqref{eq:uxxxx} and arguing the same as in Lemma \ref{lem:u4}, by the interpolation inequality $\|w\|_{L^2}^2\leq \|u_{xxx}\|_{L^2}\|w_x\|_{L^2} \leq R_4\|w_x\|_{L^2}$ we obtain the local estimate alternative to the one of Lemma \ref{lem:u4}.

\begin{equation}\label{eq:uxxxx_tan}
\frac{d}{dt}\|w\|_{L^2}^2 + \frac{2-\alpha-\beta-\gamma}{R_4^2}\|w\|_{L^2}^2 \leq  \frac{\|f_{xxx}(t)\|_{L^2}^2}{\alpha} + \frac{200R_2R_3R_4^2}{\beta}   + \frac{3^3 11^4 R_1^2R_2^2R_4^2}{2^4\gamma^3}.
\end{equation}

This gives us the following result

\begin{lemma}\label{lemma:young_x_5}
	Let $f\in L^\infty(t_1,t_2;H^3\cap H^1_0)$ satisfy $f_{xx}\in L^\infty(t_1,t_2;H^1_0)$ and let $\alpha,\beta,\gamma > 0$ be such that $\alpha+\beta+\gamma \leq 2$. Assume that the solution of the Burgers equation $u:[t_1,t_2]\to Y$ satisfies the estimates
	$$
	\|u(s)\|_{L^2} \leq R_1,\, \|u_x(s)\|_{L^2} \leq R_2,\, \|u_{xx}(s)\|_{L^2}\leq R_3,\, \|u_{xxx}\|_{L^2}\leq R_4 \quad \textrm{for}\quad s\in [t_1,t_2].
	$$
	Then
	\begin{align}
	&  \|u_{xxxx}(t)\|_{L^2}^2 \leq  \frac{D \tanh(\sqrt{CD}  (t-t_1)) + \sqrt{CD} \|u_{xxxx}(t_1)\|_{L^2}^2 }{C\tanh(\sqrt{CD}  (t-t_1)) \|u_{xxxx}(t_1)\|_{L^2}^2 + \sqrt{CD} }
	\end{align}
	with
	$$
	C = \frac{2-\alpha-\beta-\gamma}{R_4^2} \quad \textrm{and} \quad D = \frac{\|f_{xxx}\|_{L^\infty(t_0,t;L^2)}^2}{\alpha} + \frac{200R_2R_3R_4^2}{\beta}   + \frac{3^3 11^4 R_1^2R_2^2R_4^2}{2^4\gamma^3}.
	$$
\end{lemma}


\end{document}